\documentclass[11pt]{article}
\usepackage[T1]{fontenc}
\usepackage{lmodern,amsmath,amssymb,amsthm,mathtools}
\usepackage[margin=1.08in]{geometry}
\usepackage{enumitem,microtype}
\usepackage[hidelinks]{hyperref}
\numberwithin{equation}{section}
\newtheorem{theorem}{Theorem}[section]
\newtheorem{proposition}[theorem]{Proposition}
\newtheorem{lemma}[theorem]{Lemma}
\theoremstyle{remark}
\newtheorem*{remark}{Remark}
\newcommand{\R}{\mathbb R}
\newcommand{\Sph}{\mathbb S}
\newcommand{\eps}{\varepsilon}
\newcommand{\dd}{\,\mathrm d}
\newcommand{\tr}{\operatorname{tr}}
\newcommand{\diver}{\operatorname{div}}
\newcommand{\supp}{\operatorname{supp}}
\newcommand{\dist}{\operatorname{dist}}
\newcommand{\Span}{\operatorname{span}}

\newcommand{\cQ}{\mathcal Q}
\newcommand{\norm}[1]{\lVert#1\rVert}
\setlist[enumerate,1]{label=\textup{(\roman*)},leftmargin=2.3em}
\allowdisplaybreaks[2]
\title{Isolated singularities of solutions to the Yamabe equation
with boundary in higher dimensions I}
\author{Yuxuan Liao\thanks{School of Mathematical Sciences,
Beijing Normal University.
\href{mailto:yxliao@mail.bnu.edu.cn}{\texttt{yxliao@mail.bnu.edu.cn}}}}
\date{November 14, 2026}
\begin{document}
\maketitle
\begin{abstract}
We study isolated boundary singularities of positive solutions of the
scalar-flat boundary Yamabe equation. We prove the critical upper bound
for general smooth metrics in dimensions $3\le n\le17$ and for Euclidean
domains with smooth boundary in dimensions $3\le n\le22$.
For general metrics, this range extends to $n=18$ when the conformal
metric is complete near the singularity.
In every dimension, the critical upper bound implies removability or a
matching lower bound and convergence to a unique flat Fowler solution
with a positive power rate. Sharpness of these dimension ranges will be
established in a subsequent paper.

\medskip
\noindent\emph{Keywords.} Boundary Yamabe equation, isolated singularities,
Fowler solutions.

\medskip
\noindent\emph{MSC (2020).} Primary 35J61. Secondary 35B40, 35J25, 53C21.
\end{abstract}

\section{Introduction}
Let $g$ be a smooth Riemannian metric on the closed half-ball
$\overline{\mathcal B_1^+}\subset\R^n$, $n\ge3$. We study positive
solutions of the scalar-flat boundary Yamabe equation
\begin{equation}\label{eq:main}
 \begin{cases}
 -L_gu=0&\text{in }\mathcal B_1^+,\\
 B_gu=(n-2)u^{\frac{n}{n-2}}&\text{on }B_1\setminus\{0\}.
 \end{cases}
\end{equation}
Here $\mathcal B_r^+=\{x_n>0,\ |x|<r\}$ and
$B_r=\{x_n=0,\ |x'|<r\}$. The operators
\[
 L_g=\Delta_g-\frac{n-2}{4(n-1)}R_g,\qquad
 B_g=\partial_{\nu_g}+\frac{n-2}{2}h_g
\]
are the conformal Laplacian and conformal boundary operator.
Here $\Delta_g$ is the Laplace--Beltrami operator and $R_g$ is the scalar
curvature. The normal $\nu_g$ is outward, and $h_g$ is the average of the
principal curvatures, positive on a Euclidean sphere. Solutions are assumed smooth
up to the boundary away from zero.

The scalar-flat boundary Yamabe problem seeks a metric of zero scalar
curvature and constant boundary mean curvature in a prescribed conformal
class on a compact manifold with boundary \cite{Escobar,Almaraz}.
Equation~\eqref{eq:main} is the local form of this problem in the case
of positive boundary mean curvature. Allowing the conformal factor to
have isolated boundary singularities leads to the question of its growth
and asymptotic behavior near these points.

For the interior Yamabe equation with a flat background metric,
Caffarelli--Gidas--Spruck \cite{CGS} proved that every nonremovable
isolated singularity is asymptotic to a radial Fowler solution.
Korevaar et~al.~\cite{KMPS} obtained a positive power
rate and refined asymptotics. For general smooth metrics, Marques
\cite{Marques} established Fowler asymptotics in dimensions three
through five, and Xiong--Zhang \cite{XZ} treated dimension six.
The survey by Du et~al.~\cite{DJXY} reviews
these results, related equations, and the questions concerning higher
dimensions and global singular solutions.

Recently, Han et~al.~\cite{HJWXZ1} extended this
classification through dimension twenty-nine by proving the critical
upper bound for general smooth metrics. They also proved a matching
lower bound under the upper-bound hypothesis in every dimension.
Their companion paper \cite{HJWXZ2}
establishes sharpness by constructing complete singular solutions
violating the critical upper bound in every dimension at least thirty.

For the boundary equation on a flat Euclidean half-ball,
Caffarelli et~al.~\cite{CJSX} proved in every dimension that
an isolated singularity is either removable or has two-sided critical
growth, and established asymptotic symmetry and the symmetry of global
singular solutions. For smooth curved backgrounds, Liao--Ma \cite{LM}
proved the removable-or-two-sided-growth alternative and asymptotic
cylindrical symmetry in dimensions three and four. Jin--Liao--Xiong
\cite{JLXBoundary} classified the global singular solutions of the flat
boundary equation and their Jacobi fields, and proved convergence of
local singular solutions to a unique Fowler model with a positive power
rate.

We combine and adapt the methods of these works to study isolated
boundary singularities in the curved setting. We first recall the
notions of removability and completeness. We call zero a removable
singularity if $u$ extends to a positive smooth solution across zero.
The conformal metric $u^{4/(n-2)}g$ is complete at zero if every
locally rectifiable curve tending to zero has infinite length in this
metric.

\begin{theorem}\label{thm:classification}
Let $u$ be a positive smooth solution of \eqref{eq:main} on
$\overline{\mathcal B_1^+}\setminus\{0\}$. Assume that one of the
following holds.
\begin{enumerate}
\item $3\le n\le17$ and $g$ is an arbitrary smooth metric.
\item $3\le n\le18$, $g$ is an arbitrary smooth metric, and the
conformal metric $u^{4/(n-2)}g$ is complete at zero.
\item $3\le n\le22$ and $g$ is the pullback of the Euclidean metric
by a smooth boundary chart.
\end{enumerate}
Then there are $C,r_0>0$ such that
\begin{equation}\label{eq:upper-main}
 u(x)\le C |x|^{-{\frac{n-2}{2}}}\quad\text{for }0<|x|<r_0.
\end{equation}
\end{theorem}

\begin{remark}
At a nonumbilic boundary point, the upper bound \eqref{eq:upper-main}
holds in every dimension $n\ge3$
(see Section~\ref{subsec:geometric-upper}).
\end{remark}

In a subsequent paper, we will construct counterexamples showing that
each dimension range in Theorem~\ref{thm:classification} is sharp.

We write the flat Fowler models as
\begin{equation}\label{eq:model}
 U_{H,\tau}(r\omega)=r^{-{\frac{n-2}{2}}}v_H(-\log r+\tau,\omega),
 \quad\text{for }H_c\le H<0.
\end{equation}
Here $H$ is the cylindrical Hamiltonian parameter. For $H_c<H<0$,
$v_H$ is periodic and nonconstant in the logarithmic radius, with least
period $T(H)$, and phases are identified modulo $T(H)$.
At $H=H_c<0$, $v_{H_c}=v_c$ is independent of that radius and the
model is independent of $\tau$. The angular
dependence of the harmonic extension is included in $v_H$. We fix a
maximum-centered representative for every nonconstant member.
The cylindrical equation and Hamiltonian are specified in
\eqref{eq:flat-cylinder} and \eqref{eq:hamiltonian}.

\begin{theorem}\label{thm:conditional}
Let $n\ge3$ and let $g$ be smooth. Suppose a positive smooth solution
$u$ of \eqref{eq:main} satisfies $u(x)=O(|x|^{-(n-2)/2})$ as $x\to0$.
Then either zero is removable, or, near zero,
\begin{equation}\label{eq:two-sided}
 c |x|^{-{\frac{n-2}{2}}}\le u(x)\le C |x|^{-{\frac{n-2}{2}}},
\end{equation}
and
\begin{equation}\label{eq:main-relative}
 u(x)=U_{H,\tau}(x)\bigl(1+O(|x|^\beta)\bigr)
 \qquad\text{as }x\to0,
\end{equation}
for some $c,C,\beta>0$ and a unique flat Fowler solution $U_{H,\tau}$,
in fixed boundary Fermi coordinates centered at zero.
\end{theorem}

The convergence in \eqref{eq:main-relative} holds through two derivatives,
\begin{equation}\label{eq:main-asymptotics}
 \sum_{j=0}^2|x|^{\frac{n-2}{2}+j}|D^j(u-U_{H,\tau})(x)|
 \le C|x|^\beta.
\end{equation}

The upper- and lower-bound arguments adapt the methods of
\cite{HJWXZ1}, using the boundary Harnack estimates of \cite{LM}.
The Fowler asymptotics use the global classification and linearized
theory of \cite{JLXBoundary}.

Section~\ref{sec:upper} proves the geometric upper bounds, including
dimension eighteen when the conformal metric is complete at zero.
Sections~\ref{sec:lower} and \ref{sec:asymptotics} prove the lower bound
and Fowler asymptotics under \eqref{eq:upper-main}.
The appendices supply the linear estimates and geometric quadratic-form calculations.

\section{The critical upper bounds}\label{sec:upper}
We prove the upper bounds in Theorem~\ref{thm:classification}.
Proposition~\ref{prop:upper-criterion} reduces them to coercivity of
the geometric dilation form. Its adjoint argument adapts
\cite[Section 3]{HJWXZ1} to a boundary trace $L^1$ estimate.
The cases $n=3,4$ follow from \cite[Theorem 1.1]{LM}, so we assume $n\ge5$ in this section.

\subsection{A geometric criterion for the upper bound}\label{subsec:adjoint-argument}
Using the method from \cite{LM}, we first reduce the upper bound to the flat boundary.
Take boundary Fermi coordinates centered at zero, with $g(0)=I$,
$g_{nn}=1$ and $g_{in}=0$ for $i<n$. Then
$\nu_g=-\partial_{x_n}$ on $B_1$ and $d_g(x,0)\asymp|x|$ near zero.
Suppose the boundary estimate
$M=\sup_{B_\rho\setminus\{0\}}|x|^{(n-2)/2}u(x)<\infty$ holds. Set
$u_r(y)=r^{(n-2)/2}u(ry)$ and $g_r(y)=g(ry)$ for $0<r<\rho/2$,
then on $B_2\setminus\overline B_{1/2}$ the boundary equation is
\[
 \partial_{\nu_{g_r}}u_r=
 \left((n-2)u_r^{2/(n-2)}-\frac{n-2}{2}r h_g(ry)\right)u_r.
\]
The coefficient is bounded by $C(1+M^{2/(n-2)})$.
Here and below, positive constants $c,C$ may change from line to line.
Since $R_{g_r}(y)=r^2R_g(ry)$, the local Harnack argument of
\cite[proofs of Propositions 2.4 and 2.6]{LM} therefore establishes
\[
 \sup_{\mathcal B_{3/2}^+\setminus\mathcal B_{2/3}^+}u_r
 \le C\inf_{\mathcal B_{3/2}^+\setminus\mathcal B_{2/3}^+}u_r
 \le CM.
\]
Scaling back proves \eqref{eq:upper-main}. Thus it suffices to prove
the boundary upper bound.

Suppose that this boundary bound fails. The point selection in
\cite[Section 2.1.1]{LM} produces boundary points $x_k\to0$ such that
$\eps_k=u(x_k)^{-2/(n-2)}$ satisfies $\eps_k/|x_k|\to0$.
Let $\exp^{\partial}$ be the exponential map of the induced boundary
metric, and identify $y'$ with $T_{x_k}B_1$ by an orthonormal frame.
Tangential indices run from $1$ to $n-1$.
The Fermi map and the rescaled solution are
\[
 \Phi_k(y',y_n)=\exp^g_{\exp^{\partial}_{x_k}(y')}(-y_n\nu_g)\quad\text{ and }\quad
 v_k(y)=\eps_k^{(n-2)/2}u(\Phi_k(\eps_ky)).
\]
For every fixed $R$, the rescaled domains contain
$\overline{\mathcal B_{2R}^+}$ and their punctures lie outside it for
large $k$. By construction, $v_k(0)=1$ and
$v_k(y',0)\le2^{(n-2)/2}$ for $|y'|\le2R$.
Then \eqref{eq:buffered-robin-harnack} and the local estimates in
Appendix~\ref{app:local-robin} imply
\[
 C_R^{-1}\le v_k\le C_R,
 \qquad \|v_k\|_{C^j(\overline{\mathcal B_R^+})}\le C_{R,j}.
\]
After passage to a subsequence,
$v_k\to v_\infty>0$ in $C^j_{\rm loc}(\overline{\R^n_+})$.
Write $\Gamma=\partial\R^n_+$. The limiting function satisfies
\[
 \begin{cases}
 -\Delta v_\infty=0&\text{in }\R^n_+,\\
 -\partial_{x_n}v_\infty=(n-2)v_\infty^{\frac n{n-2}}
                                      &\text{on }\Gamma.
 \end{cases}
\]
By \cite[Theorem 1.2]{LiZhu}, it belongs to the bubble family below.

Set $U(x)=\bigl(|x'|^2+(1+x_n)^2\bigr)^{-(n-2)/2}$.
Its translates and dilates are
\begin{equation*}
 U_{\lambda,\xi}({x'},{x_n})=
 \left(\frac{\lambda}{|{x'}-\xi|^2+({x_n}+\lambda)^2}\right)^{\frac{n-2}{2}}.
\end{equation*}
For every $\lambda>0$ and $\xi\in\R^{n-1}$,
\[
 \begin{cases}
 -\Delta U_{\lambda,\xi}=0&\text{in }\R^n_+,\\
 \partial_\nu U_{\lambda,\xi}=(n-2)U_{\lambda,\xi}^{\frac n{n-2}}
                                      &\text{on }\Gamma.
 \end{cases}
\]

The limiting bubble has a nondegenerate boundary maximum. Smooth convergence
therefore supplies a nearby actual boundary maximum of $u$.
Replace $x_k$ by this maximum and choose $\eps_k$ from its
height. The previous scale and distance comparisons persist.
After recentering,
\begin{equation}\label{eq:rescaled-normalization}
 v_k(0)=1,\qquad\nabla_{x'}v_k(0)=0,\qquad
 v_k\longrightarrow U\quad\text{locally in }C^2.
\end{equation}
Fix a compact smooth boundary patch and an integer $N\ge1$.
\begin{lemma}\label{lem:fermi-normalization}
    Let \(n\ge5\) and let \(x_0\) lie in the patch. There is a positive
    function \(\varphi\) on a neighborhood of \(x_0\), satisfying
    \begin{equation}\label{eq:fermi-factor-normalization}
    \varphi(x_0)=1,\qquad\nabla_T\varphi(x_0)=0,
    \end{equation}
    such that \(\widetilde g=\varphi^{4/(n-2)}g\) admits boundary Fermi
    coordinates centered at \(x_0\), with geodesic normal tangential
    coordinates, in which
    \[
    \widetilde g_{nn}=1,\qquad \widetilde g_{in}=0\quad\text{for }i<n,\qquad
    \det\widetilde g=1+O(|x|^{N+1}),\qquad h_{\widetilde g}=O(|x|^N).
    \]
    The choices depend continuously on \(x_0\), with a common positive coordinate radius
    and uniform bounds for every fixed finite number of derivatives of \(\varphi\) and
    of the coordinate map.
\end{lemma}
\begin{proof}
We remove the homogeneous terms of $\log\det g$ successively,
preserving \eqref{eq:fermi-factor-normalization}. This refines the conformal
Fermi construction of \cite[Proposition 2.2]{Almaraz}.
Choose the initial boundary frames smoothly in $x_0$.

Write the conformal change as $e^{2f}g$, so that
$\varphi=e^{(n-2)f/2}$. Take $f(0)=0$ and its linear term to be
$f_1=h_g(x_0)x_n$. By the mean-curvature transformation law,
\[
 h_{e^{2f_1}g}(x_0)=h_g(x_0)-\partial_nf_1(x_0)=0.
\]
Thus \eqref{eq:fermi-factor-normalization} holds.
In Fermi coordinates with geodesic boundary coordinates, we have
$\log\det g=O(|x|^2)$.

Inductively, the current metric $g$ satisfies
\[
 \log\det g=\tau_d+O(|x|^{d+1})\qquad\text{for }2\le d\le N,
\]
where $\tau_d$ is homogeneous of degree $d$.
We seek a homogeneous polynomial $f_d$ of degree $d$.
Let $\Phi_d$ be the boundary Fermi map for $e^{2f_d}g$, using
the same origin and tangential frame.
The metric change is $O(|x|^d)$ and its first derivatives are
$O(|x|^{d-1})$. Integrating the boundary and normal geodesic equations, we obtain
\[
 \Phi_d(x)=x+X(x)+O(|x|^{d+2}),\qquad
 \Phi_d(0)=0,\qquad D\Phi_d(0)=I,
\]
where $X$ is homogeneous of degree $d+1$.
The conformal change and coordinate pullback are
\[
 g\longmapsto e^{2f_d}g\longmapsto\Phi_d^*(e^{2f_d}g).
\]
Since $g=\delta+O(|x|)$, their degree-$d$ contribution takes the form
\[
 \begin{aligned}
 \Phi_d^*(e^{2f_d}g)
     &=g+2f_d\delta+\mathcal L_X\delta+O(|x|^{d+1}),\\
 (\mathcal L_X\delta)_{ij}&=\partial_iX_j+\partial_jX_i.
 \end{aligned}
\]
The restored Fermi coordinates satisfy
\[
 \begin{aligned}
 (\Phi_d^*(e^{2f_d}g))_{nn}&=g_{nn}=1,\\
 (\Phi_d^*(e^{2f_d}g))_{in}&=g_{in}=0
     \quad\text{for }i<n.
 \end{aligned}
\]
Taking the degree-$d$ terms and using $\Phi_d(\Gamma)\subset\Gamma$, we find
\[
 \begin{cases}
 \partial_nX_n=-f_d&\text{in }\R^n_+,\\
 \partial_nX_i=-\partial_iX_n&\text{in }\R^n_+,\quad\text{for }i<n,\\
 X_n=0&\text{on }\Gamma.
 \end{cases}
\]
Integration gives
\begin{equation}\label{eq:gauge-vector}
 \begin{aligned}
 X_n&=-\int_0^{x_n}f_d(x',t)\dd t,\\
 X_T&=Y_d(x')+\int_0^{x_n}(x_n-t)\nabla_{x'}f_d(x',t)\dd t,
 \end{aligned}
\end{equation}
where $Y_d=X_T|_\Gamma$.
Gauss' lemma states that $g_{ij}(x',0)x'_j=x'_i$, hence
$(2f_d\delta+\mathcal L_{Y_d}\delta)x'=0$ on $\Gamma$.
Euler's identity implies
\[
 dY_d+\nabla_{x'}(x'\cdot Y_d)=-2f_d(x',0)x',\qquad
 x'\cdot Y_d=-\frac{|x'|^2f_d(x',0)}{d+1}.
\]
Hence
\[
 Y_d=-\frac{2x'f_d(x',0)}{d+1}
       +\frac{|x'|^2\nabla_{x'}f_d(x',0)}{d(d+1)}.
\]
The determinant transforms as
\[
 \begin{aligned}
 \log\det\bigl(\Phi_d^*(e^{2f_d}g)\bigr)
 &=2nf_d\circ\Phi_d+(\log\det g)\circ\Phi_d
       +2\log\det D\Phi_d\\
 &=\tau_d+2nf_d+2\diver X+O(|x|^{d+1}).
 \end{aligned}
\]
To remove the degree-$d$ term, we therefore impose
\[
 \tau_d+2nf_d+2\diver X=0.
\]
Since $\diver X=-f_d+\diver_{x'}X_T$, \eqref{eq:gauge-vector} yields
\[
 2(n-1)f_d+2\diver_{x'}Y_d
       +2\int_0^{x_n}(x_n-t)\Delta_{x'}f_d(x',t)\dd t=-\tau_d.
\]
On $\Gamma$ this is $L_d(f_d|_\Gamma)=-\tau_d|_\Gamma$, where
\[
 L_d=\frac2{d+1}
       \left((n-3)(d-1)+\frac{|x'|^2\Delta_{x'}}d\right).
\]
The decomposition into $|x'|^{2j}$ times harmonic polynomials of
degree $d-2j$ leads to the eigenvalues
\[
 \frac2{d+1}\left((n-3)(d-1)
      +\frac{2j(2d-2j+n-3)}d\right)>0
      \quad\text{for }0\le2j\le d.
\]
Thus $f_d|_\Gamma$ is uniquely determined.
Normal differentiation of the determinant equation on $\Gamma$ yields
$2(n-1)\partial_nf_d|_\Gamma=-\partial_n\tau_d|_\Gamma$, and
for $2\le j\le d$,
\[
 2(n-1)\partial_n^jf_d|_\Gamma
 =-\partial_n^j\tau_d|_\Gamma
       -2\Delta_{x'}\partial_n^{j-2}f_d|_\Gamma.
\]
These equations determine all remaining coefficients successively.
Their inverses depend only on $n,d$. The coefficients at each degree
are therefore polynomials in the metric jets with fixed denominators.
Since $f_d(0)=0$, $\nabla f_d(0)=0$ and $X(0)=DX(0)=0$ for $d\ge2$,
the normalization \eqref{eq:fermi-factor-normalization} is preserved.

After degrees $2,\ldots,N$, one has
$\log\det\widetilde g=O(|x|^{N+1})$, hence the required determinant
estimate. In the resulting Fermi coordinates,
\[
 h_{\widetilde g}(x',0)
 =-\frac1{2(n-1)}
       \partial_n\log\det\widetilde g(x',0)=O(|x'|^N).
\]
The finite construction and the boundary and normal exponential maps
depend smoothly on the center. Compactness of the patch ensures uniform
coordinate radii and every fixed finite derivative bound.
\end{proof}

At $x_k$, use Lemma~\ref{lem:fermi-normalization} for general metrics
and the Euclidean graph normalization of
Section~\ref{subsec:geometric-upper} for Euclidean domains.
Write $\widetilde g=\varphi^{4/(n-2)}g$ and
$\widetilde u=\varphi^{-1}u$.
By the normalization \eqref{eq:fermi-factor-normalization},
\[
 \widetilde u(x_k)=u(x_k),\qquad
 \nabla_T\widetilde u(x_k)=\nabla_Tu(x_k)=0.
\]
Let $\Phi_k$ now denote the normalized Fermi map. In each centered chart,
write $g,u$ for the pulled-back normalized metric and solution. Their
rescalings are
\[
 (g_k)_{ij}(y)=g_{ij}(\eps_k y),\qquad
 v_k(y)=\eps_k^{(n-2)/2}u(\eps_k y).
\]
The uniform bounds on these choices ensure
$g_k\to\delta$ locally and preserve $v_k\to U$ in $C^2_{\rm loc}$.

For a metric perturbation $\delta+\eta H$, we seek a correction
$U+\eta\psi_H$ satisfying the equation to first order in $\eta$.
For this flat linear problem, volume and boundary integrals are over
$\R^n_+$ and $\Gamma$, unless a domain is specified. Use the energy
\begin{equation}\label{eq:action}
 J_g(v)=\frac12\int(|\nabla v|_g^2+{\frac{n-2}{4(n-1)}}R_gv^2)\dd V_g
   +{\frac{n-2}{4}}\int h_gv^2\dd\sigma_g
   -{\frac{(n-2)^2}{2(n-1)}}\int(v_+)^{\frac{2(n-1)}{n-2}}\dd\sigma_g.
\end{equation}
Here $v_+=\max(v,0)$. Write $Q_g=J_g''(0)$ for the metric bilinear form
and $\mathcal R_g=J_g'$ for the weak residual. For compactly supported smooth tests,
\[
 \mathcal R_g(v)[\zeta]
 =-\int\zeta L_gv\,\dd V_g
  +\int\bigl(B_gv-(n-2)v_+^{n/(n-2)}\bigr)\zeta\,\dd\sigma_g.
\]
The correction belongs to the energy space
\[
 \mathcal D^{1,2}(\R^n_+)
 =\overline{C_c^\infty(\overline{\R^n_+})}^{\|\nabla\cdot\|_{L^2}},
 \qquad \|v\|_{\mathcal D^{1,2}}=\|\nabla v\|_{L^2(\R^n_+)}.
\]
For a compactly supported Fermi-preserving variation $H$, let
$J_1(H,f)=\partial_\eta J_{\delta+\eta H}(f)|_{\eta=0}$.
For polynomial $H$, $J_1$ denotes the cutoff limit whenever it exists.
The bubble Hessian is $Q_{\lambda,\xi}=J_\delta''(U_{\lambda,\xi})$, with
\[
 Q_{\lambda,\xi}(v,w)=\int\nabla v\cdot\nabla w
        -n\int_\Gamma U_{\lambda,\xi}^{\frac2{n-2}}vw.
\]
Write $Q=Q_{1,0}$ and $Q(v)=Q(v,v)$. Differentiating the equation for
$U+\eta\psi_H$ at $\eta=0$ gives
$Q(\psi_H,\zeta)=-D_fJ_1(H,U)[\zeta]$.

The bubble parameter derivatives lie in the kernel of $Q$,
\[
 \begin{aligned}
 Z_0&=\left.\partial_\lambda U_{\lambda,0}\right|_{\lambda=1}
       =-\frac{n-2}{2}U-x\cdot\nabla U,\\
 Z_i&=\left.\partial_{\xi_i}U_{1,\xi}\right|_{\xi=0}
       =-\partial_iU\quad\text{for }1\le i<n.
 \end{aligned}
\]
Write $\mathcal B_r=\{y\in\R^n:|y|<r\}$. The Cayley map converts
$Q$ to a ball Robin form,
\[
 T(x',x_n)=\left(\frac{2x'}{|x'|^2+(1+x_n)^2},
                 \frac{|x'|^2+x_n^2-1}{|x'|^2+(1+x_n)^2}\right).
\]
Since $T^*\delta=4U^{4/(n-2)}\delta$, substituting
$\phi=U(v\circ T)$ yields
\begin{equation}\label{eq:ball-robin-form}
 \begin{aligned}
 Q(\phi)&=2^{2-n}\left(\int_{\mathcal B_1}|\nabla v|^2
                         -\int_{\partial\mathcal B_1}v^2\right),\\
 \|\nabla\phi\|_2^2&=2^{2-n}\left(\int_{\mathcal B_1}|\nabla v|^2
                         +\frac{n-2}{2}\int_{\partial\mathcal B_1}v^2\right).
 \end{aligned}
\end{equation}
For a spherical harmonic $Y_\ell$ of degree $\ell$, its harmonic
extension is $r^\ell Y_\ell$ and its normal derivative on the sphere is
$\ell Y_\ell$. Thus the Steklov eigenvalues are $\ell=0,1,\ldots$,
and the Robin eigenvalues are $\ell-1$. Degree one corresponds to
$Z_0,\ldots,Z_{n-1}$. For the other degrees,
\[
 \frac{|\ell-1|}{\ell+(n-2)/2}\ge\frac2{n+2}\quad\text{for }\ell\ne1.
\]
The zero-trace part has positive Dirichlet form. Consequently
$\ker Q=\Span\{Z_0,\ldots,Z_{n-1}\}$, and $Q$ has a bounded inverse
on their Dirichlet-orthogonal complement, as in the flat case of
\cite[Proposition 2.7]{AlmarazBlowup}.

The source must therefore satisfy $D_fJ_1(H,U)[Z_\alpha]=0$.
For compactly supported $H$, this follows by differentiating
$J_1(H,U_{\lambda,\xi})=0$ in the bubble parameters.
To verify this metric identity, use
$U^{4/(n-2)}\delta=\tfrac14T^*\delta$ and the unit-ball functional
\begin{equation}\label{eq:ball-action}
 S(g)=\int_{\mathcal B_1}R_g\dd V_g
 +2(n-1)\int_{\partial\mathcal B_1}h_g\dd\sigma_g
 -2(n-2)|\partial\mathcal B_1|_g.
\end{equation}
For a metric variation $h$, its first variation is
\[
 \begin{aligned}
 S'_g[h]={}&-\int_{\mathcal B_1}
          \langle\mathrm{Ric}_g-\tfrac12R_gg,h\rangle_g\,\dd V_g\\
 &+\int_{\partial\mathcal B_1}
 \langle((n-1)h_g-(n-2))g_T-\mathrm{II}_g,h_T\rangle_g
                                                   \,\dd\sigma_g.
 \end{aligned}
\]
At $\delta$, $\mathrm{Ric}_\delta=0$, $h_\delta=1$ and
$\mathrm{II}_\delta=\gamma$, the round metric.
Both coefficient tensors vanish, so $S'_\delta[h]=0$.
By conformal covariance,
\[
 J_1(H,U)=\frac{n-2}{2^{n+1}(n-1)}
 S'_\delta\!\left[(T^{-1})^*\!\left(4U^{4/(n-2)}H\right)\right]=0.
\]
Translation and dilation imply $J_1(H,U_{\lambda,\xi})=0$ as well.

For a homogeneous Fermi-preserving $H$ of degree $d$, local
differentiation determines the volume and boundary sources
\begin{equation}\label{eq:metric-correction-source}
 \begin{aligned}
 F_H={}&\diver\bigl((\tfrac12\tr H\,I-H)\nabla U\bigr)
 -\frac{n-2}{4(n-1)}
       (\partial_i\partial_jH_{ij}-\Delta\tr H)U,\\
 F_H^{\partial}={}&\frac{n-2}{4(n-1)}(\partial_{x_n}\tr H)U.
 \end{aligned}
\end{equation}
The conormal and boundary-area variations cancel since $H_{in}=0$.
The full expansion is in Appendix~\ref{app:correction-estimates}.
The linear problem is therefore
\[
 Q(\psi_H,\zeta)=\int F_H\zeta+\int_\Gamma F_H^\partial\zeta.
\]

\begin{lemma}\label{lem:jet-inverse}
Let $H$ be a homogeneous Fermi-preserving tensor of degree
$d\ge1$, and suppose $n>2d+2$. The sources satisfy
\[
 \int_{\R^n_+}F_HZ_\alpha+\int_\Gamma F_H^\partial Z_\alpha=0
 \quad\text{for }0\le\alpha<n.
\]
The correction equation
\begin{equation}\label{eq:jet-correction}
 \begin{cases}
 -\Delta\psi_H=F_H&\text{in }\R^n_+,\\
 (\partial_\nu-nU^{\frac2{n-2}})\psi_H=F_H^{\partial}&\text{on }\Gamma
 \end{cases}
\end{equation}
has a solution in $\mathcal D^{1,2}(\R^n_+)$, unique modulo
$\Span\{Z_0,\ldots,Z_{n-1}\}$. The normalization conditions
\begin{equation}\label{eq:corrector-normalization}
 \psi_H(0)=0,\qquad\nabla_{x'}\psi_H(0)=0
\end{equation}
determine $\psi_H$ uniquely.
For every integer $j\ge0$,
\begin{equation}\label{eq:correction-decay}
 |D^j\psi_H(x)|\le C_j|H|(1+|x|)^{d+2-n-j}.
\end{equation}
For a constant symmetric tangential tensor $A$, the solution of
\eqref{eq:jet-correction} with normalization
\eqref{eq:corrector-normalization} is
\begin{equation}\label{eq:constant-correction}
 \psi_A(x)=\tfrac12A_{ij}x_i\partial_jU(x),
\end{equation}
and satisfies \eqref{eq:correction-decay} with $d=0$.

For $s>0$ and $\xi\in\R^{n-1}$, let $\psi_{H;s,\xi}$ denote the
solution obtained by replacing $U$ with $U_{s,\xi}$ in
\eqref{eq:metric-correction-source} and \eqref{eq:jet-correction},
with normalization
\begin{equation}\label{eq:parameter-corrector-normalization}
 \psi_{H;s,\xi}(\xi,0)=0,\qquad
 \nabla_{x'}\psi_{H;s,\xi}(\xi,0)=0.
\end{equation}
The decay estimates hold
uniformly after any fixed number of derivatives in $(s,\xi)$
on compact subsets of $(0,\infty)\times\R^{n-1}$.
\end{lemma}
\begin{proof}
By \eqref{eq:metric-correction-source},
\[
 |D^jF_H|\le C_j|H|(1+r)^{d-n-j},\qquad
 |D^jF_H^{\partial}|\le C_j|H|(1+r)^{d+1-n-j}.
\]
Since $d<(n-2)/2$, these bounds imply
\[
 \|F_H\|_{L^{2n/(n+2)}(\R^n_+)}
 +\|F_H^{\partial}\|_{L^{2(n-1)/n}(\Gamma)}\le C|H|.
\]
The Sobolev and trace inequalities therefore make the source functional bounded on $\mathcal D^{1,2}$.
Choose a smooth radial cutoff $\chi_R$, equal to one on
$\mathcal B_R^+$ and zero outside $\mathcal B_{2R}^+$.
Differentiating the compact-support identity after \eqref{eq:ball-action},
$J_1(\chi_RH,U_{s,\xi})=0$, at $(s,\xi)=(1,0)$ gives
\[
 \int F_{\chi_RH}Z_\alpha+\int_\Gamma F_{\chi_RH}^\partial Z_\alpha
 =-D_fJ_1(\chi_RH,U)[Z_\alpha]=0.
\]
The difference $F_{\chi_RH}-\chi_RF_H$ is supported in
$R<r<2R$, where the product rule implies
$|F_{\chi_RH}-\chi_RF_H|\le C|H|R^{d-n}$.
On $\Gamma$, $\partial_n\chi_R=0$, so
$F_{\chi_RH}^\partial=\chi_RF_H^\partial$.
Using $|Z_\alpha|\le C(1+r)^{2-n}$, we obtain
\[
 \begin{aligned}
 &\left|\int(F_{\chi_RH}-F_H)Z_\alpha\right|
 +\left|\int_\Gamma
       (F_{\chi_RH}^\partial-F_H^\partial)Z_\alpha\right|\\
 &\qquad\le C|H|\left(R^{d+2-n}
                   +\int_R^\infty r^{d+1-n}\dd r\right)
 \le C|H|R^{d+2-n}.
 \end{aligned}
\]
Since $d+2-n<0$, letting $R\to\infty$ proves the source orthogonality.

By \eqref{eq:ball-robin-form} and the spectral decomposition following it,
$Q$ is invertible on the Dirichlet-orthogonal complement of
$\Span\{Z_0,\ldots,Z_{n-1}\}$. The source orthogonality ensures the existence of a solution $\psi_H^\perp$ there, with
\[
 \int\nabla\psi_H^\perp\cdot\nabla Z_\alpha=0
 \quad\text{for }0\le\alpha<n.
\]
Every finite-energy solution lies in
$\psi_H^\perp+\Span\{Z_0,\ldots,Z_{n-1}\}$.
Local Robin regularity makes $\psi_H^\perp$ smooth up to $\Gamma$,
with the bound
\[
 |\psi_H^\perp(0)|+|\nabla_{x'}\psi_H^\perp(0)|\le C|H|.
\]
Since $Z_i(0)=\partial_iZ_0(0)=0$ for $i\ge1$ and
\[
 Z_0(0)=-\frac{n-2}{2},\qquad
 \partial_iZ_j(0)=(n-2)\delta_{ij}\quad\text{for }1\le i,j<n,
\]
imposing \eqref{eq:corrector-normalization} determines the Jacobi
coefficients, so
\[
 \psi_H=\psi_H^\perp+\frac{2\psi_H^\perp(0)}{n-2}Z_0
       -\frac1{n-2}\sum_{i=1}^{n-1}
                       \partial_i\psi_H^\perp(0)Z_i.
\]

For $d\ge1$, the Cayley map satisfies
$|T(x)-e_n|=2U(x)^{1/(n-2)}\asymp(1+r)^{-1}$,
where $e_n=(0,\ldots,0,1)$. Thus
\eqref{eq:compactified-correction-bound} implies
\[
 |\psi_H(x)|\le C|H|U(x)|T(x)-e_n|^{-d}
       \le C|H|(1+|x|)^{d+2-n}.
\]
Applying \eqref{eq:buffered-robin-estimate} on rescaled annuli
proves \eqref{eq:correction-decay} for every $j$.

For a constant $A$, $\Delta U=0$ and the tangential derivatives of
the bubble boundary equation give
$\partial_\nu\partial_jU=nU^{2/(n-2)}\partial_jU$ on $\Gamma$.
Hence \eqref{eq:constant-correction} satisfies
\[
 \begin{cases}
 -\Delta\psi_A=-A_{ij}\partial_{ij}U=F_A&\text{in }\R^n_+,\\
 (\partial_\nu-nU^{\frac2{n-2}})\psi_A=0=F_A^\partial
       &\text{on }\Gamma.
 \end{cases}
\]
Direct differentiation verifies $\psi_A(0)=0$,
$\nabla_{x'}\psi_A(0)=0$ and
$|D^j\psi_A(x)|\le C_j|A|(1+|x|)^{2-n-j}$.

Extend $H\mapsto\psi_H$ by linearity to finite sums of homogeneous
tensors. For $\theta=(s,\xi)$, the polynomial
$H_\theta(y)=H(sy'+\xi,sy_n)$ has degrees zero through $d$.
Substitute $x=(sy'+\xi,sy_n)$ in the correction equation.
Under this substitution, \eqref{eq:corrector-normalization} becomes
\eqref{eq:parameter-corrector-normalization}.
Uniqueness in Lemma~\ref{lem:jet-inverse} implies
\begin{equation}\label{eq:corrector-covariance}
 \psi_{H;s,\xi}(x)=s^{-\frac{n-2}{2}}
 \psi_{H_\theta}
       \left(\frac{x'-\xi}{s},\frac{x_n}{s}\right).
\end{equation}
The homogeneous coefficients of $H_\theta$ depend smoothly on
$\theta$. Their derivatives are bounded by $C_{\Theta,\beta}|H|$
on every compact $\Theta\subset(0,\infty)\times\R^{n-1}$.
Differentiating \eqref{eq:corrector-covariance} and applying the
estimates for degrees $0,\ldots,d$ yields
\[
 |D_\theta^\beta D_x^j\psi_{H;s,\xi}(x)|
 \le C_{\Theta,\beta,j}|H|(1+|x|)^{d+2-n-j}.
\]
Translation and dilation preserve the Dirichlet inner product and
map the bubble Jacobi spaces onto one another. Thus the covariance
formula and parameter estimates also hold for the
Dirichlet-orthogonal corrections.
\end{proof}

Take a smooth compactly supported
Fermi-preserving variation $h$ and a metric path $g_\eta$ with
$g_0=\delta$ and $\partial_\eta g_\eta|_{\eta=0}=h$, equal to
$\delta$ outside a fixed compact set.
Let $J_2(h)$ be the coefficient of $\eta^2$ in $J_{g_\eta}(U)$.
Since $J_1(k,U)=0$ for every compactly supported variation $k$,
this coefficient depends only on $h$. Let $\psi_h$ be the correction
determined by $Q(\psi_h,\zeta)=-D_fJ_1(h,U)[\zeta]$ and the center
conditions of Lemma~\ref{lem:jet-inverse}. Since
$D_fJ_1(h,U)[\psi_h]=-Q(\psi_h)$, the energy expansion is
\[
 J_{g_\eta}(U+\eta\psi_h)
 =J_\delta(U)+\eta^2\bigl(J_2(h)-\tfrac12Q(\psi_h)\bigr)+o(\eta^2).
\]
Define
\[
 \mathcal C(h):=J_2(h)-\frac12Q(\psi_h).
\]
We use symmetric bilinear polarization for $J_2$ and $\mathcal C$, so
\begin{equation}\label{eq:reduced-form}
 \mathcal C(h,k)=J_2(h,k)-\frac12Q(\psi_h,\psi_k).
\end{equation}
For polynomial variations of maximal degree $d<(n-2)/2$,
take these coefficients after cutting off at radius $R$ and then
let $R\to\infty$. By homogeneity and \eqref{eq:correction-decay}, the cutoff errors are $O(R^{2d+2-n})$, so these limits converge.

The same construction at $U_{\lambda,\xi}$ defines
$\mathcal C(h,k;\lambda,\xi)$. We write $\mathcal C(h;\lambda,\xi)=\mathcal C(h,h;\lambda,\xi)$
for the associated quadratic value.
The notation without parameters refers to $U$.
For homogeneous tensors $H_i,H_j$ of degrees $i,j$,
\eqref{eq:corrector-covariance} and $x=\lambda y$ imply
$\mathcal C(H_i,H_j;\lambda,0)=\lambda^{i+j}\mathcal C(H_i,H_j)$.
Define the dilation form by
\[
 \cQ(H_i,H_j)
 =-\left.\lambda\partial_\lambda
             \mathcal C(H_i,H_j;\lambda,0)\right|_{\lambda=1}
 =-(i+j)\mathcal C(H_i,H_j),
\]
and extend it bilinearly to sums of homogeneous terms.

Fix a small boundary disk $B_\rho$ and uniform centered coordinate
families, using Lemma~\ref{lem:fermi-normalization} for general metrics
and Section~\ref{subsec:geometric-upper} for Euclidean domains.
With the lower degrees fixed, the degree-$d$ conformal and
boundary-preserving coordinate variations have the form
\[
 \mathcal L_{Y_{d+1}}\delta+f_d\delta,
 \qquad (Y_{d+1})_n=0\quad\text{on }\Gamma,
\]
where $Y_{d+1}$ and $f_d$ are homogeneous of degrees $d+1$ and $d$.
After the chosen normalizations, let $H_d(z)$ denote the independent
degree-$d$ part of the logarithmic metric expansion modulo these
variations. All other degree-$d$ terms are polynomial expressions
in $H_1(z),\ldots,H_{d-1}(z)$, with total weight $d$ when $H_j$
is assigned weight $j$. Fix coefficient norms on these
finite-dimensional spaces.

Put $k=\lfloor(n-2)/4\rfloor$ and $L=\lfloor(n-3)/2\rfloor$.
The quadratic form retains degrees $d\le k$, and
Lemma~\ref{lem:jet-inverse} supplies the linear corrections for
$1\le d\le L$. Use the coefficient norm
\[
 |H(z)|_\eps^2=\sum_{d=1}^k\eps^{2d}|H_d(z)|^2.
\]

\begin{proposition}\label{prop:upper-criterion}
Let $n\ge5$. If there are $c,\eps_0>0$ such that, for every
$z\in B_\rho$ and $0<\eps<\eps_0$,
\begin{equation}\label{eq:geometric-upper-criterion}
 \cQ\left(\sum_{d=1}^k\eps^dH_d(z)\right)
 \ge c|H(z)|_\eps^2,
\end{equation}
then every positive solution of \eqref{eq:main} satisfies
\eqref{eq:upper-main} near zero.
\end{proposition}

\begin{proof}
Suppose the boundary upper bound fails and take the normalized
selected sequence constructed above. Write $g,v,\eps$ for
$g_k,v_k,\eps_k$, with $H_d=H_d(x_k)$ and
$|H|_\eps=|H(x_k)|_\eps$. Fix a small original radius $\delta>0$ and put
$R=\delta/\eps$. Then $v(0)=1$, $\nabla_{x'}v(0)=0$ and
$g_{ij}\to\delta_{ij}$, $v\to U$ locally in $C^2$.
Using the corrections from Lemma~\ref{lem:jet-inverse}, set
\[
 W=U+\sum_{d=1}^L\eps^d\psi_{H_d},\qquad w=v-W.
\]
Equations~\eqref{eq:rescaled-normalization} and
\eqref{eq:corrector-normalization} imply $w(0)=0$ and
$\nabla_{x'}w(0)=0$. By \eqref{eq:correction-decay}, for
$r<R$ and $\eps\le\delta$,
\begin{equation}\label{eq:corrected-bubble-positivity}
 \frac{|W-U|}{U}
 \le C\sum_{d=1}^L\eps^d|H_d|(1+r)^d\le C\delta.
\end{equation}
Choose $\delta$ so that $C\delta<1/2$. Thus $U/2\le W\le3U/2$.
Put $c_W=nW^{2/(n-2)}$ and define
\[
 N_W(w)=(n-2)\bigl((W+w)_+^{\frac n{n-2}}-W^{\frac n{n-2}}\bigr)
                    -c_Ww.
\]
By convexity, $N_W(w)\ge0$.

Define
\[
 W_\lambda(x)=U_{\lambda,0}(x)
       +\sum_{d=1}^L\eps^d\lambda^{d-\frac{n-2}{2}}
                                \psi_{H_d}(x/\lambda)
\]
and set $Z=-\partial_\lambda W_\lambda|_{\lambda=1}$.
Define the residuals by
\[
 \begin{cases}
 E=-L_gW&\text{in }\mathcal B_R^+,\\
 E^\partial=B_gW-(n-2)W^{\frac n{n-2}}&\text{on }B_R,
 \end{cases}
\]
and their pairing with a smooth test function $X$ by
\[
 \mathcal R[X]=\int_{\mathcal B_R^+}XE\,\dd V_g
                    +\int_{B_R}XE^\partial\,\dd\sigma_g.
\]
The equation for $v$ implies
\begin{equation}\label{eq:error-equations}
 \begin{cases}
 -L_gw=-E&\text{in }\mathcal B_R^+,\\
 (B_g-c_W)w=N_W(w)-E^\partial&\text{on }B_R.
 \end{cases}
\end{equation}
These equations hold weakly across the rescaled puncture by
Lemma~\ref{lem:no-source}.
The Green-weighted and energy-dual residual norms are defined in
\eqref{eq:defined-core-residual}. Lemma~\ref{lem:taylor-budget} provides the bounds
\[
 \norm{\mathcal R}_G\le C(|H|_\eps^2+\eps^{\frac{n-2}{2}}),\qquad
 \norm{\mathcal R}_*^2=o(|H|_\eps^2+\eps^{\frac{n-2}{2}}).
\]

Write $\Sigma_r=\{x_n>0,\ |x|=r\}$ for the spherical part of
$\partial\mathcal B_r^+$.
Choose $r_*$, independent of $\eps$, large enough that the exterior
Robin problem in Appendix~\ref{app:adjoint-estimates} is coercive.
The test function $Y$ constructed there satisfies
\[
 Y=Z\ \text{on }\mathcal B_{r_*/2}^+,\qquad
 Y=0\ \text{on }\Sigma_R,\qquad
 0\le Y_+\le C,\quad
 \supp Y_+\subset\overline{\mathcal B_{r_*}^+}.
\]
Set
\[
 I^\pm=\int_{B_R}Y_\pm N_W(w)\,\dd\sigma_g.
\]
Green's identity applied to \eqref{eq:error-equations}, with
$Y=0$ on $\Sigma_R$, gives
\begin{equation}\label{eq:adjoint}
 I^+-I^-=\int_{\Sigma_R}w\partial_{\nu_g}Y\,\dd\sigma_g
 +\mathcal R[Y]-\int_{\mathcal B_R^+}wL_gY\,\dd V_g
 +\int_{B_R}w(B_g-c_W)Y\,\dd\sigma_g.
\end{equation}
The function \(Y\) constructed in Appendix~\ref{app:adjoint-estimates} also satisfies
\begin{equation}\label{eq:patch-comparison}
 |\mathcal R[Y-Z]|\le C\norm{\mathcal R}_*^2+CR^{2-n}.
\end{equation}
By Lemma~\ref{lem:taylor-budget}, the assumed inequality
\eqref{eq:geometric-upper-criterion}, and
\eqref{eq:patch-comparison},
\begin{equation}\label{eq:upper-criterion}
 \begin{aligned}
 \mathcal R[Y]
 &=\mathcal R[Z]+O(\norm{\mathcal R}_*^2+R^{2-n})\\
 &=\cQ\left(\sum_{d=1}^k\eps^dH_d\right)
                +o(|H|_\eps^2+\eps^{\frac{n-2}{2}})
 \ge c|H|_\eps^2-o(|H|_\eps^2+\eps^{\frac{n-2}{2}}).
 \end{aligned}
\end{equation}
Here $R^{2-n}=O_\delta(\eps^{n-2})=o(\eps^{(n-2)/2})$.
The interior and boundary residuals of $Y$ are supported in $r<r_*$
and satisfy \eqref{eq:adjoint-local-errors}.
Since $w\to0$ uniformly on the fixed half-ball, these bounds imply
\[
 \begin{aligned}
 &\left|-\int_{\mathcal B_R^+}wL_gY\,\dd V_g
           +\int_{B_R}w(B_g-c_W)Y\,\dd\sigma_g\right|\\
 &\qquad\le C\|w\|_{C^0(\overline{\mathcal B_{r_*}^+})}
                  (\norm{\mathcal R}_G+R^{2-n})
   =o(|H|_\eps^2+\eps^{\frac{n-2}{2}}).
 \end{aligned}
\]

The outer spheres correspond to a fixed original half-annulus
separated from zero. Smooth positivity of $u$ and the uniform
normalizations ensure
\begin{equation}\label{eq:two-sided-outer-cap}
 c\eps^{\frac{n-2}{2}}\le v\le C\eps^{\frac{n-2}{2}}
                                      \quad\text{on }\Sigma_R.
\end{equation}
By \eqref{eq:corrected-bubble-positivity},
$W=O(R^{2-n})=O_\delta(\eps^{n-2})$ there.
Thus $w\ge c\eps^{(n-2)/2}$ on $\Sigma_R$ for small $\eps$.
The same construction ensures $\partial_{\nu_g}Y>0$ on $\Sigma_R$ and
\begin{equation}\label{eq:outer-flux}
 \int_{\Sigma_R}\partial_{\nu_g}Y\,\dd\sigma_g\ge c>0.
\end{equation}
Therefore
\[
 \int_{\Sigma_R}w\partial_{\nu_g}Y\,\dd\sigma_g
                         \ge c\eps^{\frac{n-2}{2}}.
\]
Combining this with \eqref{eq:adjoint} and
\eqref{eq:upper-criterion} yields
\[
 I^+-I^-\ge c'(|H|_\eps^2+\eps^{\frac{n-2}{2}})>0.
\]
The local error estimate, proved in Lemma~\ref{lem:core}, is
\begin{equation}\label{eq:core-estimate}
 \norm{w}_{L^1(B_{r_*})}
 \le C\left({\eps^{\frac{n-2}{2}}}+\norm{\mathcal R}_G+I^-\right).
\end{equation}
On $B_{r_*}$, $W$ is bounded away from zero and $w\to0$ uniformly,
so $0\le N_W(w)\le Cw^2$. Hence
\[
 I^+\le C\|w\|_{L^\infty(B_{r_*})}\|w\|_{L^1(B_{r_*})}
 \le o(1)\left(\eps^{\frac{n-2}{2}}+\norm{\mathcal R}_G+I^-\right).
\]
After absorbing the coefficient of $I^-$, we obtain
\[
 I^+-I^-\le o(1)\left(\eps^{\frac{n-2}{2}}+\norm{\mathcal R}_G\right)
                    -\tfrac12I^-.
\]
The residual estimate now implies
$I^+-I^-\le o(|H|_\eps^2+\eps^{(n-2)/2})$, contradicting the
positive lower bound above.
The boundary-to-interior reduction at the beginning
of this subsection proves \eqref{eq:upper-main}.
\end{proof}

\subsection{Verification for general metrics and Euclidean domains}\label{subsec:geometric-upper}

Let $h$ be a smooth compactly supported symmetric tensor with
$h_{in}=0$ and $\tr h=0$, and take $g_\eta=e^{\eta h}$.
Then $\det g_\eta=e^{\eta\tr h}=1$, so
$\dd V_{g_\eta}=\dd x$ and $\dd\sigma_{g_\eta}=\dd x'$.
By the Fermi conditions,
\[
 h_{g_\eta}=-\frac1{2(n-1)}
                 \partial_n\log\det g_\eta=0
                    \quad\text{on }\Gamma.
\]
Thus the boundary terms in \eqref{eq:action} are constant along
this path. The inverse metric satisfies
\[
 g_\eta^{-1}=I-\eta h+\tfrac12\eta^2h^2+O(\eta^3).
\]
Write
$R_{g_\eta}=\eta\partial_i\partial_jh_{ij}
                 +\eta^2R_2(h)+O(\eta^3)$.
The Christoffel formula yields
\[
 R_2(h)=-\tfrac14|\partial h|^2+\tfrac12|\diver h|^2
                         -\partial_i(h_{ik}\partial_jh_{jk}).
\]
This curvature expansion also appears in
\cite[proof of Proposition 3.2]{GKMW}.
Substitution in \eqref{eq:action} yields
\[
 J_2(h)=\frac14\int|h\nabla U|^2
                +\frac{n-2}{8(n-1)}\int R_2(h)U^2.
\]
On $\Gamma$, the outward normal is $-e_n$, and
\[
 h_{ik}(\partial_jh_{jk})\nu_iU^2
             =-h_{nk}(\partial_jh_{jk})U^2=0.
\]
Compact support and this boundary identity justify integration
by parts, yielding
\begin{equation}\label{eq:direct-form}
 \begin{aligned}
 J_2(h)={}&\frac14\int|h\nabla U|^2
       -\frac{n-2}{32(n-1)}\int|\partial h|^2U^2\\
 &+\frac{n-2}{16(n-1)}\int|\diver h|^2U^2
       +\frac{n-2}{8(n-1)}
                    \int h_{ik}(\partial_jh_{jk})\partial_i(U^2).
 \end{aligned}
\end{equation}
For homogeneous terms $H_d,H_e$ of degrees $d,e$, each polarized
spherical flux satisfies
\[
 \left|\int_{\Sigma_R}(H_d)_{ik}\partial_j(H_e)_{jk}\nu_iU^2\right|
                                  \le CR^{d+e+2-n}.
\]
The cutoff derivative terms on $R<r<2R$ have the same bound.
Hence the polarized formula extends to polynomial tensors with
$d+e<n-2$ by cutting off first and then letting $R\to\infty$.

For the correction energy, \eqref{eq:metric-correction-source}
reduces in $\R^n_+$ to
\[
 F_h=-h_{ij}\partial_{ij}U-(\partial_i h_{ij})\partial_jU
          -\frac{n-2}{4(n-1)}(\partial_i\partial_jh_{ij})U.
\]
On $\Gamma$, $F_h^\partial=0$ because $\tr h=0$.

By Lemma~\ref{lem:fermi-normalization}, $h_g(0)=0$ and the boundary coordinates are geodesic.
Put $A=\mathrm{II}^{\circ}(0)$. For $i,j<n$, the Fermi expansion is
\[
 g_{ij}(x)=\delta_{ij}-2x_nA_{ij}+O(|x|^2).
\]
Thus the first logarithmic metric coefficient is $H_1=-2x_nA$.
Let $I_n=\int_{\R^n_+}U^2$. Since $U(x)\asymp(1+r)^{2-n}$,
\[
 \int_{\R^n_+\cap\{r>1\}}U^2
       \asymp\int_1^\infty r^{3-n}\dd r
       =\frac1{n-4}<\infty\quad\text{for }n>4.
\]
For $n=5$, $k=0$, so Proposition~\ref{prop:upper-criterion}
already establishes the upper bound. For $6\le n\le9$, one has $k=1$,
and the criterion contains only $H_1$.
In the explicit calculations below, write
$D=|x'|^2+(1+x_n)^2$.

\begin{lemma}
For every $n\ge5$, the first coefficient $H_1=-2x_nA$ satisfies
\begin{equation}\label{eq:first-jet-sign}
 \cQ(H_1)=-2\mathcal C(H_1)
 \ge\frac{(n-2)(n-4)}{4n(n-1)}I_n|A|^2.
\end{equation}
\end{lemma}
\begin{proof}
If $A=0$, then $H_1=\psi_{H_1}=0$ and
\eqref{eq:first-jet-sign} holds. Assume $A\ne0$.
Since $\diver H_1=0$ and $|\partial H_1|^2=4|A|^2$,
\eqref{eq:direct-form} applies with two nonzero terms.
Tangential rotational symmetry and \eqref{eq:beta-rule} imply
\[
 \int_{\R^n_+}x_n^2\partial_iU\,\partial_jU
       =\frac{I_n}{2(n-1)}\delta_{ij}\quad\text{for }i,j<n.
\]
Contracting with $(A^2)_{ij}$ yields
\[
 J_2(H_1)=\frac{I_n|A|^2}{2(n-1)}
           -\frac{n-2}{8(n-1)}I_n|A|^2
         =\frac{6-n}{8(n-1)}I_n|A|^2.
\]

Let $\psi=\psi_{H_1}$. Formula~\eqref{eq:metric-correction-source}
gives
\[
 F_{H_1}=2n(n-2)x_n({x'}^tAx')D^{-(n+2)/2},
\]
and $F_{H_1}^\partial=0$ on $\Gamma$.
To bound $Q(\psi)$ from below, use the test function
\[
 \phi=x_n({x'}^tAx')D^{-n/2},
\]
with the same tangential angular factor.
It has zero trace on $\Gamma$ and finite energy, since
$|\nabla\phi|=O(r^{2-n})$ at infinity.
Differentiating explicitly, we find
\[
 -\Delta\phi=2n(2x_n+1)({x'}^tAx')D^{-(n+2)/2}.
\]
The equation and Robin boundary condition for $\psi$, with
$U^{2/(n-2)}=D^{-1}$, give
\[
 Q(\psi,\phi)
 =-\int\phi\Delta\psi
       +\int_\Gamma\phi(\partial_\nu\psi-nD^{-1}\psi)
 =\int F_{H_1}\phi.
\]
Since $\phi|_\Gamma=0$, integration by parts and
\eqref{eq:beta-rule} give
\[
 Q(\phi)=\int|\nabla\phi|^2=\int\phi(-\Delta\phi)
                =\frac{nI_n|A|^2}{2(n-1)(n-2)^2}.
\]
The same integral formula applied to the source yields
\[
 \int F_{H_1}\phi
                  =\frac{I_n|A|^2}{(n-1)(n-2)}.
\]

The identity $\Delta_{x'}({x'}^tAx')=2\tr A=0$ shows that the
source has tangential harmonic degree two. Angular projection commutes
with the correction equation, so uniqueness in Lemma~\ref{lem:jet-inverse}
forces the same degree for $\psi$. The test function $\phi$ has the same
angular factor.
Under the Cayley map this sector contains only spherical degrees
$\ell\ge2$ on $\partial\mathcal B_1$.
For its ball functions $v$, the Steklov decomposition in
\eqref{eq:ball-robin-form}, including the Dirichlet energy of the
zero-trace part, implies
\[
 \int_{\mathcal B_1}|\nabla v|^2-\int_{\partial\mathcal B_1}v^2
 \ge\frac12\int_{\mathcal B_1}|\nabla v|^2>0\qquad\text{for }v\ne0.
\]
Thus $Q$ is positive definite on the subspace containing $\psi$ and $\phi$.
By Cauchy--Schwarz,
\[
 Q(\psi)\ge\frac{Q(\psi,\phi)^2}{Q(\phi)}
          =\frac{(\int F_{H_1}\phi)^2}{Q(\phi)}
          =\frac{2I_n|A|^2}{n(n-1)}.
\]
By \eqref{eq:reduced-form},
$\cQ(H_1)=Q(\psi)-2J_2(H_1)$.
Substituting the two estimates proves \eqref{eq:first-jet-sign}.
\end{proof}

For $6\le n\le9$, $|H_1(z)|\asymp|A(z)|$, so the lemma gives
\[
 \cQ(\eps H_1(z))
     =\eps^2\cQ(H_1(z))
     \ge c\eps^2|H_1(z)|^2=c|H(z)|_\eps^2.
\]
This is \eqref{eq:geometric-upper-criterion}, uniformly at nearby
centers. Proposition~\ref{prop:upper-criterion} establishes the upper bound
in these dimensions.

For $n\ge10$, the same estimate applies at a nonumbilic point.
If $A(0)\ne0$, then $|A(z)|\ge c>0$ at nearby centers and the
higher coefficients are uniformly bounded. Hence
\[
 \cQ\left(\sum_{d=1}^k\eps^dH_d(z)\right)
 \ge c\eps^2|A(z)|^2-C\eps^3
 \ge c'\eps^2\ge c''|H(z)|_\eps^2.
\]
Proposition~\ref{prop:upper-criterion} applies here as well.

When $A(0)=0$, higher jets enter the dilation form. We identify
$\mathcal C$ with the constrained ball Hessian to determine its
sign and kernel.

Fix $\theta=(\lambda,\xi)$ and write
$U_\theta=U_{\lambda,\xi}$,
$\psi_{H;\theta}=\psi_{H;\lambda,\xi}$.
Set
\[
 T_\theta(x)=T\left(\frac{x'-\xi}{\lambda},\frac{x_n}{\lambda}\right),
 \qquad T_\theta^*\delta=4U_\theta^{4/(n-2)}\delta.
\]
For compactly supported Fermi-preserving $H$, consider the metric path
\[
 \widehat g_\eta=4(T_\theta^{-1})^*
 \bigl[(U_\theta+\eta\psi_{H;\theta})^{4/(n-2)}
                                      (\delta+\eta H)\bigr].
\]
It satisfies $\widehat g_0=\delta$, and its first variation is
\[
 \bar h_\theta=(T_\theta^{-1})^*
 \left[4U_\theta^{4/(n-2)}
 \left(H+\frac4{n-2}\frac{\psi_{H;\theta}}{U_\theta}\delta\right)\right].
\]
By \eqref{eq:jet-correction} at $U_\theta$,
$DR_\delta[\bar h_\theta]=Dh_\delta[\bar h_\theta]=0$.
The conformal energy identity is
\[
 J_{\delta+\eta H}(U_\theta+\eta\psi_{H;\theta})
       =\frac{n-2}{2^{n+1}(n-1)}S(\widehat g_\eta).
\]
Write $S''$ for the Hessian at the Euclidean metric,
\[
 S''[p,q]=\left.
 \frac{\partial^2}{\partial s\,\partial t}
 S(\delta+sp+tq)\right|_{s=t=0}.
\]
Since $S'_\delta=0$, matching the $\eta^2$ coefficients yields
\begin{equation}\label{eq:compactified-coefficient}
 \mathcal C(H;\theta)
 =\frac{n-2}{2^{n+2}(n-1)}S''[\bar h_\theta,\bar h_\theta].
\end{equation}
The bilinear identity follows by polarization.

For a smooth symmetric tensor $h$ on $\overline{\mathcal B_1}$, let
$\sigma$ be the normalized solution of \eqref{eq:scalar-completion}
and put $\bar h=h+\sigma\delta$. Let $X$ be the normalized solution
of \eqref{eq:tangent-vector-system}, and put
$k=\bar h-\mathcal L_X\delta$, $\tau=\tr k$ and
$k_{\nu\nu}=k(\nu,\nu)$.

\begin{lemma}\label{lem:ball-hessian}
The tensors $\bar h$ and $k$ satisfy
\[
 \begin{cases}
 DR_\delta[\bar h]=0,\quad\diver k=0&\text{in }\mathcal B_1,\\
 Dh_\delta[\bar h]=0,\quad k(\nu,T)=0&\text{on }\partial\mathcal B_1,
 \end{cases}
\]
where $T$ is any tangent vector and $X\cdot\nu=0$.
The Hessian of $S$ satisfies
\begin{equation}\label{eq:ball-hessian}
    \begin{aligned}
 S''[\bar h,\bar h]
& =-\frac12\int_{\mathcal B_1}|\nabla k|^2
   +\frac12\int_{\partial\mathcal B_1}
                  (3k_{\nu\nu}\tau-(n+2)k_{\nu\nu}^2)\\
&\le
 -\frac12\|\nabla k\|_2^2
 -\frac{n+2}{2}\|k_{\nu\nu}\|_{L^2(\partial\mathcal B_1)}^2
 \le-c\|k\|_{H^1(\mathcal B_1)}^2.
 \end{aligned}
\end{equation}
The completed quadratic form $h\mapsto S''[\bar h,\bar h]$
vanishes precisely for
\[
 h=f\delta+\mathcal L_Y\delta,\qquad
 Y\cdot\nu=0\quad\text{on }\partial\mathcal B_1.
\]
The identities and kernel description extend continuously to
$h\in H^1$, with $f\in H^1$ and $Y\in H^2$.
\end{lemma}

\begin{proof}
Appendix~\ref{app:ball-projections} constructs the corrections
and proves their Sobolev bounds.
Since $R_\delta=0$ and $h_\delta=1$, a boundary-tangent Lie
derivative preserves both linearized curvature constraints.
Diffeomorphism invariance and $S'_\delta=0$ therefore imply
$S''[\bar h,\bar h]=S''[k,k]$.
We write $A:B=\sum_{i,j}A_{ij}B_{ij}$ in an orthonormal frame,
using only tangential indices for tensors on the boundary.
For $DR_\delta[k]=Dh_\delta[k]=0$, differentiating the first
variation of $S$ yields
\begin{equation}\label{eq:ball-stress-hessian}
 S''[k,k]=-\int_{\mathcal B_1}D\mathrm{Ric}_\delta[k]:k
         +\int_{\partial\mathcal B_1}
                       (k_T-D\mathrm{II}_\delta[k]):k_T.
\end{equation}
The equations $\diver k=0$ and $DR_\delta[k]=0$ imply
\[
 D\mathrm{Ric}_\delta[k]
       =-\tfrac12\Delta k-\tfrac12\nabla^2\tau,
 \qquad \Delta\tau=0.
\]

On the unit sphere, $\nu=x$ and $\gamma=\delta_T$.
Use a tangential orthonormal frame $e_a$ which is geodesic at the
chosen boundary point and parallel along the normal.
Since $k(\nu,T)=0$, the normal variation is
$D\nu[k]=-\tfrac12k_{\nu\nu}\nu$.
Differentiating $\mathrm{II}(T_1,T_2)
=\langle\nabla_{T_1}\nu,T_2\rangle$ gives
\[
 D\mathrm{II}_\delta[k]
 =k_T+\tfrac12(\nabla_\nu k)_T-\tfrac12k_{\nu\nu}\gamma.
\]
Taking the trace, including the inverse-metric variation, yields
\[
 2(n-1)Dh_\delta[k]
 =\partial_\nu\tau-\partial_\nu k_{\nu\nu}
                         -(n-1)k_{\nu\nu}=0.
\]
The normal divergence equation is
\[
 \begin{aligned}
 0=(\diver k)(\nu)
 &=\partial_\nu k_{\nu\nu}
                   +\sum_a(\nabla_{e_a}k)(e_a,\nu)\\
 &=\partial_\nu k_{\nu\nu}+(n-1)k_{\nu\nu}-\tr_\gamma k_T\\
 &=\partial_\nu k_{\nu\nu}+nk_{\nu\nu}-\tau.
 \end{aligned}
\]
Thus
\begin{equation}\label{eq:normal-trace-relations}
 \partial_\nu k_{\nu\nu}=\tau-nk_{\nu\nu},\qquad
 \partial_\nu\tau=\tau-k_{\nu\nu}.
\end{equation}

Integration by parts in the interior term of
\eqref{eq:ball-stress-hessian}, using $\diver k=0$, yields
\[
 \begin{aligned}
 -\int_{\mathcal B_1}D\mathrm{Ric}_\delta[k]:k
 ={}&-\frac12\int_{\mathcal B_1}|\nabla k|^2\\
 &+\frac12\int_{\partial\mathcal B_1}
 \bigl(k_T:\nabla_\nu k_T
       +k_{\nu\nu}\partial_\nu k_{\nu\nu}
       +k_{\nu\nu}\partial_\nu\tau\bigr).
 \end{aligned}
\]
The boundary term is
\[
 \int_{\partial\mathcal B_1}
          (k_T-D\mathrm{II}_\delta[k]):k_T
 =-\frac12\int_{\partial\mathcal B_1}k_T:\nabla_\nu k_T
   +\frac12\int_{\partial\mathcal B_1}
                       k_{\nu\nu}(\tau-k_{\nu\nu}).
\]
The tangential derivative terms cancel.
Substituting \eqref{eq:normal-trace-relations} proves
\eqref{eq:ball-hessian}.

Since $\diver(kx)=\tau$ and $\tau$ is harmonic,
\[
 \int_{\partial\mathcal B_1}k_{\nu\nu}
       =\int_{\mathcal B_1}\tau
       =\frac1n\int_{\partial\mathcal B_1}\tau.
\]
Also, integrating the second relation in
\eqref{eq:normal-trace-relations} yields
\[
 0=\int_{\partial\mathcal B_1}\partial_\nu\tau
   =\int_{\partial\mathcal B_1}\tau
                    -\int_{\partial\mathcal B_1}k_{\nu\nu}.
\]
Hence both boundary means are zero.
For a harmonic function with zero boundary mean, the Steklov inequality reads
$\int_{\partial\mathcal B_1}\tau^2
 \le\int_{\mathcal B_1}|\nabla\tau|^2$.
Therefore
\[
 \int_{\partial\mathcal B_1}k_{\nu\nu}\tau
 =\int_{\partial\mathcal B_1}
                    \tau(\tau-\partial_\nu\tau)
 =\int_{\partial\mathcal B_1}\tau^2
                  -\int_{\mathcal B_1}|\nabla\tau|^2\le0.
\]
This proves the asserted upper bound for the Hessian.
For constant symmetric matrices, the map
$P\mapsto(x\mapsto x^tPx|_{\partial\mathcal B_1})$ is injective.
The componentwise Poincar\'e and trace inequalities consequently imply
\[
 \|k\|_{H^1}^2\le C\bigl(\|\nabla k\|_2^2
                    +\|k_{\nu\nu}\|_{L^2(\partial\mathcal B_1)}^2\bigr),
\]
which proves the coercive bound.

For $h_j\to h$ in $H^1$, the bounded corrections imply
$\bar h_j\to\bar h$ and $k_j\to k$ in $H^1$.
The boundary trace theorem passes the identity and inequalities
to the limit.
Vanishing of the Hessian forces $k=0$, hence
$h=-\sigma\delta+\mathcal L_X\delta$.
Conversely, if $h=f\delta+\mathcal L_Y\delta$ with $Y$ tangent,
the scalar equation implies $\sigma+f=a\cdot x$ for some $a\in\R^n$.
The field
\[
 Y_a(x)=(1+|x|^2)a-2(a\cdot x)x
\]
is tangent on the sphere and satisfies
$\mathcal L_{Y_a}\delta=-4(a\cdot x)\delta$.
Thus $\bar h=\mathcal L_{Y-Y_a/4}\delta$, and the vector
minimization implies $k=0$.
\end{proof}

For polynomial $H$ of maximal degree $d<(n-2)/2$, cut off its
compactification near the omitted boundary point to pass the
identity to the limit.
If $r$ is distance to that point, the transformed tensor and its
first derivative have sizes
$O(r^{-d})$ and $O(r^{-d-1})$.
The squared $H^1$ cutoff error near that point is bounded by
\begin{equation}\label{eq:compactification-tail}
 C\int_0^\rho r^{n-2d-3}\dd r
        =\frac{C\rho^{n-2d-2}}{n-2d-2}.
\end{equation}
Letting $\rho\downarrow0$ and using the bounded corrections extends
\eqref{eq:compactified-coefficient} to these tensors, uniformly on
compact bubble-parameter sets.

Let $X_\theta$ solve \eqref{eq:tangent-vector-system} for
$\bar h_\theta$, and set
$k_\theta=\bar h_\theta-\mathcal L_{X_\theta}\delta$.
The degree-one scalar freedom is a boundary-tangent Lie
derivative, as in the proof of Lemma~\ref{lem:ball-hessian}.
That lemma and \eqref{eq:compactified-coefficient} give
\[
 \mathcal C(H;\theta)
       \le-c\|k_\theta\|_{H^1(\mathcal B_1)}^2.
\]

\begin{lemma}\label{lem:homogeneous-kernel}
Let $H_d$ be a homogeneous Fermi-preserving symmetric tensor of degree
$1\le d<(n-2)/2$. Then $\mathcal C(H_d;\theta)\le0$ for every bubble
parameter $\theta$. Equality at one parameter holds if and only if
\begin{equation}\label{eq:removable-homogeneous-jet}
 H_d=\mathcal L_{Y_{d+1}}\delta+f_d\delta,\qquad
 (Y_{d+1})_n=0\quad\text{on }\Gamma,
\end{equation}
where $Y_{d+1}$ and $f_d$ are homogeneous of degrees $d+1$ and $d$.
In this case equality holds at every bubble parameter.
\end{lemma}
\begin{proof}
Nonpositivity follows from Lemma~\ref{lem:ball-hessian} and
\eqref{eq:compactified-coefficient}. If equality holds, then
$k_\theta=0$ and $\bar h_\theta=\mathcal L_{X_\theta}\delta$.
Pull back $X_\theta$ by $T_\theta$ to obtain a field $Y$ tangent to
$\Gamma$. By conformal covariance,
\[
 H_d+\frac4{n-2}\frac{\psi_{H_d;\theta}}{U_\theta}\delta
 =\mathcal L_Y\delta+\frac4{n-2}Y(\log U_\theta)\delta.
\]
The correction and the vector equation have smooth data near
$T_\theta(0)$. Local boundary regularity \cite{ADN} makes $Y$ and
$f=\frac4{n-2}(Y(\log U_\theta)-\psi_{H_d;\theta}/U_\theta)$ smooth
near zero. Taking their homogeneous Taylor terms proves
\eqref{eq:removable-homogeneous-jet}.

Conversely, assume \eqref{eq:removable-homogeneous-jet} and set
$\psi=Y_{d+1}(U_\theta)-\frac{n-2}{4}f_dU_\theta$. Then
\[
 U_\theta^{4/(n-2)}
 \left(H_d+\frac4{n-2}\frac{\psi}{U_\theta}\delta\right)
 =\mathcal L_{Y_{d+1}}(U_\theta^{4/(n-2)}\delta).
\]
The correction equation follows from diffeomorphism and conformal covariance.
Since $|\nabla\psi|=O(|x|^{d+1-n})$, the condition $2d<n-2$ ensures
$\psi\in\mathcal D^{1,2}$.
Imposing \eqref{eq:parameter-corrector-normalization} changes $\psi$
by a bubble Jacobi field and leaves the correction energy unchanged.
After compactification, the completed variation is
$\mathcal L_{X_\theta}\delta$, with
$X_\theta=(T_\theta)_*Y_{d+1}$ and
$|D^jX_\theta|=O(r^{1-d-j})$ for $j\le2$ near the omitted point.
The cutoff estimate \eqref{eq:compactification-tail} therefore implies
\[
 \|\mathcal L_{\chi X_\theta}\delta-
             \mathcal L_{X_\theta}\delta\|_{H^1(\mathcal B_1)}\to0
\]
for smooth cutoffs vanishing near that point and tending to one.
The fields $\chi X_\theta$ are boundary tangent, so diffeomorphism
invariance and continuity imply
$S''[\mathcal L_{X_\theta}\delta,\mathcal L_{X_\theta}\delta]=0$.
Equation~\eqref{eq:compactified-coefficient} proves equality for every
$\theta$.
\end{proof}

For several homogeneous degrees, bilinearity and
$\cQ(H_i,H_j)=-(i+j)\mathcal C(H_i,H_j)$ yield
\[
 \cQ\left(\sum_{d=1}^k\eps^dH_d\right)
 =-\sum_{i,j=1}^k(i+j)\eps^{i+j}\mathcal C(H_i,H_j).
\]
In the normalized logarithmic expansion, $A=\mathrm{II}^{\circ}$
and $\overline W$ is the boundary Weyl tensor.
Modulo the determined products of lower-degree coefficients, the
homogeneous tensors are
\begin{equation}\label{eq:cubic-fermi-slice}
\begin{aligned}
 H_1&=-2x_nA,\\
 H_2&\equiv H_{\overline W}(x')-2x_nT_1(x')-2x_n^2S,\\
 H_3&\equiv C(x')-2x_nP(x')-2x_n^2T_2(x')-2x_n^3A_3.
\end{aligned}
\end{equation}
Here $H_{\overline W}$ and $C$ are the quadratic and cubic
boundary normal-coordinate terms, with
$H_{\overline W}x'=Cx'=0$.
The normal Taylor coefficients $T_1,T_2$ are linear in $x'$,
$P$ is quadratic, and $S,A_3$ are constant matrices.
All these tensors are tangential, symmetric and trace free.
The trace normalization and the conformal and boundary-coordinate
freedom are fixed as in Lemma~\ref{lem:fermi-normalization}.
The determined terms have coefficient sizes $O(|A|^2)$ in degree
two and $O(|A||H_2|+|A|^3)$ in degree three.
The norms below measure the independent coefficients, in particular
\[
 |H_2|^2=|\overline W|^2+|T_1|^2+|S|^2.
\]

Since $k=\lfloor(n-2)/4\rfloor$, one has $k=2$
for $10\le n\le13$ and $k=3$ for $14\le n\le17$.
The cubic estimate below also holds at $n=18$, where $k=4$.
The degree-four terms are treated in
Section~\ref{subsec:eighteen-upper}.

\begin{proposition}\label{prop:metric-coercivity}
There are $c,\eps_0>0$, uniform at the chosen centers, such that
the complete normalized Fermi tensors, with independent parts as in
\eqref{eq:cubic-fermi-slice}, satisfy the following estimates for
$0<\eps<\eps_0$. For $10\le n\le13$,
\begin{equation}\label{eq:quadratic-coercivity}
 \cQ(\eps H_1+\eps^2H_2)\ge
 c(\eps^2|A|^2+\eps^4|H_2|^2).
\end{equation}
For $14\le n\le18$,
\begin{equation}\label{eq:cubic-coercivity}
 \cQ(\eps H_1+\eps^2H_2+\eps^3H_3)
 \ge c\sum_{j=1}^3\eps^{2j}|H_j|^2.
\end{equation}
\end{proposition}
\begin{proof}
\emph{Quadratic estimate, $10\le n\le13$.}
First consider the independent parts. The degree-two terms are
$H_{\overline W},-2x_nT_1,-2x_n^2S$.
Reflection $x'\mapsto-x'$ preserves the correction equation and
the normalization \eqref{eq:corrector-normalization}.
Uniqueness in Lemma~\ref{lem:jet-inverse}
therefore makes the corrections inherit the parity of their tensors.
The $T_1$ term is odd and the other terms are even, so
\[
 \mathcal C(x_nT_1,H_{\overline W})
 =\mathcal C(x_nT_1,x_nA)
 =\mathcal C(x_nT_1,x_n^2S)=0.
\]
For the Weyl term $H_{\overline W}$,
\[
 \tr H_{\overline W}=0,\qquad
 \diver_{x'}H_{\overline W}=0,\qquad
 H_{\overline W}x'=0,\qquad \partial_nH_{\overline W}=0.
\]
Thus \eqref{eq:metric-correction-source} implies
$F_{H_{\overline W}}=F_{H_{\overline W}}^\partial=0$, hence
$\psi_{H_{\overline W}}=0$.
The vanishing contractions of the Weyl tensor then imply
$\mathcal C(H_{\overline W},x_nA)
 =\mathcal C(H_{\overline W},x_n^2S)=0$.
For the Weyl term, Lemma~\ref{lem:ball-hessian} and
\eqref{eq:removable-homogeneous-jet} give
$\mathcal C(H_{\overline W})<0$ when $\overline W\ne0$.

For the linear symmetric trace-free tensor, write
$(T_1)_{ij}(x')=(T_1)_{ij,k}x'_k$.
Applying \eqref{eq:beta-rule} to \eqref{eq:direct-form} cancels the
terms containing $|\diver_{x'}T_1|^2$ and yields
\[
 \frac{J_2(-2x_nT_1)}{I_n}
 =\frac{n-2}{8(n-5)(n-1)}
       \left(-|T_1|^2+\frac4{n-6}
                     \sum(T_1)_{ij,k}(T_1)_{ik,j}\right).
\]
By Cauchy--Schwarz,
\[
    \sum_{i,j,k}(T_1)_{ij,k}(T_1)_{ik,j}
    \le \Bigl(\sum_{i,j,k}(T_1)_{ij,k}^2\Bigr)^{1/2}
    \Bigl(\sum_{i,j,k}(T_1)_{ik,j}^2\Bigr)^{1/2}
    =|T_1|^2.
\]
Hence
\[
-|T_1|^2+\frac4{n-6}\sum_{i,j,k}(T_1)_{ij,k}(T_1)_{ik,j}
\le \Bigl(-1+\frac4{n-6}\Bigr)|T_1|^2
=-\frac{n-10}{n-6}|T_1|^2\le0
\]
for \(n\ge10\). Therefore
\[
J_2(-2x_nT_1)\le0,
\]
with equality at \(n=10\) only if
\((T_1)_{ij,k}=(T_1)_{ik,j}\).
The Jacobi form is nonnegative on odd functions by
\eqref{eq:ball-robin-form}, hence
$\mathcal C(-2x_nT_1)\le J_2(-2x_nT_1)<0$ for $n>10$ and $T_1\ne0$.
At $n=10$, equality in the direct term requires
$(T_1)_{ij,k}=(T_1)_{ik,j}$.
Together with symmetry and trace freeness, this makes
\[
 Q_{T_1}=(T_1)_{ij,k}x'_ix'_jx'_k,\qquad
 \Delta_{x'}Q_{T_1}=0.
\]
If $T_1\ne0$, then $Q_{T_1}\ne0$ and
\[
 F_{-2x_nT_1}
       =2n(n-2)x_nQ_{T_1}D^{-(n+2)/2}\ne0.
\]
The tangential harmonic degree is three. The spectral decomposition
following \eqref{eq:ball-robin-form} makes $Q$ positive definite there. Thus $Q(\psi_{-2x_nT_1})>0$, establishing strict negativity
also at $n=10$.
Compactness of the unit coefficient spheres makes both estimates uniform.

It remains to control the mixed terms in $A$ and $S$.
Their corrections have tangential harmonic degree two.
Appendix~\ref{app:general-forms} proves the upper estimate
\begin{equation}\label{eq:first-correction-upper}
 Q(\psi_{H_1})\le
       \frac{n^2-n-14}{4(n-5)(n-1)}I_n|A|^2.
\end{equation}
Put $\psi_1=\eps\psi_{H_1}$ and
$\psi_2=\eps^2\psi_{-2x_n^2S}$.
The correction part of the dilation form satisfies
\[
 Q(\psi_1)+3Q(\psi_1,\psi_2)+2Q(\psi_2)
=2Q(\psi_2+\tfrac34\psi_1)-\tfrac18Q(\psi_1)
 \ge-\tfrac18\eps^2Q(\psi_{H_1}).
\]
Combining this bound with \eqref{eq:direct-form} yields
\[
 \frac{\cQ(\eps H_1-2\eps^2x_n^2S)}{I_n}
 \ge M_{11}\eps^2|A|^2
       +2M_{12}\eps^3 A:S+M_{22}\eps^4|S|^2,
\]
where
\begin{equation*}
 M=
 \begin{pmatrix}
 \dfrac{7n^2-87n+254}{32(n-5)(n-1)}&
 \dfrac{3(n-8)}{4(n-5)(n-1)}\\[4pt]
 \dfrac{3(n-8)}{4(n-5)(n-1)}&
 \dfrac{4(n-8)}{(n-6)(n-5)(n-1)}
 \end{pmatrix}.
\end{equation*}
Its first diagonal entry is positive, and
\[
 \det M=
 \frac{(n-8)(n-2)(5n-38)}
      {16(n-6)(n-5)^2(n-1)^2}>0.
\]
This proves \eqref{eq:quadratic-coercivity} for the independent coefficients.

\emph{Cubic estimate, $14\le n\le18$.}
The form $\cQ$ is invariant under the tangential
orthogonal group on $V=\R^{n-1}$.
It therefore pairs only equivalent irreducible representations,
so positivity can be checked on each isotypic component.
Here $[\lambda]$ denotes the trace-free representation with Young
diagram $\lambda$.
The linear coefficients $T_1,T_2$ have the decomposition
\[
 \mathrm{Sym}^2_0V\otimes V=[3]\oplus[2,1]\oplus[1].
\]
The quadratic coefficient $P$ belongs to
\[
 \mathrm{Sym}^2_0V\otimes\mathrm{Sym}^2V
 =[4]\oplus[3,1]\oplus[2,2]\oplus2[2]\oplus[1,1]\oplus[0].
\]
The intrinsic cubic term satisfies
\[
 \{C=C^t,\ \deg C=3,\ Cx'=0,\ \tr C=0\}
 =[3,2]\oplus[2,1].
\]
In addition, $A,S,A_3$ have type $[2]$ and
$H_{\overline W}$ has type $[2,2]$.
Consequently the representations occurring in more than one
homogeneous degree are exactly
\[
\begin{array}{c|c|c}
 \text{representation}&\text{degrees}&\text{multiplicities}\\\hline
 {[2]}&1,2,3&1,1,3\\
 {[2,2]}&2,3&1,1\\
 {[3]}&2,3&1,1\\
 {[2,1]}&2,3&1,2\\
 {[1]}&2,3&1,1.
\end{array}
\]
On each row, the invariant form is a finite matrix on the copies
of that representation, tensored with its irreducible inner product.

All other representations occur only in degree three. On such a
component, \(\cQ=-6\mathcal C\ge0\) by Lemma~\ref{lem:ball-hessian}
and \eqref{eq:compactified-coefficient}. If equality holds, then
Lemma~\ref{lem:ball-hessian} and
\eqref{eq:removable-homogeneous-jet} imply
\[
H_3=\mathcal L_{Y_4}\delta+f_3\delta,\qquad
(Y_4)_n=0\ \text{on }\Gamma.
\]
The trace, Fermi and boundary normal-coordinate conditions yield
the homogeneous equations of Lemma~\ref{lem:fermi-normalization}
for \(f_3/2\) and \(Y_4\), with zero determinant datum.
Invertibility of \(L_3\) and the normal recursion force \(f_3=0\),
and \eqref{eq:gauge-vector} and its boundary value give \(Y_4=0\).
Hence \(H_3=0\), contradicting the nonzero representative.
Thus \(\cQ>0\) on every nonzero vector in these components.
Finite-dimensional compactness ensures a positive lower bound for
\(\cQ\) on their direct sum.

The five $[2]$ copies are $A,S,A_3$ and the two $[2]$ components
of $P$, in degrees $1,2,3,3,3$. Their corrections consist of explicit
interior particular solutions and harmonic Robin solutions.
For boundary data $q$ smooth near zero, with
$\int_\Gamma Dq^2<\infty$ and $\int_\Gamma qZ_\alpha=0$ for
$0\le\alpha<n$, let $h_q$ be the finite-energy solution of
\begin{equation}\label{eq:harmonic-robin}
 \begin{cases}
 \Delta h_q=0&\text{in }\R^n_+,\\
 (\partial_\nu-nU^{2/(n-2)})h_q=q&\text{on }\Gamma.
 \end{cases}
\end{equation}
with $h_q(0)=0$ and $\nabla_{x'}h_q(0)=0$, as in
\eqref{eq:corrector-normalization}.
In the $[2]$ component the compactified source degrees are $j+2$,
$j\ge0$. Let $\Pi_j$ be the corresponding orthogonal projections
in $L^2(\Gamma,D\,\dd x')$ and put
$\Pi_{\le N}=\sum_{j=0}^N\Pi_j$.
The inverse eigenvalues in \eqref{eq:ball-robin-form} are
$1/[2(j+1)]$. Hence
\begin{equation}\label{eq:robin-tail-control}
 0\le R_N(q)
 :=Q(h_q)-Q(h_{\Pi_{\le N}q})
 =\sum_{j>N}\frac{\int_\Gamma D|\Pi_jq|^2}{2(j+1)}
 \le\frac{\int_\Gamma D|q-\Pi_{\le N}q|^2}{2(N+2)}.
\end{equation}
Let $q_i$ be the Robin datum of the degree-$i$ contribution,
including its factor $\eps^i$.
The dilation weights contribute the coefficients $i+j$ in the mixed
remainder terms. To bound this remainder by a finite quadratic form,
use the following identity for a nonnegative quadratic form $R$,
\begin{equation}\label{eq:spectral-completion}
 \begin{split}
 &R(q_1)+2R(q_2)+3R(q_3)
       +3R(q_1,q_2)+4R(q_1,q_3)+5R(q_2,q_3)\\
 &=3R(q_3+\tfrac23q_1+\tfrac56q_2)
                       -\tfrac1{12}R(2q_1+q_2).
 \end{split}
\end{equation}
For the $[2]$ block, take $R=R_N$. Dropping the nonnegative square
in \eqref{eq:spectral-completion} bounds the spectral remainder below by
$-\frac1{12}R_N(2q_1+q_2)$. Estimate~\eqref{eq:robin-tail-control}
controls this term. Combining it with the particular solutions in
Appendix~\ref{app:general-forms} proves \eqref{eq:five-copy-lower}.
The same appendix supplies the matrices for the other repeated
blocks \(([2,2])\), \(([3])\), \(([2,1])\) and \(([1])\). Their
Schur pivots satisfy \eqref{eq:metric-schur-bounds} and the
smaller-matrix estimates in that appendix. Together with the
single-degree components, this proves \eqref{eq:cubic-coercivity}
for the independent coefficients.

It remains to include the determined products. Write
\[
H_d=H_d^{\mathrm{ind}}+H_d^{\mathrm{det}}\qquad\text{for }d=1,2,3,
\]
where \(H_d^{\mathrm{ind}}\) is the independent part and
\(H_d^{\mathrm{det}}\) is the product of lower-degree coefficients
fixed by the normalizations. By the coefficient bounds preceding
the proposition, on the fixed bounded coefficient set,
\[
\eps^d|H_d^{\mathrm{det}}|
\le C\eps\Bigl(\sum_{j=1}^{d-1}\eps^{2j}|H_j|^2\Bigr)^{1/2}
\qquad\text{for }d=2,3.
\]
Consequently, for \(s=2\) or \(3\),
\[
    \Bigl|\cQ\Bigl(\sum_{d=1}^s\eps^dH_d^{\mathrm{det}},
    \sum_{d=1}^s\eps^dH_d^{\mathrm{ind}}\Bigr)\Bigr|
    +\Bigl|\cQ\Bigl(\sum_{d=1}^s\eps^dH_d^{\mathrm{det}},
    \sum_{d=1}^s\eps^dH_d^{\mathrm{det}}\Bigr)\Bigr|
    \le C\eps\sum_{j=1}^s\eps^{2j}|H_j|^2.
\]
Adding this to the independent-coefficient bound yields
\[
\cQ\Bigl(\sum_{d=1}^s\eps^dH_d\Bigr)
\ge(c-C\eps)\sum_{j=1}^s\eps^{2j}|H_j|^2.
\]
After decreasing \(\eps_0\), both \eqref{eq:quadratic-coercivity} and
\eqref{eq:cubic-coercivity} follow for the full normalized tensors.
\end{proof}

For $5\le n\le17$, the preceding estimates give
\eqref{eq:geometric-upper-criterion}.
The critical linear terms at $n=6,8$ vanish by the cutoff calculation
in Lemma~\ref{lem:taylor-budget}.
Proposition~\ref{prop:upper-criterion} therefore proves the
general-metric upper estimate.

For Euclidean domains, a single boundary graph function determines
the metric variation. Its shape form extends the upper estimate to
$18\le n\le22$.
At each selected center, Euclidean motions put the boundary in the local graph form
$x_n=f(x')$ with $f(0)=0$ and $\nabla f(0)=0$.
Write $f_2=f_2^\circ+c|x'|^2$, with $\Delta f_2^\circ=0$.
To remove the quadratic trace, use the special conformal map
\[
 M_c(x)=\frac{x-c|x|^2e_n}{1-2cx_n+c^2|x|^2}
\]
which has $M_c(0)=0$, $DM_c(0)=I$, and
\[
 M_c^*\delta=\varphi_c^{4/(n-2)}\delta,\qquad
 \varphi_c(x)=(1-2cx_n+c^2|x|^2)^{-(n-2)/2}.
\]
It changes the quadratic graph term to
$\widetilde f_2=f_2-c|x'|^2=f_2^\circ$.
For the Euclidean representative of the solution, put
\[
 u_c(M_cx)=\varphi_c(x)^{-1}u(x).
\]
Conformal covariance preserves \eqref{eq:main} on the transformed
Euclidean domain.
Since $\varphi_c(0)=1$ and $\nabla_T\varphi_c(0)=0$, the selected
value $u(0)$ and the condition $\nabla_Tu(0)=0$ are preserved.
The maps and their inverses have uniform smooth bounds at nearby
centers. Retain $f,u$ for the normalized graph and solution.

Let $\Phi$ be the boundary Fermi parametrization of this domain,
so that $g=\Phi^*\delta$.
Let $\gamma$ and $\mathrm{II}$ be the induced boundary metric and
second fundamental form.
The normal collar and the Gauss equation imply
\begin{equation}\label{eq:euclidean-collar}
 \begin{aligned}
 g_{x_n}&=\gamma-2x_n\mathrm{II}
                    +x_n^2\mathrm{II}\gamma^{-1}\mathrm{II},\\
 R^\gamma_{ijkl}
        &=\mathrm{II}_{ik}\mathrm{II}_{jl}
                         -\mathrm{II}_{il}\mathrm{II}_{jk}.
 \end{aligned}
\end{equation}
For the graph $x_n=\eta f(x')$,
$\gamma_\eta=\delta+\eta^2\dd f\otimes\dd f$ and
$\left.\partial_\eta\mathrm{II}_\eta\right|_{\eta=0}=D^2f$.
Differentiating the collar formula produces the first metric variation
\begin{equation}\label{eq:graph-metric-variation}
 (h_f)_{ij}=
 \begin{cases}
 -2x_n\partial_{ij}f,&i,j<n,\\
 0,&i=n\text{ or }j=n.
 \end{cases}
\end{equation}
It also satisfies $h_f=\mathcal L_{Y_f}\delta$, where
$Y_f(x)=(-x_n\nabla f(x'),f(x'))$.
Thus a graph term $f_{d+1}$ of degree $d+1$ induces a metric term
$h_{f_{d+1}}$ of degree $d$.

Since $Y_f|_\Gamma=fe_n$, the normal conformal factor of the Cayley
map $T_\theta$, for $\theta=(\lambda,\xi)$, determines the outward normal
component $-\widehat f_\theta$, where
\[
 \widehat f_\theta(z)
       =\frac{2}{1+|z|^2}\lambda^{-1}f(\xi+\lambda z).
\]

Let $\Lambda$ be the Dirichlet-to-Neumann map of $\mathcal B_1$ and define
\begin{equation}\label{eq:shape-form}
 Q_\partial(f,g)
 =\langle f,\Lambda(-\Delta_{\Sph^{n-1}}-(n-1))g\rangle,
\end{equation}
with the $L^2(\Sph^{n-1})$ inner product.
To compute the constrained Hessian in terms of this form, take a
smooth normal displacement $\zeta$ of the unit sphere,
choose a smooth field $X$ with $X=\zeta\nu$ on the boundary and put
$h=\mathcal L_X\delta+\sigma\delta$.
The curvature constraints $DR_\delta[h]=Dh_\delta[h]=0$ are
\[
 \begin{cases}
 \Delta\sigma=0&\text{in }\mathcal B_1,\\
 -(\Delta_{\Sph^{n-1}}+(n-1))\zeta
       +\dfrac{n-1}{2}(\partial_\nu\sigma-\sigma)=0
                                     &\text{on }\Sph^{n-1}.
 \end{cases}
\]
For the spherical harmonic components $\zeta_\ell,\sigma_\ell$,
$\Lambda\sigma_\ell=\ell\sigma_\ell$ and
$-\Delta_{\Sph^{n-1}}\zeta_\ell=\ell(\ell+n-2)\zeta_\ell$.
The boundary equation therefore reads
\[
 (\ell-1)\left((\ell+n-1)\zeta_\ell
                         +\frac{n-1}{2}\sigma_\ell\right)=0.
\]
Fix the degree-one Robin kernel by $\sigma_1=0$.
For $\ell\ne1$,
$\sigma_\ell=-2(\ell+n-1)\zeta_\ell/(n-1)$.
Two extensions of $X$ differ by a field $Y$ vanishing on the boundary.
Diffeomorphism invariance and $S'_\delta=0$ imply
$S''[\mathcal L_Y\delta,h]=0$, so the Hessian depends only on
$\zeta$.

On the unit sphere, $\mathrm{II}_\delta=\gamma=\delta_T$, so
$h_T=(2\zeta+\sigma)\gamma$.
Differentiating $(n-1)h_g=\tr_{g_T}\mathrm{II}_g$ gives
\[
 \tr_\gamma D\mathrm{II}_\delta[h]-\tr_\gamma h_T
       =(n-1)Dh_\delta[h]=0.
\]
Since $h_T$ is a scalar multiple of $\gamma$, this trace identity
makes the boundary term in \eqref{eq:ball-stress-hessian} zero.
Since $\mathrm{Ric}_\delta=0$ and $\Delta\sigma=0$,
$D\mathrm{Ric}_\delta[h]=-\frac{n-2}{2}\nabla^2\sigma$.
Integrating by parts, we find
\[
 \begin{aligned}
 S''[h,h]
 &=\frac{n-2}{2}\int_{\mathcal B_1}
                  \nabla^2\sigma:(\mathcal L_X\delta+\sigma\delta)
 =(n-2)\int_{\mathcal B_1}
                         \partial_{ij}\sigma\,\partial_iX_j\\
 &=(n-2)\int_{\Sph^{n-1}}\nabla^2\sigma(\nu,X)
 =(n-2)\sum_\ell\ell(\ell-1)
                                  \langle\sigma_\ell,\zeta_\ell\rangle.
 \end{aligned}
\]
A degree-one addition to $\sigma$ contributes zero to this sum.
Substitution of $\sigma_\ell$ yields
\begin{equation}\label{eq:shape-hessian}
 S''[h,h]=-\frac{2(n-2)}{n-1}Q_\partial(\zeta,\zeta).
\end{equation}
The eigenvalue of $Q_\partial$ in degree $\ell$ is
$\ell(\ell-1)(\ell+n-1)$.
Hence its kernel consists of constants and first spherical harmonics,
and
$|Q_\partial(f,g)|\le C\|f\|_{H^{3/2}}\|g\|_{H^{3/2}}$.
Density extends \eqref{eq:shape-hessian} to $H^{3/2}(\Sph^{n-1})$.

For a graph polynomial of maximal degree $d+1$ with $2d<n-2$,
set $X=(T_\theta)_*Y_f$.
If $r$ is distance to the omitted point, then
$|D^jX|=O(r^{1-d-j})$ for $j=0,1,2$.
Thus $X\in H^2(\mathcal B_1)$ and its normal trace
$\zeta=-\widehat f_\theta$ belongs to $H^{3/2}(\Sph^{n-1})$.
Since $h_f=\mathcal L_{Y_f}\delta$, its compactified and
scalar-corrected variation, with
$DR_\delta[\bar h_\theta]=Dh_\delta[\bar h_\theta]=0$, has the form
\[
 \bar h_\theta=\mathcal L_X\delta+\sigma\delta.
\]
Since $S'_\delta=0$, \eqref{eq:compactified-coefficient} and
\eqref{eq:shape-hessian} imply
\begin{equation}\label{eq:shape-energy-coefficient}
 \mathcal C(h_f;\theta)
 =-\frac{(n-2)^2}{2^{n+1}(n-1)^2}
                    Q_\partial(\widehat f_\theta,\widehat f_\theta).
\end{equation}
At the unit bubble, write
$\widehat f_{d+1}=2f_{d+1}/(1+|z|^2)$.
With $k$ as in the upper criterion, put
\[
 \mathcal Q_\partial(f_2,\ldots,f_{k+1})
 =\frac{(n-2)^2}{2^{n+1}(n-1)^2}
       \sum_{i,j=1}^k(i+j)
             Q_\partial(\widehat f_{i+1},\widehat f_{j+1}).
\]

Polarization and $\cQ(H_i,H_j)=-(i+j)\mathcal C(H_i,H_j)$ give
\begin{equation}\label{eq:graph-metric-dilation}
 \cQ\left(\sum_{d=1}^k\eps^d h_{f_{d+1}}\right)
       =\mathcal Q_\partial(\eps f_2,\ldots,\eps^kf_{k+1}).
\end{equation}
\begin{proposition}\label{prop:euclidean-shape-form}
Let $18\le n\le22$, with $k=4$ for $n\le21$ and $k=5$ for $n=22$.
The homogeneous Taylor terms of a normalized graph with
$f(0)=0$, $\nabla f(0)=0$ and $\Delta f(0)=0$
satisfy, for every $\eps>0$,
\begin{equation}\label{eq:euclidean-coercivity}
 \mathcal Q_\partial(\eps f_2,\ldots,\eps^kf_{k+1})
 \ge c\sum_{d=1}^k\eps^{2d}\|f_{d+1}\|^2.
\end{equation}
Here each norm is a fixed coefficient norm on the corresponding
homogeneous polynomial space, and $c>0$ depends only on these norms
and $n$.
\end{proposition}
\begin{proof}
By the spectrum following \eqref{eq:shape-hessian}, the inverse images
of $\ker Q_\partial$ under $f\mapsto2f/(1+|z|^2)$ form
$\Span\{1,z_1,\ldots,z_{n-1},|z|^2\}$.
The graph normalizations remove these directions.

Apply the Fischer decomposition to the retained graph polynomial,
\[
 \sum_{d=1}^kf_{d+1}(z)
       =\sum_{\ell,j}|z|^{2j}H_{\ell j}(z),\qquad
 \Delta H_{\ell j}=0,\quad \deg H_{\ell j}=\ell.
\]
Orthogonality on $\Sph^{n-2}$ and rotational invariance of
$Q_\partial$ eliminate pairings with different $\ell$.
The degree restrictions
$2\le\ell+2j\le k+1$, excluding the quadratic trace
$(\ell,j)=(0,1)$, determine the admissible indices
\begin{equation}\label{eq:euclidean-radial-indices}
\begin{array}{c|cc}
 \ell&j\text{ for }18\le n\le21&j\text{ for }n=22\\\hline
 0&2&2,3\\
 1&1,2&1,2\\
 2&0,1&0,1,2\\
 3&0,1&0,1\\
 4&0&0,1\\
 5&0&0\\
 6&\varnothing&0
\end{array}
\end{equation}
Fix a basis vector $H_\ell$ in an orthonormal harmonic basis.
Let $j_0<\cdots<j_r$ be its allowed radial indices and set
$d_i=\ell+2j_i-1$.
For $\widehat f_{\ell j}=2|z|^{2j}H_\ell/(1+|z|^2)$, the Gram entries are
\[
 G^{(\ell)}_{ij}
       =Q_\partial(\widehat f_{\ell j_i},\widehat f_{\ell j_j}).
\]
The dilation matrix is
$M^{(\ell)}_{ij}=(d_i+d_j)G^{(\ell)}_{ij}$.
The kernel description ensures $G^{(\ell)}_{ii}>0$, which settles
the one-dimensional matrices.
For a two-dimensional matrix, positivity is equivalent to
\begin{equation}\label{eq:shape-correlation-test}
 \frac{(G^{(\ell)}_{01})^2}{G^{(\ell)}_{00}G^{(\ell)}_{11}}
       <\frac{4d_0d_1}{(d_0+d_1)^2}.
\end{equation}

The explicit Gram formula \eqref{eq:euclidean-gram-evaluator} implies,
for $18\le n\le22$ and $\ell=1,2,3$,
\begin{equation}\label{eq:shape-correlation-bounds}
 \frac{(G^{(\ell)}_{01})^2}{G^{(\ell)}_{00}G^{(\ell)}_{11}}
       <\frac34.
\end{equation}
The corresponding thresholds in \eqref{eq:shape-correlation-test}
are $8/9,3/4,8/9$.
At $n=22$, the additional pairs $\ell=0,4$ have degrees $(3,5)$.
Their correlations are less than $4/5$ and $3/4$, respectively,
both below the threshold $15/16$.
The only three-dimensional matrix has $(n,\ell)=(22,2)$ and
degrees $(1,3,5)$.
After normalization by $G^{(2)}_{00}$, its Schur pivots satisfy
$p_1=2$, $p_2>3$ and $p_3>3/2$, again by
\eqref{eq:euclidean-gram-evaluator}.
Every radial matrix is therefore positive definite.
Finite dimensionality and the orthogonal angular decomposition
establish \eqref{eq:euclidean-coercivity}.
\end{proof}

To apply Proposition~\ref{prop:upper-criterion}, retain graph terms
through degree $k+1$ in the quadratic form.
The linear correction uses graph terms through degree $L+1$,
where $L=\lfloor(n-3)/2\rfloor$, because
\eqref{eq:graph-metric-variation} lowers the degree by one.

For the full normalized metric coefficients $H_d$,
\eqref{eq:euclidean-collar} identifies the linear part $h_{f_{d+1}}$, modulo
the variations in \eqref{eq:removable-homogeneous-jet}.
Every additional coefficient is a product of lower graph
coefficients with total metric weight $d$, assigning weight $i$
to $f_{i+1}$.
On the fixed bounded family, the weighted size of these products
is at most $C\eps$ times the weighted size of the independent graph
coefficients.
Consequently \eqref{eq:graph-metric-dilation} implies
\begin{equation*}
 \left|\cQ\left(\sum_{d=1}^k\eps^dH_d\right)
       -\mathcal Q_\partial(\eps f_2,\ldots,\eps^kf_{k+1})\right|
 \le C\eps\sum_{d=1}^k\eps^{2d}\|f_{d+1}\|^2.
\end{equation*}
The same coefficient bounds give
$|H|_\eps^2\le C\sum_{d=1}^k\eps^{2d}\|f_{d+1}\|^2$.
For sufficiently small $\eps$, Proposition~\ref{prop:euclidean-shape-form}
therefore implies
\[
 \cQ\left(\sum_{d=1}^k\eps^dH_d\right)
 \ge(c-C\eps)\sum_{d=1}^k\eps^{2d}\|f_{d+1}\|^2
 \ge c'|H|_\eps^2.
\]
Proposition~\ref{prop:upper-criterion} now establishes
\eqref{eq:upper-main} for $18\le n\le22$, completing the Euclidean
range in Theorem~\ref{thm:classification}.

\subsection{The eighteenth-dimensional endpoint}\label{subsec:eighteen-upper}

For general metrics in dimension eighteen, $k=4$ and $L=7$.
The cubic estimate leaves the fourth jet, whose square has weight
$\eps^8=\eps^{(n-2)/2}$, the order of the outer boundary term in
\eqref{eq:adjoint}.
Fix a uniform conformal Fermi normalization through degree eight
on one original neighborhood of radius $\delta_0$, using
Lemma~\ref{lem:fermi-normalization}.
Write $H_d(z)$ for its independent logarithmic coefficient fields.
With the center argument suppressed, use
\[
 |H|_{\eps,3}^2=\sum_{d=1}^3\eps^{2d}|H_d|^2.
\]
In the expansion of $\cQ(\sum_{d=1}^4\eps^dH_d)$, the mixed
term between degree four and degrees one through three is
\[
 M_\eps=2\sum_{d=1}^3\eps^{d+4}\cQ(H_d,H_4).
\]
Boundedness of the bilinear form on the coefficient spaces yields
\[
 |M_\eps|
 \le C\eps^4|H_4|\sum_{d=1}^3\eps^d|H_d|
 \le C\eps^4|H_4||H|_{\eps,3}.
\]

\begin{proposition}
In dimension eighteen, \eqref{eq:upper-main} holds under any of
the following conditions.
\begin{enumerate}
\item $H_4(0)=0$ or $(H_1,H_2,H_3)(0)\ne0$.
\item $\displaystyle\lim_{x\to0}u(x)=\infty$.
\item The conformal metric $u^{4/(n-2)}g$ is complete at zero.
\end{enumerate}
Every normalized selected sequence $x_k\to0$ violating the upper
bound has a subsequence for which
$(H_1,H_2,H_3)(0)=0$, $H_4(0)\ne0$, and
\[
 \lim_{k\to\infty}
       \frac{|H(x_k)|_{\eps_k,3}^2}{\eps_k^8}\in(0,\infty).
\]
\end{proposition}
\begin{proof}
Suppose the upper bound fails.
Use the normalized selected sequence and the functions
$W,w,Y,I^\pm$ from the proof of
Proposition~\ref{prop:upper-criterion}, with $n=18$,
$R=\delta/\eps$, $L=7$ and $0<\delta\le\delta_0$.
In particular, $H_d=H_d(x_k)$, $w(0)=\nabla_{x'}w(0)=0$ and
$w\to0$ locally.
For fixed $\delta$, \eqref{eq:two-sided-outer-cap} and
\eqref{eq:corrected-bubble-positivity} imply, for late indices,
\[
 c_\delta\eps^8\le w\le C_\delta\eps^8
                                  \quad\text{on }\Sigma_R.
\]
The eighth-order Taylor estimates in
Appendix~\ref{app:taylor-estimates} give
\begin{equation}\label{eq:eighteen-residual-budget}
 \begin{aligned}
 \|\mathcal R\|_G&\le C(|H|_{\eps,3}^2+\eps^8),\\
 \|\mathcal R\|_*^2
 &\le C\eps^2|H|_{\eps,3}^2
                +C\eps^{16}(1+|\log\eps|)^2
 =o_\delta(|H|_{\eps,3}^2+\eps^8).
 \end{aligned}
\end{equation}
By \eqref{eq:adjoint-local-errors} and local convergence,
\[
 \left|-\int_{\mathcal B_R^+}wL_gY\,\dd V_g
       +\int_{B_R}w(B_g-c_W)Y\,\dd\sigma_g\right|
       =o_\delta(|H|_{\eps,3}^2+\eps^8).
\]
The absorption argument following \eqref{eq:core-estimate}, using
the residual bound \eqref{eq:eighteen-residual-budget}, yields
\begin{equation}\label{eq:eighteen-adjoint-upper}
 I^+-I^-\le o_\delta(|H|_{\eps,3}^2+\eps^8)-\tfrac12I^-
 \le o_\delta(|H|_{\eps,3}^2+\eps^8).
\end{equation}

The corrections through degree $L=7$ leave $H_8$ as the first
uncancelled linear term.
Its energy coefficient satisfies $J_1(H_8,U)=0$ by the cutoff
identity in Lemma~\ref{lem:taylor-budget}.
The remaining projection, computed in the eighteen-dimensional part of
Appendix~\ref{app:taylor-estimates}, is
\begin{equation}\label{eq:eighteen-projection}
 \mathcal R[Y]
 =\cQ\left(\sum_{d=1}^3\eps^dH_d\right)
       +\eps^8\cQ(H_4)+M_\eps
       +o_\delta(|H|_{\eps,3}^2+\eps^8).
\end{equation}
By \eqref{eq:cubic-coercivity},
$\cQ(\sum_{d=1}^3\eps^dH_d)\ge c|H|_{\eps,3}^2$.
Since degree $4$ is below $(n-2)/2=8$, Lemma~\ref{lem:ball-hessian} and
\eqref{eq:compactified-coefficient} imply
\[
 \cQ(H_4)=-8\mathcal C(H_4)\ge0.
\]
The lower bound for $w$ on $\Sigma_R$ and \eqref{eq:outer-flux} give
\[
 \int_{\Sigma_R}w\partial_{\nu_g}Y\,\dd\sigma_g
                                  \ge c_\delta\eps^8.
\]
Combining these bounds with \eqref{eq:adjoint} and
\eqref{eq:eighteen-adjoint-upper} yields
\[
 c|H|_{\eps,3}^2+c_\delta\eps^8
 \le C|H_4|\eps^4|H|_{\eps,3}
                  +o_\delta(|H|_{\eps,3}^2+\eps^8).
\]
For fixed $\delta$, absorb the remainder into the left side.
Then $|H|_{\eps,3}>0$ for late indices, and division yields
\begin{equation}\label{eq:eighteen-balance}
 \frac{|H|_{\eps,3}}{\eps^4}
       +\frac{\eps^4}{|H|_{\eps,3}}
 \le C_\delta|H_4|.
\end{equation}
The right side is bounded, so a subsequence has the stated
positive finite limit.
Moreover,
\[
 |H_d(x_k)|\le C_\delta\eps_k^{4-d}
                \quad\text{for }d=1,2,3.
\]
By continuity, $(H_1,H_2,H_3)(0)=0$.
The left side of \eqref{eq:eighteen-balance} is at least two,
so $H_4(0)\ne0$. Thus condition~(i) implies the upper bound.

Suppose now that $u(x)\to\infty$ uniformly.
Keep the choices in Lemma~\ref{lem:fermi-normalization} fixed while
varying $\delta$.
Young's inequality in \eqref{eq:eighteen-projection} yields
\begin{equation}\label{eq:eighteen-loss}
 \begin{aligned}
 \mathcal R[Y]&-\int_{\mathcal B_R^+}wL_gY\,\dd V_g
                    +\int_{B_R}w(B_g-c_W)Y\,\dd\sigma_g\\
 &\ge c_1|H|_{\eps,3}^2-C_1\eps^8
                    -o_\delta(|H|_{\eps,3}^2+\eps^8).
 \end{aligned}
\end{equation}
The constants $c_1,C_1$ depend on the fixed coefficient bounds on
the original neighborhood, uniformly for $0<\delta\le\delta_0$.
Let $c_0>0$ be the uniform lower flux bound in
\eqref{eq:uniform-outer-flux}, and choose $M>(C_1+1)/c_0$.
Uniform divergence and the uniform center gauges allow a fixed
$\delta>0$ such that $v\ge M\eps^8$ on $\Sigma_R$ for late indices.
Since $W=O_\delta(\eps^{16})$ there,
\[
 \int_{\Sigma_R}w\partial_{\nu_g}Y\,\dd\sigma_g
              \ge c_0M\eps^8-O_\delta(\eps^{16}).
\]
Equations \eqref{eq:adjoint}, \eqref{eq:eighteen-adjoint-upper}
and \eqref{eq:eighteen-loss} imply
\[
 c_1|H|_{\eps,3}^2+(c_0M-C_1)\eps^8
       \le o_\delta(|H|_{\eps,3}^2+\eps^8),
\]
which is a contradiction for this fixed $\delta$.

If $u^{4/(n-2)}g$ is complete at zero, use the positive gauge constructed in
Appendix~\ref{app:green-kernels}.
The conformal metric satisfies
\[
 (\varphi^{-1}u)^{4/(n-2)}
                  \bigl(\varphi^{4/(n-2)}g\bigr)
       =u^{4/(n-2)}g.
\]
Proposition~\ref{prop:potential-completeness} implies
$\varphi^{-1}u\to\infty$, hence $u\to\infty$.
The preceding case applies.
This proves Theorem~\ref{thm:classification}(ii) in dimension eighteen.
\end{proof}

\section{The matching lower bound and removability}\label{sec:lower}

Like in the lower dimension cases, we establish the following two-sided estimates for positive solutions of
\eqref{eq:main} near an isolated boundary singularity.
\begin{proposition}\label{prop:matching-lower}
    Let $n\ge3$ and let $u$ be a positive solution of \eqref{eq:main}
    satisfying \eqref{eq:upper-main}. Then $u$ either extends smoothly
    across zero or satisfies the lower bound in \eqref{eq:two-sided}.
\end{proposition}
Under the critical upper bound \eqref{eq:upper-main}, the solution either
extends smoothly across the singularity or has two-sided critical
growth. The proof follows the Pohozaev and discrete-scale argument of
\cite[Sections 4--9]{HJWXZ1}, with the annular comparison estimates
adapted to the boundary equation.

\subsection{The Pohozaev limit and neck geometry}

Work in fixed smooth boundary Fermi coordinates at zero.
Let $S=\Sph^{n-1}_+$, with outward normal $\nu_S$, and set
$r=e^{-t}$ and
\[
 V(t,\omega)=r^{\frac{n-2}{2}}u(r\omega).
\]
Its spherical mean is $\bar V(t)=|S|^{-1}\int_SV(t,\omega)\dd\omega$.
Thus \eqref{eq:upper-main} becomes $0<V\le C$ for large $t$.
Write the equation as
\begin{equation}\label{eq:forced-cylinder}
 \begin{cases}
 V_{tt}+\Delta_SV-\frac{(n-2)^2}{4}V=\mathcal E
 &\text{in }(0,\infty)\times S,\\
 \partial_{\nu_S}V=(n-2)V^{\frac n{n-2}}+\mathcal E^{\partial}
 &\text{on }(0,\infty)\times\partial S.
 \end{cases}
\end{equation}
The terms $\mathcal E,\mathcal E^\partial$ arise from
$g-I=O(r)$ and its derivatives.

To track the neck limits, use the flat Pohozaev functional
\begin{equation}\label{eq:pohozaev}
 P(r,f)=\int_{\Sigma_r}
 \left({\frac{n-2}{2}}f f_r-\frac r2|\nabla f|^2+rf_r^2\right)
       +{\frac{(n-2)^2r}{2(n-1)}}\int_{\partial B_r}|f|^{\frac{2(n-1)}{n-2}}.
\end{equation}
When the limit exists, write
$\mathcal P(u)=\lim_{r\downarrow0}P(r,u)$.

\begin{lemma}\label{lem:cylindrical-estimates}
For all sufficiently large $t$ and $\omega\in\overline S$,
\begin{enumerate}
\item
\begin{equation}\label{eq:relative-harnack}
 C^{-1}\bar V(t)\le V(t,\omega)\le C\bar V(t).
\end{equation}
For each $j\ge1$, all order-$j$ derivatives in $(t,\omega)$ satisfy
$|D^jV(t,\omega)|\le C_j\bar V(t)$.
\item The right-hand sides of \eqref{eq:forced-cylinder} satisfy,
for every fixed $k\ge0$,
\[
 \|\mathcal E\|_{C^k([t-1,t+1]\times\overline S)}
 +\|\mathcal E^{\partial}\|_{C^k([t-1,t+1]\times\partial S)}
 \le C_ke^{-t}{\bar V}(t).
\]
The spherical mean obeys
\[
 |\bar V''-\tfrac{(n-2)^2}{4}\bar V|
       \le C(\bar V^{\frac2{n-2}}+e^{-t})\bar V.
\]
\item The Pohozaev functional satisfies
\begin{equation}\label{eq:relative-pohozaev-drift}
 \left|\frac d{dt}P(e^{-t},u)\right|\le Ce^{-t}\bar V(t)^2.
\end{equation}
In particular, $\mathcal P(u)$ exists and is finite.
\end{enumerate}
The constants depend only on $n$, the fixed metric, the upper-bound
constant and the indicated derivative order.
\end{lemma}
\begin{proof}
\emph{Step 1. Proof of (i).}
For $r=e^{-t}$, put $u_r(y):=r^{(n-2)/2}u(ry)$ and $g_r(y):=g(ry)$. On $\{1/4<|y|<4,\,y_n>0\}$,
\[
 \begin{cases}
 -L_{g_r}u_r=0&\text{for }y_n>0,\\
 \partial_{\nu_{g_r}}u_r
 =\left((n-2)u_r^{\frac2{n-2}}-\frac{n-2}{2}r\,h_g(ry)\right)u_r
     &\text{for }y_n=0.
 \end{cases}
\]
Since $u_r\le C$, the Robin coefficient is uniformly bounded. We have
\[
 u_r(\omega)=V(t,\omega)\quad\text{on }S,\qquad
 \bar V(t)=\frac1{|S|}\int_S u_r\dd\omega.
\]
Thus \eqref{eq:buffered-robin-harnack} and
\eqref{eq:buffered-robin-estimate}, applied to $u_r/\bar V(t)$, imply
\[
 C^{-1}\le\frac{u_r}{\bar V(t)}\le C,\qquad
 \left\|\frac{u_r}{\bar V(t)}\right\|_{C^k}\le C_k
\]
on fixed smaller half-annuli. The relation
\[
 u_r(e^{-\tau}\omega)=e^{(n-2)\tau/2}V(t+\tau,\omega)
 \quad\text{for }|\tau|\le1
\]
then yields \eqref{eq:relative-harnack} and $|D^j V|\le C_j\bar V(t)$.

\emph{Step 2. Proof of (ii).}
By direct coordinate transformation,
\[
 \begin{aligned}
 \mathcal E(t,\omega)
 &=\left[(\Delta-\Delta_{g_r})u_r
       +\frac{n-2}{4(n-1)}R_{g_r}u_r\right](\omega),\\
 \mathcal E^\partial(t,\omega)
 &=-\frac{n-2}{2}r\,h_g(r\omega)V(t,\omega).
 \end{aligned}
\]
The rescaled coefficients satisfy
\[
 \|g_r-I\|_{C^{k+1}}\le C_k r,\qquad
 R_{g_r}(y)=r^2R_g(ry).
\]
Together with the bounds in Step 1, these imply
\[
 \|\mathcal E\|_{C^k([t-1,t+1]\times\overline S)}
 +\|\mathcal E^\partial\|_{C^k([t-1,t+1]\times\partial S)}
 \le C_k e^{-t}\bar V(t).
\]
Integrating \eqref{eq:forced-cylinder} over $S$ yields
\[
 |S|\left(\bar V''-\frac{(n-2)^2}{4}\bar V\right)
 = \int_S\mathcal E - (n-2)\int_{\partial S}V^{\frac{n}{n-2}} - \int_{\partial S}\mathcal E^\partial,
\]
whence \eqref{eq:relative-harnack} and the source bounds prove the mean inequality.

\emph{Step 3. Proof of (iii).}
In cylindrical variables,
\[
P(e^{-t},u) = \frac12\int_S\left(V_t^2-|\nabla_S V|^2-\frac{(n-2)^2}{4}V^2\right) +
\frac{(n-2)^2}{2(n-1)}\int_{\partial S}V^{\frac{2(n-1)}{n-2}}.
\]
Differentiating in $t$, integrating by parts on $S$, and using \eqref{eq:forced-cylinder} yields
\[
\begin{aligned}
 \frac{\dd}{\dd t}P(e^{-t},u)
&= \int_S V_t\left(V_{tt}+\Delta_S V-\frac{(n-2)^2}{4}V\right) - \int_{\partial S}
V_t\left(\partial_{\nu_S}V-(n-2)V^{\frac{n}{n-2}}\right)\\
 &= \int_S V_t\mathcal E - \int_{\partial S}V_t\mathcal E^\partial.
\end{aligned}
\]
Using (i)--(ii) proves \eqref{eq:relative-pohozaev-drift}. Since $\bar V(t)\le C$,
$\int_{t_0}^\infty e^{-t}\bar V(t)^2\dd t<\infty$, so $\mathcal P(u)\in\R$ exists.
\end{proof}

Suppose that zero is nonremovable and the lower bound in
\eqref{eq:two-sided} fails. Then there are $s_j\to\infty$ with
$\min_{\omega\in\overline S}V(s_j,\omega)\to0$.
By \eqref{eq:relative-harnack}, $\bar V(s_j)\to0$.
The cylindrical expression for $P$ and the relative derivative
estimates imply
\[
 |P(e^{-s_j},u)|\le C\bar V(s_j)^2
       +C\bar V(s_j)^{2(n-1)/(n-2)}\longrightarrow0.
\]
Hence $\mathcal P(u)=0$.

If $\bar V(t)\to0$, the mean estimate implies, for every
$0<\eta<(n-2)/2$ and all $t\ge T_\eta$,
\[
 \bar V''\ge\left(\frac{n-2}{2}-\eta\right)^2\bar V.
\]
Comparison on $[T_\eta,T]$, followed by $T\to\infty$, gives
\[
 \bar V(t)\le\bar V(T_\eta)
       e^{-(\frac{n-2}{2}-\eta)(t-T_\eta)}.
\]
Thus $u(x)\le C_\eta|x|^{-\eta}$ by relative Harnack.
Fix $p>n-1$ and choose $\eta$ small enough that
$pn\eta/(n-2)<n-1$. Then
\[
 \int_{B_{r_0}}u^{pn/(n-2)}\dd\sigma_g
 \le C_\eta\int_0^{r_0}r^{n-2-pn\eta/(n-2)}\dd r<\infty.
\]
The Green representation \eqref{eq:green-representation} and its
kernel bound imply
\[
 \|u\|_{L^\infty(\mathcal B_{r_0/2}^+)}
 \le C\left(1+\|u^{n/(n-2)}\|_{L^p(B_{r_0})}\right)<\infty.
\]
The mixed inverse in Appendix~\ref{app:green-kernels} provides the
$H^1$ weak extension for this boundary datum.
The local boundary estimates in Appendix~\ref{app:local-robin}
then establish a continuous extension and $C^1$ regularity.
The maximum and boundary-point principles ensure $u(0)>0$.
Elliptic regularity now provides a smooth extension, contradicting
nonremovability. Consequently,
\begin{equation}\label{eq:neck-alternative}
 \liminf_{t\to\infty}{\bar V}(t)=0<\limsup_{t\to\infty}{\bar V}(t),
 \qquad \mathcal P(u)=0.
\end{equation}

To locate the peaks, consider $V(s_j+\cdot,\cdot)$ for
$s_j\to\infty$. The relative estimates ensure compactness on every
fixed closed cylinder. Each
subsequential limit $v\ge0$ satisfies
\begin{equation}\label{eq:flat-cylinder}
 \begin{cases}
 v_{tt}+\Delta_Sv-\frac{(n-2)^2}{4}v=0&\text{in }\R\times S,\\
 \partial_{\nu_S}v=(n-2)v^{\frac n{n-2}}&\text{on }\R\times\partial S.
 \end{cases}
\end{equation}
If $v\ne0$, the strong maximum and boundary-point principles ensure
$v>0$ on $\R\times\overline S$.
The cylindrical expression for $P$ passes to the limit, with
$P(e^{-(s_j+t)},u)\to\mathcal P(u)=0$ at every fixed $t$.
By \cite[Theorem 2.1]{JLXBoundary}, the inverse transform is either
a regular bubble or has a nonremovable singularity and is axially
symmetric. The latter case has strictly negative Pohozaev invariant
by \cite[Theorem 1.1]{JLXBoundary}, so every nonzero limit is the
cylindrical transform of a regular bubble.

For $\theta=(\lambda,\xi)$, write
$v_\theta(t,\omega)=r^{(n-2)/2}U_{\lambda,\xi}(r\omega)$.
Direct differentiation yields
\begin{equation}\label{eq:bubble-fixed-origin-dilation}
 \begin{aligned}
 \partial_tv_\theta
 &=-r^{(n-2)/2}
       \left(\frac{n-2}{2}U_{\lambda,\xi}
                   +y\cdot\nabla U_{\lambda,\xi}\right)_{y=r\omega}\\
 &=-\frac{n-2}{2}
       \frac{\lambda^2+|\xi|^2-r^2}
            {|r\omega'-\xi|^2+(r\omega_n+\lambda)^2}\,v_\theta.
 \end{aligned}
\end{equation}
Its sign is independent of $\omega$, so the only critical time
of $\bar v_\theta$ occurs at $r^2=\lambda^2+|\xi|^2$.
At that time,
\[
 \partial_{tt}v_\theta
 =-\frac{(n-2)r^2v_\theta}
          {|r\omega'-\xi|^2+(r\omega_n+\lambda)^2}<0.
\]
Thus the mean has a unique nondegenerate maximum.

For uniform comparison, translate each bubble profile so that its mean maximum occurs
at $t=0$, equivalently $\lambda^2+|\xi|^2=1$.
Let $C_0$ be the cylindrical upper bound inherited by all limits.
Evaluating this bound on the boundary gives
\[
 (|\xi|/\lambda)^{(n-2)/2}
       =|\xi|^{(n-2)/2}U_{\lambda,\xi}(\xi,0)\le C_0.
\]
Hence $1=\lambda^2+|\xi|^2\le(1+C_0^{4/(n-2)})\lambda^2$.
The normalized parameters lie in a compact set $\Theta$ with
$\lambda\ge c(C_0)>0$.
Their mean maxima therefore satisfy
\[
 c_0=\min_{\theta\in\Theta}\bar v_\theta(0)>0.
\]
The explicit formula also implies
$|D^kv_\theta(t,\omega)|\le C_ke^{-(n-2)|t|/2}$ uniformly on $\Theta$.

Choose $0<\eta_0<\min(c_0/4,\limsup\bar V/2)$ sufficiently small.
For large $t$, the mean estimate implies
\[
 \bar V''\ge\frac{(n-2)^2}{8}\bar V>0
\]
whenever $\bar V\le2\eta_0$.
By \eqref{eq:neck-alternative}, the late components
$(t_j^-,t_j^+)$ of $\{\bar V<\eta_0\}$ have finite endpoints and
occur infinitely often. Enumerate them in increasing order.
Their endpoints satisfy $\bar V(t_j^\pm)=\eta_0$, and strict convexity
ensures a unique minimum $b_j$, with $\delta_j=\bar V(b_j)$.
If $\delta_j\ge\delta>0$ along a subsequence, a limit of
$V(b_j+\cdot,\cdot)$ is nonzero and hence is a bubble.
Passing the minimum condition and the convexity estimate to the limit gives
\[
 \bar v_\infty'(0)=0,\qquad
 \bar v_\infty''(0)\ge\frac{(n-2)^2}{8}\bar v_\infty(0)>0,
\]
contrary to the strict negativity at the sole critical time of a
bubble mean. Therefore $\delta_j\to0$.

For the exact annular increments in
Proposition~\ref{prop:annular-increment}, use the metric Pohozaev
functional at the adjacent minima. Write $\gamma$ for the induced face metric,
take $\nu=x/r$ on $\Sigma_r$, and use Euclidean measures on the
sphere and its edge. Define
\begin{equation}\label{eq:metric-pohozaev}
\begin{aligned}
 P_g(r,f)&={}\int_{\Sigma_r}\left\{
  (x\cdot\nabla f+\tfrac{n-2}{2}f)g^{ij}f_i\nu_j
  -\frac r2\left(|\nabla f|_g^2+
          \frac{n-2}{4(n-1)}R_gf^2\right)\right\}\sqrt{\det g}\\
 &+r\int_{\partial B_r}\left\{
  \frac{(n-2)^2}{2(n-1)}|f|^{\frac{2(n-1)}{n-2}}
                -\frac{n-2}{4}h_gf^2\right\}\sqrt{\det\gamma}.
\end{aligned}
\end{equation}
By Lemma~\ref{lem:cylindrical-estimates} and $g=I+O(r)$,
\[
 |P_g(r,u)-P(r,u)|\le Cr\bar V(-\log r)^2.
\]
Thus $P_g(r,u)\to\mathcal P(u)$ as $r\to0$.

Choose $a_j$ to maximize $\bar V$ on $[t_{j-1}^+,t_j^-]$ and set
\[
 \bar\rho_j=e^{-b_j}<\rho_j=e^{-a_j}<\bar\rho_{j-1}.
\]
In the rescaling
\[
 v_j(y)=\rho_j^{\frac{n-2}{2}}u(\rho_jy),
\]
the annulus between the neighboring minima has radii
\[
 \sigma_j=\frac{\bar\rho_j}{\rho_j}=e^{a_j-b_j},\qquad
 R_j=\frac{\bar\rho_{j-1}}{\rho_j}=e^{a_j-b_{j-1}}.
\]
Write $\mu_j=-P_g(\bar\rho_j,u)$ for the negative of the
Pohozaev value at the inner minimum.

\begin{lemma}\label{lem:exact-neck-shape}
For the sequence constructed above, the following hold for all
sufficiently large $j$, with constants independent of $j$.
\begin{enumerate}
\item The maximizing time $a_j$ is unique and nondegenerate.
\item For $y\in\overline{\R^n_+}$ with
$\sigma_j/2\le r=|y|\le2R_j$,
\begin{equation}\label{eq:annular-bounds}
 \begin{gathered}
 v_j(y)\asymp(1+r)^{2-n},\\
 r|\nabla v_j(y)|+r^2|D^2v_j(y)|\le C(1+r)^{2-n}.
 \end{gathered}
\end{equation}
\item The quantities $\mu_j=-P_g(\bar\rho_j,u)$ are positive and satisfy
\begin{equation}\label{eq:capacity-matching}
 \begin{aligned}
 \mu_j\asymp\sigma_j^{n-2},\qquad
 \mu_{j-1}\asymp R_j^{2-n}.
 \end{aligned}
\end{equation}
Moreover, $\mu_j\asymp(\rho_{j+1}/\rho_j)^{(n-2)/2}\to0$.
\item The outer radii satisfy
\begin{equation}\label{eq:neck-scale-summability}
 \sum_j\rho_jR_j<\infty.
\end{equation}
In particular, $\rho_jR_j\to0$.
\end{enumerate}
\end{lemma}
\begin{proof}
\emph{Step 1. Proof of (i).}
Set $I_j:=[t_{j-1}^+,t_j^-]$. On this interval,
\[
 \bar V\ge\eta_0,\qquad
 \bar V(t_{j-1}^+)=\bar V(t_j^-)=\eta_0.
\]
If $|I_j|\to\infty$ along a subsequence, the classification
following \eqref{eq:neck-alternative} gives a limit $v_\infty\ge0$
of $V(t_{j-1}^++\cdot,\cdot)$. Since $\bar v_\infty(s)\ge\eta_0$
for all $s\ge0$, this limit is a bubble $v_\theta$.
Equation~\eqref{eq:bubble-fixed-origin-dilation} implies
$\bar v_\theta(s)\to0$ as $s\to\infty$, a contradiction.
Consequently, $|I_j|\le C$.

If $\max_{I_j}\bar V\le2\eta_0$, then
\[
 \bar V''\ge\frac{(n-2)^2}{8}\bar V>0\quad\text{on }I_j.
\]
Strict convexity and the endpoint values force $\bar V<\eta_0$
on $\operatorname{int}I_j$, contradicting $\bar V\ge\eta_0$.
Hence $\max_{I_j}\bar V>2\eta_0$, and $a_j\in\operatorname{int}I_j$.
Any subsequential limit $v_\infty$ of $V(a_j+\cdot,\cdot)$ is a
bubble $v_\theta$ satisfying
\[
 \bar v_\infty(0)\ge c_0>4\eta_0,\qquad
 \bar v_\infty'(0)=0,\qquad \bar v_\infty''(0)<0.
\]
By compactness of $\Theta$ and $C^2$ convergence on $I_j$, there
exist fixed $s_0,c>0$ such that for all large $j$,
\[
 \begin{aligned}
 \bar V''(t)&\le-c\quad\text{for }|t-a_j|\le s_0,\\
 |\bar V'(t)|&\ge c\quad\text{for }t\in I_j,\ |t-a_j|\ge s_0,
 \end{aligned}
\]
proving (i) and $|a_j-t_{j-1}^+|+|t_j^--a_j|\le C$.

\emph{Step 2. Two-exponential comparison and proof of (ii).}
On $[t_j^-,t_j^+]$, $\bar V''\ge\frac{(n-2)^2}{8}\bar V$ and $\bar V(t_j^\pm)=\eta_0$ imply
\[
 \bar V(t) \le \eta_0\left(e^{-\frac{n-2}{2\sqrt{2}}(t-t_j^-)}+e^{-\frac{n-2}{2\sqrt{2}}(t_j^+-t)}\right).
\]
Writing $\bar V''-\frac{(n-2)^2}{4}\bar V=e\bar V$, Lemma~\ref{lem:cylindrical-estimates}(ii)
implies $\|e\|_{L^1(t_j^-,t_j^+)}\le C\int_{t_j^-}^{t_j^+}(\bar V^{2/(n-2)}+e^{-t})\dd t\le
C\eta_0^{2/(n-2)}+Ce^{-t_j^-}$. Comparing $\bar V$ with the solution $h_j''=\frac{(n-2)^2}{4}h_j$,
$h_j(t_j^\pm)=\eta_0$, namely
\[
 h_j(t) := \eta_0\frac{\cosh(\frac{n-2}{2}(t-(t_j^-+t_j^+)/2))}{\cosh(\frac{n-2}{4}(t_j^+-t_j^-))},
\]
via the Dirichlet Green function on $[t_j^-,t_j^+]$ yields
\[
 |\bar V(t)-h_j(t)|+|\bar V'(t)-h_j'(t)| \le C\int_{t_j^-}^{t_j^+}e^{-\frac{n-2}{2}|t-s|}|e(s)|\bar V(s)\dd s.
\]
Using $h_j(s)\le e^{\frac{n-2}{2}|t-s|}h_j(t)$ and choosing $\eta_0$ small and $j$ large so that
$C\|e\|_{L^1}<1/2$ yields $|\bar V-h_j|+|\bar V'-h_j'|\le C(\eta_0^{2/(n-2)}+e^{-t_j^-})h_j$. In
particular,
\begin{equation}\label{eq:neck-two-exponentials}
 \bar V(t)\asymp\eta_0
 \left(e^{-\frac{n-2}{2}(t-t_j^-)}+e^{-\frac{n-2}{2}(t_j^+-t)}\right)
 \quad\text{for }t_j^-\le t\le t_j^+.
\end{equation}
Since $\bar V'(b_j)=0$, the derivative comparison gives
\[
 \frac{n-2}{2}\left|\tanh\left(\frac{n-2}{2}
       \left(b_j-\frac{t_j^-+t_j^+}{2}\right)\right)\right|
 =\frac{|h_j'(b_j)|}{h_j(b_j)}
 \le C(\eta_0^{2/(n-2)}+e^{-t_j^-}).
\]
Hence $b_j=(t_j^-+t_j^+)/2+O(1)=(a_j+a_{j+1})/2+O(1)$.
Evaluating \eqref{eq:neck-two-exponentials} at $b_j$ yields
$\delta_j^2\asymp e^{-(n-2)(a_{j+1}-a_j)/2}$. Thus
\begin{equation}\label{eq:neck-minimum-scales}
 \delta_j^2 \asymp e^{-(n-2)(b_j-a_j)} = \sigma_j^{n-2},\qquad
 \delta_{j-1}^2 \asymp e^{-(n-2)(a_j-b_{j-1})} = R_j^{2-n}.
\end{equation}
Since $|I_j|\le C$, \eqref{eq:neck-two-exponentials} yields $\bar V(t)\asymp e^{(n-2)(t-a_j)/2}$
for $b_{j-1}\le t\le a_j$ and $\bar V(t)\asymp e^{-(n-2)(t-a_j)/2}$ for $a_j\le t\le b_j$. By
Lemma~\ref{lem:cylindrical-estimates}(i) and $v_j(r\omega)=r^{-(n-2)/2}V(a_j-\log r,\omega)$, we
obtain \eqref{eq:annular-bounds} on $\sigma_j\le r\le R_j$.

To extend the estimates across $r=\sigma_j,R_j$, apply
Lemma~\ref{lem:cylindrical-estimates} to
$\delta_j^{-1}V(b_j+s,\omega)$. This sequence is bounded in
$C^2_{\mathrm{loc}}(\R\times\overline S)$. Its geometric errors
are $O(e^{-b_j})$, and its nonlinear boundary coefficient is
$\delta_j^{2/(n-2)}$, both tending to zero.
Every subsequential limit $w>0$ therefore solves
\[
 \begin{cases}
 w_{ss}+\Delta_S w-\dfrac{(n-2)^2}{4}w=0
     &\text{in }\R\times S,\\
 \partial_{\nu_S}w=0&\text{on }\R\times\partial S.
 \end{cases}
\]
Even reflection across $\partial S$ and inversion of the cylindrical
transform yield a positive harmonic function on $\R^n\setminus\{0\}$.
By \cite[Corollary 3.14]{ABR}, this function has the form
\[
 c_1|x|^{2-n}+c_2,\qquad c_1,c_2\ge0.
\]
The normalization at $s=0$ implies
\[
 \begin{aligned}
 w(0)&=c_1+c_2=1,\\
 w_s(0)&=\frac{n-2}{2}(c_1-c_2)=0.
 \end{aligned}
\]
Thus $c_1=c_2=1/2$, and the full sequence satisfies
\begin{equation}\label{eq:neck-minimum-limit}
\delta_j^{-1}V(b_j+s,\omega)\longrightarrow \cosh\left(\frac{n-2}{2}s\right)\qquad\text{in
}C^2_{\mathrm{loc}}(\R\times\overline S).
\end{equation}
Applying \eqref{eq:neck-minimum-limit} on $|s|\le\log 4$ at $b_{j-1}$ and $b_j$ extends
\eqref{eq:annular-bounds} to $\sigma_j/2\le r\le 2R_j$.

\emph{Step 3. Proof of (iii) and (iv).}
At $t=b_j$, \eqref{eq:neck-minimum-limit} implies $V/\delta_j\to 1$, $V_t/\delta_j\to 0$, and $\nabla_S V/\delta_j\to 0$, so
\[
 \frac{P(e^{-b_j},u)}{\delta_j^2} = -\frac{(n-2)^2}{8}|S|+o(1)+O(\delta_j^{\frac{2}{n-2}}).
\]
By the comparison following \eqref{eq:metric-pohozaev},
$|P_g(e^{-b_j},u)-P(e^{-b_j},u)|\delta_j^{-2}\le Ce^{-b_j}\to 0$, whence $\mu_j =
\frac{(n-2)^2}{8}|S|\delta_j^2(1+o(1))>0$. By \eqref{eq:neck-minimum-scales},
$\mu_j\asymp\delta_j^2\asymp\sigma_j^{n-2}\asymp(\rho_{j+1}/\rho_j)^{(n-2)/2}\to 0$ and
$\mu_{j-1}\asymp\delta_{j-1}^2\asymp R_j^{2-n}$, proving \eqref{eq:capacity-matching}. Finally,
since $\rho_{j+1}/\rho_j\to 0$, choosing $J$ with $\rho_{j+1}\le\rho_j/2$ for $j\ge J$ and using
$b_{j-1}>a_{j-1}$ yields
\[
\sum_{j\ge J+1}\rho_j R_j = \sum_{j\ge J+1}\bar\rho_{j-1} \le \sum_{j\ge J+1}\rho_{j-1} \le
\rho_J\sum_{m=0}^\infty 2^{-m} = 2\rho_J < \infty,
\]
proving \eqref{eq:neck-scale-summability}.
\end{proof}

\subsection{The geometric increment and the lower bound}

Under \eqref{eq:neck-alternative}, the positive quantities
$\mu_j\to0$ from Lemma~\ref{lem:exact-neck-shape} are incompatible
with positive increments up to summable relative errors.
To compute the leading increment, we remove the homogeneous
metric terms in \eqref{eq:removable-homogeneous-jet}, whose
reduced coefficient vanishes.

\begin{lemma}\label{lem:finite-jet-removal}
Let $n\ge3$ and $\ell=\lceil(n-2)/2\rceil$.
After finitely many smooth positive conformal changes and
boundary-preserving coordinate changes, the metric in boundary
Fermi coordinates has one of the following expansions.
\begin{enumerate}
\item For an integer $1\le d<(n-2)/2$,
\[
 g=\delta+H_d+O(r^{d+1}),
\]
where $H_d$ is a degree-$d$ representative of the first nonzero
class modulo \eqref{eq:removable-homogeneous-jet}. In particular,
$\log g=H_d+O(r^{d+1})$.
For every compact set $\Theta\subset(0,\infty)\times\R^{n-1}$
of bubble parameters,
\begin{equation}\label{eq:first-quotient}
 -\mathcal C(H_d;\lambda,\xi)\ge c_\Theta|H_d|^2>0
                      \quad\text{for }(\lambda,\xi)\in\Theta.
\end{equation}
\item The metric satisfies
\[
 g=\delta+O(r^\ell).
\]
For $n=3,4$, this alternative holds with $\ell=1$.
\end{enumerate}
Each displayed remainder $F=O(r^q)$ satisfies
$|D^jF|\le C_jr^{q-j}$ for every fixed $j\ge0$.
\end{lemma}
\begin{proof}
For $n=3,4$, $\ell=1$ and Fermi coordinates imply (ii).
For $n\ge5$, assume inductively that degrees below $d$ have been
removed, where $1\le d<\ell$, and write
\[
 g=\delta+H_d+O(r^{d+1}).
\]
By Lemma~\ref{lem:homogeneous-kernel}, $\mathcal C(H_d;\theta)\le0$
for all $\theta=(\lambda,\xi)$. Equality at one parameter value,
and hence at every parameter value, holds if and only if
\[
 H_d=\mathcal L_{Y_{d+1}}\delta+f_d\delta,\qquad
 (Y_{d+1})_n=0\quad\text{on }\Gamma.
\]
If equality holds, let $\Phi_d$ be the time-one flow of $-Y_{d+1}$.
Since $2d\ge d+1$, the coordinate and conformal changes satisfy
\[
 \begin{aligned}
 \Phi_d(x)&=x-Y_{d+1}(x)+O(r^{2d+1}),\\
 \Phi_d^*g&=\delta+f_d\delta+O(r^{d+1}),\\
 \widehat g&:=e^{-f_d}\Phi_d^*g=\delta+O(r^{d+1}).
 \end{aligned}
\]
The Christoffel symbols of $\widehat g$ are $O(r^d)$, and its
inward normal on $\Gamma$ is $e_n+O(r^{d+1})$.
Integrating the boundary and normal geodesic equations, we obtain
a Fermi map satisfying
\[
 \Psi_d(x)=x+O(r^{d+2}),\qquad
 \Psi_d^*\widehat g=\delta+O(r^{d+1}).
\]

Iterate for $d=1,\dots,\ell-1$. If a first term with
$\mathcal C(H_d;\theta)<0$ occurs, continuity and compactness of
$\Theta$ imply \eqref{eq:first-quotient}, proving alternative (i).
Otherwise all degrees below $\ell$ are removed, proving
alternative (ii). At each step,
\[
 \begin{gathered}
 \Phi_d(0)=\Psi_d(0)=0,\qquad D\Phi_d(0)=D\Psi_d(0)=I,\\
 e^{(n-2)f_d(0)/4}=1.
 \end{gathered}
\]
Thus critical bounds and removability are preserved by
\eqref{eq:covariance}.
\end{proof}

Apply Lemmas~\ref{lem:cylindrical-estimates} and
\ref{lem:exact-neck-shape} to the transformed metric and solution.
All neck quantities and Pohozaev values below refer to this pair.

In alternative (i), put $g_\rho(y)=g(\rho y)$.
On the $j$th annulus, with $r=|y|$, the derivative bounds imply
\[
 |g_{\rho_j}-\delta|+r|Dg_{\rho_j}|+r^2|D^2g_{\rho_j}|
 \le C(\rho_jr)^d\le C(\rho_jR_j)^d.
\]
We write
\[
 \eta_j=(\rho_jR_j)^d=\bar\rho_{j-1}^{\,d}
\]
for this geometric error factor in the annular comparison.

\begin{proposition}\label{prop:annular-increment}
In alternative (i) of Lemma~\ref{lem:finite-jet-removal}, for all
sufficiently large $j$,
\begin{equation}\label{eq:increment}
 \mu_j-\mu_{j-1}\ge\tfrac c2\rho_j^{2d}
                -C\eta_j(\mu_j+\mu_{j-1}).
\end{equation}
Moreover, $\sum_j\eta_j<\infty$, and
\begin{equation}\label{eq:increment-recursion}
 \mu_{j-1}\le(1+C'\eta_j)\mu_j
\end{equation}
for all sufficiently large $j$.
\end{proposition}

\begin{proof}[Proof of Proposition~\ref{prop:annular-increment}]
\emph{Step 1. Annular Pohozaev identity.}
For $0<\sigma<R$, set
\[
 \mathcal A_{\sigma,R}^+:=\{x_n>0,\ \sigma<|x|<R\},\qquad
 \Gamma_{\sigma,R}:=B_R\setminus\overline{B_\sigma}.
\]
For a positive $f\in C^\infty(\overline{\mathcal A_{\sigma,R}^+})$,
let $J_{g_\rho,\mathcal A_{\sigma,R}^+}(f)$ be the energy
\eqref{eq:action} restricted to this annulus and its flat face.
Define
\[
 \mathcal V_\rho(f):=-\rho\partial_\rho
                  J_{g_\rho,\mathcal A_{\sigma,R}^+}(f),
\]
where $f,\sigma,R$ are held fixed. Integrating by parts in
\eqref{eq:metric-pohozaev}, we obtain
\begin{align*}
 &P_{g_\rho}(R,f)-P_{g_\rho}(\sigma,f)\\
 &\quad=\mathcal V_\rho(f)
 +\int_{\mathcal A_{\sigma,R}^+}(x\cdot\nabla f+\tfrac{n-2}{2}f)L_{g_\rho}f\dd V_{g_\rho}\\
 &\qquad-\int_{\Gamma_{\sigma,R}}(x\cdot\nabla f+\tfrac{n-2}{2}f)\bigl(B_{g_\rho}f-(n-2)f^{\frac{n}{n-2}}\bigr)\dd\sigma_{g_\rho}.
\end{align*}
For $v_j(y)=\rho_j^{(n-2)/2}u(\rho_j y)$, both residuals vanish.
The scaling identity
\[
 P_{g_{\rho_j}}(r,v_j)=P_g(\rho_j r,u)
\]
therefore implies
\begin{equation}\label{eq:increment-dilation}
 \mathcal V_{\rho_j}(v_j)
 =P_g(\bar\rho_{j-1},u)-P_g(\bar\rho_j,u)
 =\mu_j-\mu_{j-1}.
\end{equation}
Hereafter we suppress $j$ and write
\[
 \sigma=\sigma_j\to0,\qquad R=R_j\to\infty,\qquad
 \rho R=\rho_jR_j\to0.
\]

\emph{Step 2. Corrected bubble and annular error.}
After passing to a subsequence, the compactness and classification
argument following \eqref{eq:neck-alternative} yields
$v\to U_{\theta_*}$ locally smoothly away from $0$.
Fix a compact parameter set $\Theta\subset(0,\infty)\times\R^{n-1}$
containing all possible $\theta_*$ in its interior.
For $\theta=(\lambda,\xi)\in\Theta$, let $\psi_\theta$ be the correction
for $H_d$ at $U_\theta$ from Lemma~\ref{lem:jet-inverse}, normalized
by \eqref{eq:parameter-corrector-normalization}, and set
\[
 \Phi := U_\theta+\rho^d\psi_\theta,\qquad w := v-\Phi.
\]
Since
\[
 \left|\frac{\rho^d\psi_\theta}{U_\theta}\right|
 \le C[\rho(1+R)]^d=o(1),
\]
we have $\Phi>0$ and $\Phi\asymp(1+r)^{2-n}$ on
$\mathcal A_{\sigma,R}^+$.
Fix a nonnegative, nonzero $\chi\in C_c^\infty(\Gamma\setminus\{0\})$.
The matrix $(\int_\Gamma\chi Z_\alpha Z_\beta)$ is uniformly positive
definite on $\Theta$. Thus the implicit function theorem selects
$\theta\in\Theta$ near $\theta_*$ such that
\[
 \int_\Gamma\chi Z_\alpha w=0\qquad\text{for }0\le\alpha<n.
\]
Write $U:=U_\theta$ and $\psi:=\psi_\theta$, and define
\[
 E:=-L_{g_\rho}\Phi,\qquad
 E^\partial:=B_{g_\rho}\Phi-(n-2)\Phi^{n/(n-2)}.
\]
Equation~\eqref{eq:jet-correction} cancels the order-$\rho^d$ terms.
Taylor expansion of $L_{g_\rho}$ and $B_{g_\rho}$ then yields the bounds
\begin{equation}\label{eq:leading-annular-residual}
    \begin{aligned}
 |E| &\le C\bigl\{\rho^{d+1}(1+r)^{d+1-n}+\rho^{2d}(1+r)^{2d-n}\bigr\},\\
|E^\partial| &\le C\bigl\{\rho^{d+1}(1+r)^{d+2-n}+\rho^{2d}(1+r)^{2d+1-n}\bigr\},
    \end{aligned}
\end{equation}
on $\mathcal A_{\sigma,R}^+$ and $\Gamma_{\sigma,R}$, respectively.
Since $2d\ge d+1$ and $\rho(1+R)\to0$, the second term in each
bound is controlled by the first. Also, $d<(n-2)/2$ implies
$n\ge5$ and $d+1<n-2$.
By \eqref{eq:annular-bounds},
\[
 |w|\le C\quad\text{on }\Sigma_\sigma,\qquad
 |w|\le CR^{2-n}\quad\text{on }\Sigma_R.
\]
Define
\[
 c:=n\int_0^1((1-s)\Phi+sv)^{2/(n-2)}\dd s.
\]
Then $0\le c\le C(1+r)^{-2}$, and $w$ solves
\begin{equation}\label{eq:lower-remainder-equations}
 \begin{cases}
 -\Delta_{g_\rho}w = -E-\dfrac{n-2}{4(n-1)}R_{g_\rho}w&\text{in }\mathcal A_{\sigma,R}^+,\\
 \partial_{\nu_{g_\rho}}w = cw-\dfrac{n-2}{2}h_{g_\rho}w-E^\partial&\text{on }\Gamma_{\sigma,R}.
 \end{cases}
\end{equation}
By Appendix~\ref{app:annular-comparison},
\begin{equation}\label{eq:lower-error}
 |w(y)|+r|\nabla w(y)|\le CW(r)
 \quad\text{for }\sigma\le r=|y|\le R,
\end{equation}
where
\[
 W(r):=\rho^{d+1}(1+r)^{d+3-n}
       +\sigma^{n-2}r^{2-n}+R^{2-n}.
\]

\emph{Step 3. Virial evaluation and comparison.}
Expand
\[
 \log g=\sum_{k=d}^{2d}H_k+O(r^{2d+1}).
\]
Since $2d<n-2$, \eqref{eq:full-first-action-variation} and cutoff
limits imply
\[
 J_1(H_k,U)=0\quad\text{for }d\le k\le2d.
\]
The correction equation implies
\[
 Q(\psi)=-D_fJ_1(H_d,U)[\psi].
\]
By \eqref{eq:reduced-form}, the $\rho^{2d}$ coefficient in
$\mathcal V_\rho(\Phi)$ is
\[
 \begin{aligned}
 -d\,D_fJ_1(H_d,U)[\psi]-2d\,J_2(H_d)
 &=-2d\bigl(J_2(H_d)-\tfrac12Q(\psi)\bigr)\\
 &=-2d\,\mathcal C(H_d;\theta).
 \end{aligned}
\]
The cutoff errors in $\rho^kJ_1(H_k,U)$ from $r<\sigma$ and $r>R$
have size
\[
 O(\rho^k\sigma^{n+k-2}+\rho^kR^{k+2-n})
 \quad\text{for }d\le k\le2d.
\]
The quadratic contributions from these regions have size
\[
 O(\rho^{2d}\sigma^{n-1}+(\rho R)^{2d}R^{2-n}).
\]
Both errors are bounded by $C(\rho R)^d(\sigma^{n-2}+R^{2-n})$.
The Taylor remainder is $o(\rho^{2d})$, so
\begin{equation}\label{eq:lower-profile-dilation}
 \mathcal V_\rho(\Phi)
 =-2d\,\rho^{2d}\mathcal C(H_d;\theta)+o(\rho^{2d})
 +O\bigl((\rho R)^d(\sigma^{n-2}+R^{2-n})\bigr).
\end{equation}
By Appendix~\ref{app:annular-comparison},
\begin{equation}\label{eq:actual-virial-comparison}
 \begin{aligned}
 |\mathcal V_\rho(v)-\mathcal V_\rho(\Phi)|
 \le C(\rho R)^d(\sigma^{n-2}+R^{2-n})+o(\rho^{2d}).
 \end{aligned}
\end{equation}

\emph{Step 4. Deduction of \eqref{eq:increment} and \eqref{eq:increment-recursion}.}
Combining \eqref{eq:first-quotient}, \eqref{eq:lower-profile-dilation},
and \eqref{eq:actual-virial-comparison}, we obtain
\[
 \mathcal V_\rho(v) \ge c\rho^{2d}-o(\rho^{2d})-C(\rho R)^d(\sigma^{n-2}+R^{2-n}).
\]
By \eqref{eq:capacity-matching},
\[
 \sigma_j^{n-2}+R_j^{2-n}\le C(\mu_j+\mu_{j-1}).
\]
Substituting this and \eqref{eq:increment-dilation}, and absorbing
$o(\rho_j^{2d})$, proves \eqref{eq:increment}.
Since $d\ge1$ and $\rho_jR_j<1$ for large $j$,
\eqref{eq:neck-scale-summability} implies
\[
 \sum_{j\gg1}\eta_j\le\sum_{j\gg1}\rho_jR_j<\infty.
\]
For $C\eta_j\le1/2$, \eqref{eq:increment} implies
\[
 \mu_{j-1}\le\frac{1+C\eta_j}{1-C\eta_j}\mu_j
             \le(1+4C\eta_j)\mu_j,
\]
proving
\eqref{eq:increment-recursion}.
\end{proof}

\begin{proof}[Proof of Proposition~\ref{prop:matching-lower}]
Suppose $0$ is nonremovable and the lower bound in \eqref{eq:two-sided} fails. Then $\inf_{t\ge
t_0}\bar V(t)=0$, so \eqref{eq:neck-alternative} and Lemma~\ref{lem:exact-neck-shape} apply in the
gauge of Lemma~\ref{lem:finite-jet-removal}.

\emph{Case (i).}
Iterating \eqref{eq:increment-recursion} from $j$ to $J\ge j$ yields
\[
 0 < \mu_{j-1} \le \mu_J\prod_{k=j}^J(1+C'\eta_k) \le \mu_J\exp\left(C'\sum_{k\ge j}\eta_k\right).
\]
Since $\sum_k\eta_k<\infty$ and $\mu_J\to 0$ by \eqref{eq:capacity-matching}, letting $J\to\infty$
forces $\mu_{j-1}=0$, contradicting $\mu_{j-1}>0$.

\emph{Case (ii).}
Here $g=\delta+O(r^\ell)$ with $\ell\ge(n-2)/2$, so $|\mathcal E|+|\mathcal E^\partial|\le
Ce^{-\ell t}\bar V(t)$, whence $|\frac{\dd}{\dd t}P(e^{-t},u)|\le Ce^{-\ell t}\bar V(t)^2$. Since
$\mathcal P(u)=0$ and $\bar V\le C$, integrating from $t$ to $\infty$ gives $|P(e^{-t},u)|\le
Ce^{-\ell t}$. By \eqref{eq:neck-minimum-limit},
$P(e^{-b_j},u)=-\frac{(n-2)^2}{8}|S|\delta_j^2+o(\delta_j^2)$, so integrating the Pohozaev drift
from $b_j$ to $t_j^+$ yields
\begin{equation}\label{eq:effective-flat-neck}
 c\delta_j^2 \le C\int_{b_j}^{t_j^+}e^{-\ell t}\bar V(t)^2\dd t + Ce^{-\ell t_j^+}.
\end{equation}
On $[b_j,t_j^+]$, \eqref{eq:neck-two-exponentials} and \eqref{eq:neck-minimum-scales} imply $\bar
V(t)\le C\delta_j e^{(n-2)(t-b_j)/2}$, and $t_j^+\le 2b_j-t_j^-+C$. Since $\ell\ge(n-2)/2$,
\[
 \delta_j^{-2}\int_{b_j}^{t_j^+}e^{-\ell t}\bar V(t)^2\dd t
 \le Ce^{-\frac{n-2}{2}b_j}\int_0^{t_j^+-b_j}e^{\frac{n-2}{2}s}\dd s
 \le Ce^{-\frac{n-2}{2}(2b_j-t_j^+)} \le Ce^{-\frac{n-2}{2}t_j^-} \longrightarrow 0.
\]
Absorbing this into \eqref{eq:effective-flat-neck} yields $\delta_j^2\le Ce^{-\ell t_j^+}$. From
Lemma~\ref{lem:exact-neck-shape} ($\delta_j^2\asymp e^{-(n-2)(t_j^+-t_j^-)/2}$), we obtain
\[
 t_j^+-t_j^- \ge \frac{4}{n-2}\log(1/\delta_j)-C \ge \frac{2\ell}{n-2}t_j^+-C.
\]
Consequently,
\[
 t_j^-\le\left(1-\frac{2\ell}{n-2}\right)t_j^++C\le C,
\]
contradicting $t_j^-\to\infty$. This proves Proposition~\ref{prop:matching-lower}.
\end{proof}

\section{Fowler asymptotics}\label{sec:asymptotics}
We prove the Fowler asymptotics in Theorem~\ref{thm:conditional}.
For a nonremovable singularity satisfying \eqref{eq:upper-main},
Proposition~\ref{prop:matching-lower} ensures $0<c\le V\le C$ for
large $t$, in the original boundary Fermi coordinates.
The Pohozaev limit determines the Fowler orbit. We then prove
exponential convergence to this orbit with a limiting phase.

Normalize the cylindrical Pohozaev functional by
\begin{equation}\label{eq:hamiltonian}
 \mathcal H(v;t)=\frac1{|\Sph^{n-2}|}\Bigl\{
 \frac12\int_S\left(v_t^2-|\nabla_Sv|^2-\frac{(n-2)^2}{4}v^2\right)+\frac{(n-2)^2}{2(n-1)}
                  \int_{\partial S}v^{2(n-1)/(n-2)}\Bigr\}.
\end{equation}
Then $\mathcal H(V;t)=P(e^{-t},u)/|\Sph^{n-2}|$.
Put $H_\infty=\mathcal P(u)/|\Sph^{n-2}|$.
Lemma~\ref{lem:cylindrical-estimates} implies
\[
 |\mathcal H(V;t)-H_\infty|\le Ce^{-t}.
\]
For every $t_j\to\infty$, a subsequence of $V(t_j+\cdot)$
converges smoothly on compact cylinders to an entire flat solution
$v$ with $0<c\le v\le C$ and $\mathcal H(v)=H_\infty$.
The symmetry theorem \cite[Theorem 1.3]{CJSX}, with boundary dimension
$n-1$ and order $1/2$, and the classification
\cite[Theorem 1.1]{JLXBoundary} give
\[
 H_\infty\in[H_c,0),\qquad
 v=v_{H_\infty}(\cdot+\tau).
\]
Set $H=H_\infty$ and
$\mathcal O_H=\{v_H(\cdot+\tau):\tau\in\R\}$.
Compactness now ensures convergence to $\mathcal O_H$ on every fixed
cylinder. To obtain a rate and a limiting phase, put $T=T(H)$ and
$J=[-2T,2T]$ when $H_c<H<0$.
At $H=H_c$, use $J=[-2,2]$ and $\mathcal O_{H_c}=\{v_c\}$.
Fix $0<\alpha<1$. For sufficiently large $s,t$, define
\[
 d(s)=\dist_{C^{2,\alpha}(J\times\overline S)}
              (V(s+\cdot),\mathcal O_H),\qquad
 D(t)=\sup_{s\ge t}d(s).
\]
The supremum $D$ is nonincreasing and tends to zero.
The following contraction extends the phase argument of
\cite[Section 4]{JLXBoundary} to the decaying sources in
\eqref{eq:forced-cylinder}.

\begin{lemma}\label{lem:tail-contraction}
Let $V$ solve \eqref{eq:forced-cylinder} with the two-sided bounds
and the $O(e^{-t})$ source and Pohozaev estimates above.
Suppose $D(t)\to0$.
For $H_c<H<0$ and every $0<\gamma<1$, there are $L,t_0>0$ such that,
whenever $t\ge t_0$ and $D(t)>e^{-\gamma t}$,
\begin{equation}\label{eq:tail-halving}
 D(t+L)\le\tfrac12D(t).
\end{equation}
At $H=H_c$, the conclusion holds for $0<\gamma<1/2$.
\end{lemma}
\begin{proof}
\emph{Step 1. Phase selection.}
For $H_c<H<0$, set $T:=T(H)$. Since $v_H$ is nonconstant and $J$
contains full periods,
\[
 \inf_{\tau\in\R/T\mathbb Z}
 \|v_{H,t}(\cdot+\tau)\|_{L^2(J\times S)}>0.
\]
When $d(s)$ is small, minimize
\[
 \frac12\|V(s+\cdot)-v_H(\cdot+\tau)\|_{L^2(J\times S)}^2
 \quad\text{over }\tau\in\R/T\mathbb Z.
\]
At a minimizing phase, the second derivative satisfies
\[
 \|v_{H,t}\|_2^2-\langle V-v_H,v_{H,tt}\rangle
 \ge\frac12\|v_{H,t}\|_2^2>0.
\]
Compactness of the phase circle and the least-period condition
separate the other phases. The implicit function theorem therefore
yields a unique minimizing phase $\tau(s)\in\R/T\mathbb Z$ satisfying
\begin{equation}\label{eq:phase-condition}
 \bigl\langle V(s+\cdot)-v_H(\cdot+\tau(s)),\ v_{H,t}(\cdot+\tau(s))\bigr\rangle_{L^2(J\times S)} = 0,
\end{equation}
and
\[
 d(s)\le
 \|V(s+\cdot)-v_H(\cdot+\tau(s))\|_{C^{2,\alpha}(J\times\overline S)}
 \le Cd(s).
\]
If $|h|\le T$, the intervals centered at $s$ and $s+h$ overlap in
at least one full period. Comparing both models with $V$ on the
overlap, we obtain
\begin{equation}\label{eq:phase-overlap}
 \dist_{\R/T\mathbb Z}\bigl(\tau(s+h),\ \tau(s)+h\bigr) \le C\bigl(d(s)+d(s+h)\bigr).
\end{equation}
Partitioning the interval between $s$ and $s+h$ into subintervals
of length at most $T$ and summing \eqref{eq:phase-overlap} yields
\[
 \dist_{\R/T\mathbb Z}(\tau(s+h),\tau(s)+h)
 \le C(1+|h|)D(t)\quad\text{for }s,s+h\ge t.
\]
For $H=H_c$, set $\tau(s):=0$.

\emph{Step 2. Jacobi limit.}
Suppose \eqref{eq:tail-halving} fails. Choose $L_j\to\infty$,
$t_j\to\infty$, and $s_j$ such that
\[
 s_j\ge t_j+L_j,\qquad D_j:=D(t_j)>e^{-\gamma t_j},\qquad
 d(s_j)>D_j/2.
\]
Set
\[
 \Psi_j(t,\omega) := \frac{V(s_j+t,\omega)-v_H(t+\tau(s_j),\omega)}{D_j}.
\]
Since $s_j-t_j\to\infty$, we can cover $[h-1,h+1]$ by translates
of $J$ and apply Step 1. For each fixed $h$ and all large $j$,
\[
 \|\Psi_j\|_{C^0([h-1,h+1]\times\overline S)}
 \le\begin{cases}
 C(1+|h|),&\text{if }H_c<H<0,\\
 C,&\text{if }H=H_c,
 \end{cases}
\]
with $C$ independent of $h$. Subtracting \eqref{eq:flat-cylinder}
from \eqref{eq:forced-cylinder} and dividing by $D_j$, we obtain
\[
 \begin{cases}
 (\Psi_j)_{tt}+\Delta_S\Psi_j-\dfrac{(n-2)^2}{4}\Psi_j = \dfrac{\mathcal E(s_j+\cdot)}{D_j}&\text{in }(-s_j,\infty)\times S,\\
 \partial_{\nu_S}\Psi_j = c_j\Psi_j + \dfrac{\mathcal E^\partial(s_j+\cdot)}{D_j}&\text{on }(-s_j,\infty)\times\partial S,
 \end{cases}
\]
where
\[
 c_j:=n\int_0^1\bigl((1-\sigma)v_H(t+\tau(s_j),\omega)
              +\sigma V(s_j+t,\omega)\bigr)^{2/(n-2)}\dd\sigma.
\]
Since
\[
 \frac{e^{-s_j}}{D_j}
 \le e^{-(1-\gamma)t_j}e^{-(s_j-t_j)}\longrightarrow0,
\]
the local Robin estimates in Appendix~\ref{app:local-robin} provide
uniform $C^{3,\alpha}$ bounds on compact cylinders.
Passing to a subsequence with $\tau(s_j)\to\tau_*$, we have
\[
 \Psi_j\to\Psi\quad\text{in }
 C^{2,\alpha}_{\mathrm{loc}}(\R\times\overline S).
\]
The limit solves the linearized equation at $v_H(\cdot+\tau_*)$
and satisfies $|\Psi(t,\omega)|\le C(1+|t|)$. Moreover,
\[
 \begin{aligned}
 \|\Psi\|_{C^{2,\alpha}(J\times\overline S)}
 &=\lim_{j\to\infty}\|\Psi_j\|_{C^{2,\alpha}(J\times\overline S)}\\
 &\ge\liminf_{j\to\infty}\frac{d(s_j)}{D_j}\ge\frac12.
 \end{aligned}
\]

\emph{Step 3. Pohozaev exclusion.}
If $H_c<H<0$, \cite[Theorem 1.3]{JLXBoundary} and the linear
growth bound imply
\[
 \Psi=c\,v_{H,t}(\cdot+\tau_*)
       +e\,\partial_H v_H(\cdot+\tau_*).
\]
Expanding $\mathcal H(V;s_j)-H$ at $v_H(\cdot+\tau(s_j))$
on $t=0$, we obtain
\[
 \begin{aligned}
 D\mathcal H(v_H(\cdot+\tau(s_j)))[\Psi_j](0)
 &=\frac{\mathcal H(V;s_j)-H}{D_j}+O(D_j)\\
 &\longrightarrow0.
 \end{aligned}
\]
Since
\[
 D\mathcal H(v_H)[v_{H,t}]=0,\qquad
 D\mathcal H(v_H)[\partial_H v_H]=1,
\]
we obtain $e=0$. Passing \eqref{eq:phase-condition} divided by
$D_j$ to the limit yields
\[
 0=\langle\Psi,v_{H,t}(\cdot+\tau_*)\rangle
   =c\|v_{H,t}\|_{L^2(J\times S)}^2.
\]
Thus $c=0$ and $\Psi\equiv0$, contradicting
$\|\Psi\|_{C^{2,\alpha}}\ge1/2$.

If $H=H_c$, then $v_{c,t}=0$ and $D\mathcal H(v_c)=0$.
The second-order expansion is
\[
 \tfrac12 D^2\mathcal H(v_c)[\Psi_j,\Psi_j](0) = \frac{\mathcal H(V;s_j)-H_c}{D_j^2} + O(D_j).
\]
Because $\gamma<1/2$,
\[
 \frac{e^{-s_j}}{D_j^2}
 \le e^{-(1-2\gamma)t_j}e^{-(s_j-t_j)}\longrightarrow0,
\]
so $D^2\mathcal H(v_c)[\Psi,\Psi]=0$.
By \cite[Theorem 1.3]{JLXBoundary},
\[
 \Psi(t,\omega)=\phi_c(\omega)
       \bigl(c_1\cos(\kappa_c t)+c_2\sin(\kappa_c t)\bigr),
 \qquad \|\phi_c\|_{L^2(S)}=1,
\]
where
\[
 \begin{cases}
 (-\Delta_S+\tfrac{(n-2)^2}{4})\phi_c=-\kappa_c^2\phi_c
     &\text{in }S,\\
 \partial_{\nu_S}\phi_c=nv_c^{\frac2{n-2}}\phi_c
     &\text{on }\partial S.
 \end{cases}
\]
Substituting into \eqref{eq:hamiltonian} yields
\[
 \begin{aligned}
 \tfrac12D^2\mathcal H(v_c)[\Psi,\Psi]
 &=\frac{\|\Psi_t\|_{L^2(S)}^2
            +\kappa_c^2\|\Psi\|_{L^2(S)}^2}{2|\Sph^{n-2}|}\\
 &=\frac{\kappa_c^2}{2|\Sph^{n-2}|}(c_1^2+c_2^2)>0.
 \end{aligned}
\]
Since $\Psi\not\equiv0$, this contradicts
$D^2\mathcal H(v_c)[\Psi,\Psi]=0$.
\end{proof}

Monotonicity of $D$ and Lemma~\ref{lem:tail-contraction} imply
$D(t+L)\le\frac12D(t)+e^{-\gamma t}$. Choose
$0<\beta<\min(\gamma,\log2/L)$ and set $t_k=t_0+kL$. Iteration yields
\[
 D(t_k)\le2^{-k}D(t_0)+
 C\sum_{j=0}^{k-1}2^{-(k-1-j)}e^{-\gamma t_j}
 \le C_\beta e^{-\beta t_k}.
\]
Monotonicity extends this bound to all sufficiently large $t$.
For $H_c<H<0$, put $s_k=s_0+kT$. Equation~\eqref{eq:phase-overlap}
allows successive real phase lifts with
\[
 |(\tau(s_{k+1})-s_{k+1})-(\tau(s_k)-s_k)|
 \le Ce^{-\beta s_k}.
\]
Their increments are summable. Thus a limiting phase $\tau_\infty$
exists, with $|\tau(s_k)-s_k-\tau_\infty|\le Ce^{-\beta s_k}$.
Another application of \eqref{eq:phase-overlap} controls the intermediate
times. At $H=H_c$, take $\tau_\infty=0$. In both cases,
\begin{equation}\label{eq:cylindrical-asymptotics}
 \|V(t+\cdot)-v_H(t+\cdot+\tau_\infty)\|
 _{C^{2,\alpha}([0,1]\times\overline S)}\le Ce^{-\beta t}.
\end{equation}
The difference equation has smooth Robin coefficients and sources
$O(e^{-t})$, so local regularity yields the same rate in each fixed
$C^k$ norm. Returning to $x=r\omega$ through \eqref{eq:model} and
$r\partial_r=-\partial_t$ proves \eqref{eq:main-asymptotics}.
The positive lower bound of $v_H$ implies \eqref{eq:main-relative}.
The limiting energy fixes $H$. Two limiting phases give the same
profile on a full period and therefore agree modulo $T(H)$.
At $H=H_c$ the model is independent of phase. Together with
Proposition~\ref{prop:matching-lower}, this proves
Theorem~\ref{thm:conditional}.

\appendix

\section{Linear metric correction estimates}\label{app:correction-estimates}

We expand \eqref{eq:action} to obtain the sources used in
Lemma~\ref{lem:jet-inverse}.
Let $H$ be a homogeneous Fermi-preserving tensor of degree $d$, and
write $\tau_H=\tr H$. The inverse metric density and boundary area
density have first variations $\frac12\tau_H\delta^{ij}-H_{ij}$
and $\frac12\tr_TH$, respectively. The curvature variations are
\[
 DR_\delta[H]=\partial_i\partial_jH_{ij}-\Delta\tau_H,\qquad
 Dh_\delta[H]=-\frac{\partial_{x_n}\tau_H}{2(n-1)}.
\]

For a positive function $f$, the first energy variation is
\begin{equation}\label{eq:full-first-action-variation}
\begin{split}
 J_1(H,f)={}&\frac12\int\left[
       (\tfrac12\tau_H\delta^{ij}-H_{ij})f_if_j
       +\frac{n-2}{4(n-1)}DR_\delta[H]f^2\right]\\
 &+\frac{n-2}{4}\int_\Gamma Dh_\delta[H]f^2
       -\frac{(n-2)^2}{4(n-1)}\int_\Gamma(\tr_TH)f^{\frac{2(n-1)}{n-2}}.
\end{split}
\end{equation}
Differentiation in $f$ at $U$ recovers \eqref{eq:metric-correction-source}.
The conormal and boundary-area terms cancel as follows.
On $\Gamma$, $H_{in}=0$ and $\tau_H=\tr_TH$, so
\[
 \bigl((\tfrac12\tau_H I-H)\nabla U\bigr)\cdot\nu
       -\frac{n-2}{2}\tau_H U^{\frac n{n-2}}
 =\frac12\tau_H\bigl(\partial_\nu U-(n-2)U^{\frac n{n-2}}\bigr)=0.
\]

\begin{proof}[Decay estimate for Lemma~\ref{lem:jet-inverse}]
Let $1\le d<(n-2)/2$, as in Lemma~\ref{lem:jet-inverse}, and
write $\psi_H(x)=U(x)v(T(x))$.
Under conformal covariance, the equation becomes
\[
 \begin{cases}
 -\Delta v=\tfrac14
       (U^{-\frac{n+2}{n-2}}F_H)\circ T^{-1}
       &\text{in }\mathcal B_1,\\
 (\partial_\nu-1)v=\tfrac12
       (U^{-\frac n{n-2}}F_H^\partial)\circ T^{-1}
       &\text{on }\partial\mathcal B_1\setminus\{e_n\}.
 \end{cases}
\]
The omitted point is $e_n$, and
\[
 |T(x)-e_n|=\frac2{\sqrt{|x'|^2+(1+x_n)^2}}
       \asymp(1+|x|)^{-1}.
\]
Thus, for $\rho=|y-e_n|<1/4$, the volume and boundary sources
are bounded by $C|H|\rho^{-d-2}$ and
$C|H|\rho^{-d-1}$, respectively.
The inverse Robin kernel on the complement of the degree-one
harmonics is bounded by $C|y-z|^{2-n}+C$.
The degree-one harmonics added to impose \eqref{eq:corrector-normalization}
are bounded by $C|H|$.
The source is locally integrable at $e_n$ because $d<n-2$.
Split its potential into $|z-e_n|<\rho/2$, the intervening region,
and $|z-e_n|>2\rho$. In the intervening region the source is
$O(|H|\rho^{-d-2})$ and the kernel integral is $O(\rho^2)$.
The three contributions combine to yield
\begin{equation}\label{eq:compactified-correction-bound}
 |v(y)|\le C|H|\left(
 \rho^{2-n}\int_0^{\rho/2}s^{n-d-3}\dd s
 +\rho^{-d}+\int_{2\rho}^1s^{-d-1}\dd s\right)
 \le C|H|\rho^{-d}.
\end{equation}
The boundary potential has the same powers of $s$ in the inner and
outer integrals, and its middle part is also bounded by
$C|H|\rho^{-d}$.
\end{proof}

\section{Analytic estimates and constructions}

\subsection{Mixed Green kernels and weak extension}\label{app:green-kernels}
For \(v|_{\Sigma_r}=0\), scaling Poincaré and trace inequalities imply
\[
\int_{\mathcal B_r^+}v^2\le Cr^2\int_{\mathcal B_r^+}|\nabla v|^2,
\qquad \int_{B_r}v^2\le Cr\int_{\mathcal B_r^+}|\nabla v|^2.
\]
Consequently the mixed form satisfies
\begin{equation}\label{eq:small-mixed-coercivity}
    Q_g(v,v)\ge(c-Cr-Cr^2)\|\nabla v\|_2^2.
\end{equation}
The estimates below are used for both expanding balls and
half-annuli. All coefficient bounds are uniform after pulling
back to one fixed coordinate neighborhood. Thus
\(g_\rho(x)=g(\rho x)\) on \(\mathcal B_R^+\), with \(\rho R\le\delta_0\),
has ellipticity constants independent of \(R\) and \(\rho\).
The constants below are allowed to depend on the fixed metric.

Let \(G_0\) be the Green kernel for \(-\Delta_{g_\rho}\) on a
half-ball or half-annulus, with Neumann condition on the flat
face and Dirichlet condition on the spherical boundaries.
Write \(C^{ij}=\sqrt{\det g_\rho}\,(g_\rho)^{ij}\). Reflect the
tangential block and the normal coefficient \(C^{nn}\) evenly, and
the normal-tangential coefficients \(C^{in}\) oddly. The Fermi
conditions make the reflected divergence-form coefficients
Lipschitz and uniformly elliptic. The reflected domain has
Dirichlet condition on its spherical boundaries. The Dirichlet
Green estimates \cite{LSW,GW} establish
\begin{equation}\label{eq:mixed-green-upper}
    0<G_0(x,y)\le C|x-y|^{2-n}.
\end{equation}
For the two-sided estimates, take the half-ball \(\mathcal B_R^+\).
In the reflected ball put \(d_R(x)=\dist(x,\partial \mathcal B_R)\) and
\[
K_R(x,y)=|x-y|^{2-n}
\left(1\wedge\frac{d_R(x)}{|x-y|}\right)
\left(1\wedge\frac{d_R(y)}{|x-y|}\right).
\]
The interior Green lower bound of \cite{LSW,GW}, followed by Harnack
chains at distance comparable to \(|x-y|\) from the spherical boundary,
gives the lower estimate away from the spherical boundary. To reach
the spherical boundary \(\Sigma_R\), let \(\ell\) be the distance to
\(\Sigma_R\) in a patch of radius \(r\le R/8\) containing no pole.
The coefficients are Lipschitz, \(|\nabla\ell|=1\), and
\(|D^2\ell|\le C/R\). For fixed large \(K\), on a smaller patch one has
\[
-\Delta_{g_\rho}(\ell+K\ell^2/r)<0
<-\Delta_{g_\rho}(\ell-K\ell^2/r).
\]
Indeed, \(|\Delta_{g_\rho}\ell|\le C(R^{-1}+\rho)\) and
\(|\nabla\ell|_{g_\rho}^2\ge c\). On \(\ell\le r/(4K)\), the quadratic
term of size \(K/r\) dominates when \(K\) is large and
\(\rho R\le\delta_0\). Comparison, interior Harnack and scaling yield
\[
\frac{w(x)}{w(x_r)}\asymp\frac{d_R(x)}r
\]
for a positive harmonic function vanishing on the spherical boundary.
Here \(x_r\) is an interior point at distance comparable to \(r\).
Apply this first in \(x\) and then in \(y\), using symmetry. Thus
\(G_0\asymp K_R\), with constants independent of \(\rho,R\) when
\(\rho R\le\delta_0\).

By even reflection of the Neumann condition, the half-ball kernel is
the sum of the reflected Dirichlet kernels with poles at \(y\) and its
reflection. For a face pole the two terms coincide. The
boundary-distance comparison yields
\begin{subequations}\label{eq:green-three}
    \begin{equation}
        G_0(0,y)\asymp|y|^{2-n}-R^{2-n}.
    \end{equation}
    For three points, it implies
    \begin{equation}
        \frac{G_0(x,z)G_0(z,y)}{G_0(x,y)}
        \le C\bigl(|x-z|^{2-n}+|z-y|^{2-n}\bigr).
    \end{equation}
\end{subequations}
The boundary pole is the limit of interior poles after reflection.
The second inequality follows from the two-sided Dirichlet estimate
by splitting according to which of \(|x-z|,|z-y|\) is at least
\(|x-y|/2\) and retaining the boundary-distance factors on that side.

The same estimates hold for the metric pair \((-L_{g_\rho},B_{g_\rho})\)
on a sufficiently small coordinate neighborhood before rescaling. In fact
\(|{\frac{n-2}{4(n-1)}}R_{g_\rho}|\le C\rho^2\) and
\(|{\frac{n-2}{2}}h_{g_\rho}|\le C\rho\).
By \eqref{eq:green-three},
\begin{align*}
    &\sup_{x,y}\frac{
        \int G_0(x,z)|{\frac{n-2}{4(n-1)}}R_{g_\rho}(z)|G_0(z,y)\dd z
        +\int_\Gamma G_0(x,z)|{\frac{n-2}{2}}h_{g_\rho}(z)|G_0(z,y)\dd z}
    {G_0(x,y)}\\
    &\qquad\le C\sup_x\left(
    \rho^2\int_{\mathcal B_R^+}|x-z|^{2-n}\dd z
    +\rho\int_{B_R}|x-z|^{2-n}\dd z\right)
    \le C((\rho R)^2+\rho R).
\end{align*}
Choose \(\delta_0\) so that the last bound is less than \(1/4\).
Iterating this bound in the signed volume/face resolvent series gives
\[
 |G-G_0|\le\sum_{m=1}^\infty4^{-m}G_0=\tfrac13G_0.
\]
Thus \(0<G\asymp G_0\). The same smallness makes the mixed form
coercive. This argument handles the signed scalar and mean-curvature
potentials introduced by the finite conformal gauges.
On a half-annulus, coercivity again ensures a positive mixed kernel.
Comparison on the annulus with the half-ball kernel implies
\eqref{eq:mixed-green-upper} for this metric kernel as well.

The elementary potential bounds needed below are as follows. For
\(1\le d<n-2\) put \(\phi_d=(1+r)^{d+2-n}\). Then
\begin{equation}\label{eq:potential-powers}
    \begin{gathered}
        \int_{\R^n_+}|x-y|^{2-n}(1+|y|)^{d-n}\dd y
        \le C\phi_d(x),\\
        \int_{\R^{n-1}}|x-z|^{2-n}(1+|z|)^{d+1-n}\dd z
        \le C\phi_d(x).
    \end{gathered}
\end{equation}
For \(r=|x|\ge2\), the regions \(|y|<r/2\), \(r/2<|y|<2r\) and
\(|y|>2r\) yield the combined volume bound
\[
C r^{2-n}\int_0^{r/2}(1+s)^{d-n}s^{n-1}\dd s
+Cr^{d-n}r^2+C\int_{2r}^\infty s^{d+1-n}\dd s
\le Cr^{d+2-n}.
\]
For the face the middle kernel integral is \(O(r)\), leading to the same
power. The integral over $s>2r$ converges because \(d<n-2\).
Local kernel integrability treats \(r<2\).

Put \(K_0=(1+r)^{2-n}\).
For any weight
\begin{equation}\label{eq:four-annular-weights}
    \phi\in\{\phi_d,r^{2-n},1,K_0\},
\end{equation}
the same splitting on \(\mathcal B_R^+\) yields
\begin{equation}\label{eq:metric-potential-small}
    G_0[\rho^2\phi]+G_{0,\Gamma}[\rho\phi]
    \le C((\rho R)^2+\rho R)\phi.
\end{equation}
Finally, if
\(0\le c\le C(1+r)^{-2}\), its face potential restricted to
\(r<\eta\) or \(r>L\) has operator norm at most
\(C(\eta+L^{-1})\) on each weight in
\eqref{eq:four-annular-weights}. Near zero this follows by
scaling the order-one face potential. At infinity its coefficient
is bounded by \(Cr^{-2}\), so scaling produces \(L^{-1}\).
The three-region calculation above establishes these bounds uniformly
also when the point \(x\) lies outside the source region.

On the original half-ball, choose \(r\) so that
\eqref{eq:small-mixed-coercivity} is coercive. The Lax--Milgram
theorem gives a unique weak solution of
\[
\begin{cases}
    -L_g\varphi=0&\text{in }\mathcal B_r^+,\\
    B_g\varphi=0&\text{on }B_r,\\
    \varphi=1&\text{on }\Sigma_r.
\end{cases}
\]
Testing with \(\varphi_-\) and using coercivity forces \(\varphi_-=0\),
hence \(\varphi\ge0\). The strong maximum principle and boundary
point lemma ensure \(\varphi>0\). Write \(G[f]\) and
\(G_\Gamma[f^{\partial}]\) for the volume and face Green
potentials. The Green representation with the kernel just constructed yields
\begin{equation}\label{eq:positive-gauge-bound}
    \varphi-1=-G\left[\frac{n-2}{4(n-1)}R_g\right]
    -G_\Gamma\left[\frac{n-2}{2}h_g\right].
\end{equation}
From the kernel upper bound,
\(\int G\le Cr^2\) and \(\int G_\Gamma\le Cr\). Hence
\[
\|\varphi-1\|_\infty\le C(r^2\|R_g\|_\infty+r\|h_g\|_\infty).
\]
After a further shrinking, \(1/2\le\varphi\le3/2\).
Local boundary regularity makes $\varphi$ smooth on every smaller closed half-ball. For
\(\widetilde g=\varphi^{\frac{4}{n-2}}g\) and
\(\widetilde u=\varphi^{-1}u\), conformal covariance yields
\begin{equation}\label{eq:covariance}
    \begin{cases}
        -L_{\widetilde g}(\varphi^{-1}v)=-\varphi^{-{\frac{n+2}{n-2}}}L_gv
        &\text{in }\mathcal B_r^+,\\
        B_{\widetilde g}(\varphi^{-1}v)=\varphi^{-{\frac{n}{n-2}}}B_gv
        &\text{on }B_r.
    \end{cases}
\end{equation}
Thus this gauge is scalar flat with minimal boundary. Fixed positive
gauges preserve critical growth, removability, and relative asymptotics.

For \(n\ge3\), take \(r\) sufficiently small. Let \(G\) be the positive
Green kernel of \((-L_g,B_g)\) on \(\mathcal B_r^+\) with zero Dirichlet
boundary value on \(\Sigma_r\), and let \(h_0\) be the
\((-L_g,B_g)\)-harmonic lift of \(u|_{\Sigma_r}\), with zero
boundary datum on \(B_r\).

\begin{lemma}\label{lem:no-source}
Every positive punctured solution satisfies
\[
 u^{\frac{n}{n-2}}\in L^1_{\mathrm{loc}}(B_r),
\]
and its weak equation extends across the puncture with zero point
source. More precisely, every
$\eta\in C_c^\infty(\mathcal B_r^+\cup B_r)$ satisfies
\[
 \int_{\mathcal B_r^+}\left(
 \langle\nabla u,\nabla\eta\rangle_g
       +\frac{n-2}{4(n-1)}R_gu\eta\right)\dd V_g
       +\frac{n-2}{2}\int_{B_r}h_gu\eta\dd\sigma_g
 =(n-2)\int_{B_r}u^{\frac n{n-2}}\eta\dd\sigma_g.
\]
On a smaller half-ball,
\begin{equation}\label{eq:green-representation}
 u(x)={h_0}(x)+{(n-2)}\int_{B_r}G(x,z)u(z)^{\frac{n}{n-2}}\dd\sigma_g(z).
\end{equation}
\end{lemma}

\begin{proof}[Proof of Lemma~\ref{lem:no-source}]
Apply the gauge $\varphi$ from Appendix~\ref{app:green-kernels} so that
$R_{\widetilde g}=0$ and $h_{\widetilde g}=0$. Reflect
$\widetilde u:=\varphi^{-1}u$ evenly across $\{x_n=0\}$ to $V$, and
reflect $C^{ij}:=\sqrt{\det\widetilde g}\,\widetilde g^{ij}$ as above.
Away from zero,
\[
 -\diver(C\nabla V)=\mu_0,
\]
where $\mu_0$ is defined on the punctured reflected face by
\[
 \mu_0:=2(n-2)\sqrt{\det\widetilde g_T}\,
 V^{\frac{n}{n-2}}\mathcal H^{n-1}\!\lfloor\{x_n=0\}.
\]
Let $G_0$ be the Dirichlet Green kernel of $-\diver(C\nabla\cdot)$ on
a smaller reflected ball $\mathcal B$. Put
$\mu_j:=\mathbf 1_{\{|z|>1/j\}}\mu_0$ and $P_j:=G_0\mu_j$.
For fixed $j$, $P_j$ is bounded near zero and
$G_0(x,0)\asymp|x|^{2-n}$. Comparison on punctured balls gives
\[
 P_j\le V+\eta G_0(\cdot,0)\qquad\text{for every }\eta>0.
\]
Letting $\eta\downarrow0$ yields $P_j\le V$.
For a concentric ball $\mathcal B_r\Subset\mathcal B$, evaluation
at a fixed point $x_*$ away from the face yields
\[
 c\mu_j(\mathcal B_r)\le P_j(x_*)\le V(x_*).
\]
Monotone convergence and smoothness of the density away from zero yield
\[
 \mu_0(\mathcal B)<\infty,\qquad
 \mu_0(\{0\})=0.
\]
Set $P:=G_0\mu_0$. The inequality $P_j\le V$ passes to the limit, so $P\le V$.

The reflected coefficients are Lipschitz and uniformly elliptic.
The nonnegative function $w:=V-P$ solves
$-\diver(C\nabla w)=0$ in $\mathcal B\setminus\{0\}$.
By \cite[Theorem 1]{Serrin}, either $w$ is removable at zero or
$w(x)\asymp|x|^{2-n}$. The latter would imply
\[
 \int_{|z|<r}V(z,0)^{\frac{n}{n-2}}\dd z
 \ge c\int_0^r t^{-2}\dd t=\infty,
\]
contradicting $\mu_0(\mathcal B)<\infty$.
Hence $w$ is weakly $C$-harmonic across zero.
By \cite[Theorems 5.1 and 6.1]{LSW} (or \cite{GW}), for
$1\le p<n/(n-1)$,
\[
 \|\nabla P\|_{L^p(\mathcal B_{r/2})}
 \le\mu_0(\mathcal B)
       \sup_y\|\nabla_xG_0(\cdot,y)\|_{L^p(\mathcal B_{r/2})}<\infty.
\]
Thus $V\in W^{1,p}_{\mathrm{loc}}$ and
\[
 -\diver(C\nabla V)=-\diver(C\nabla P)=\mu_0
 \quad\text{in }\mathcal D'(\mathcal B).
\]
Both $P$ and $w$ are even. On the half-ball they are the face-source
potential and the harmonic lift of the spherical data by
\cite[Theorem 6.1]{LSW}. The conformal identities
\[
 \begin{aligned}
 G_g(x,y)&=\varphi(x)\varphi(y)G_{\widetilde g}(x,y),\\
 \widetilde u^{\frac{n}{n-2}}\dd\sigma_{\widetilde g}
     &=\varphi u^{\frac{n}{n-2}}\dd\sigma_g
 \end{aligned}
\]
then establish the weak extension and \eqref{eq:green-representation}.
\end{proof}

\subsection{Boundary potentials and completeness}

Work in a positive scalar-flat, minimal-boundary gauge of
\eqref{eq:main}, so that $R_g=0$ and $h_g=0$.
The boundary density
\[
 f(y)=(n-2)u(y,0)^{n/(n-2)}\sqrt{\det g_T(y)}\ge0
\]
belongs to $L^1(B_\rho)$ by Lemma~\ref{lem:no-source}.
Define
\[
 V_f(r)=\int_{B_{\rho/2}}(r+|y|)^{2-n}f(y)\dd y.
\]

\begin{proposition}\label{prop:potential-completeness}
For $n\ge3$ in this gauge and $0<r<r_0<\rho/4$,
\begin{equation}\label{eq:potential-minimum}
 c(1+V_f(r))
 \le\min_{\substack{|x|=r\\x_n\ge0}}u(x)
 \le u(0,r)\le C(1+V_f(r)).
\end{equation}
The conformal metric $u^{4/(n-2)}g$ is complete at zero if and only if
\begin{equation}\label{eq:potential-completeness}
 \int_0^{r_0}(1+V_f(r))^{2/(n-2)}\dd r=\infty.
\end{equation}
Either $u\to\infty$ uniformly at zero, or there is $L\in(0,\infty)$
such that, for every $0<\eta<1$,
\[
 \lim_{\substack{x\to0\\x_n\ge\eta|x|}}u(x)=L.
\]
In particular, completeness of $u^{4/(n-2)}g$ at zero implies
$u\to\infty$ uniformly.
\end{proposition}
\begin{proof}
By \eqref{eq:green-representation},
\[
 u(x)=h_0(x)+\int_{B_\rho}G(x,y)f(y)\dd y,
 \qquad c\le h_0\le C,
\]
where the bounds on $h_0$ follow from the maximum principle.
For $|x|<\rho/4$ and $|y|\ge\rho/2$, we have $|x-y|\ge\rho/4$, so
\[
 \int_{B_\rho\setminus B_{\rho/2}}G(x,y)f(y)\dd y
 \le C\rho^{2-n}\|f\|_{L^1(B_\rho)}.
\]
For $|x|=r<\rho/4$ and $|y|<\rho/2$,
\[
 |x-y|\le r+|y|\asymp|(0,r)-y|.
\]
The Green bounds \eqref{eq:mixed-green-upper} and
\eqref{eq:green-three} therefore imply \eqref{eq:potential-minimum}.

By \eqref{eq:potential-minimum}, the inward normal ray
$\{x=(0,r):0<r<r_0\}$ has conformal length at most
\[
 C\int_0^{r_0}(1+V_f(r))^{\frac{2}{n-2}}\dd r.
\]
Conversely, for any locally rectifiable curve
$\alpha:[0,T)\to\overline{\mathcal B_{r_0}^+}\setminus\{0\}$ with $|\alpha(0)|=r_0$ and
$\lim_{t\to T^-}\alpha(t)=0$, setting $r(t):=|\alpha(t)|$, ellipticity and
\eqref{eq:potential-minimum} imply
\[
 L_{u^{4/(n-2)}g}(\alpha|_{[0,t]})
 \ge c\int_0^t\bigl(1+V_f(r(s))\bigr)^{\frac{2}{n-2}}|r'(s)|\dd s
 \ge c\int_{r(t)}^{r_0}\bigl(1+V_f(s)\bigr)^{\frac{2}{n-2}}\dd s.
\]
Letting $t\to T^-$ proves \eqref{eq:potential-completeness}.

By monotone convergence,
\[
 \lim_{r\downarrow0}V_f(r)
 =\int_{B_{\rho/2}}|y|^{2-n}f(y)\dd y.
\]
If this limit is infinite, \eqref{eq:potential-minimum} gives
$u(x)\to\infty$ uniformly as $x\to0$.
If it is finite, then $\int_{B_\rho}|y|^{2-n}f(y)\dd y<\infty$.
On $\{x_n\ge\eta|x|\}$,
\[
 G(x,y)\le C_\eta|y|^{2-n},\qquad
 G(x,y)\longrightarrow G(0,y)\quad\text{for }y\ne0.
\]
By dominated convergence,
\[
 u(x)\longrightarrow
 h_0(0)+\int_{B_\rho}G(0,y)f(y)\dd y=:L\in(0,\infty).
\]
In this case the normal ray has finite conformal length, so completeness
forces $u(x)\to\infty$ uniformly as $x\to0$.
\end{proof}

\subsection{Local Robin and critical trace estimates}\label{app:local-robin}

Let $D'\Subset D\cup\Gamma_D$ be nested half-balls or cylinder strips,
where $\Gamma_D$ is the smooth boundary part. Their separation from
$\partial D\setminus\Gamma_D$ is fixed and positive, and $D$ avoids
the edge. Assume uniform ellipticity and fixed smooth coefficient
bounds. The constants below depend on these bounds and the separation.

For the linear problem
\[
 \begin{cases}
 Aw=f&\text{in }D,\\
 \partial_\nu w+\kappa w=f^{\partial}&\text{on }\Gamma_D,
 \end{cases}
\]
with smooth bounded $\kappa$, the local oblique estimates of
\cite{ADN,Lieberman} establish, for every $D'\Subset D\cup\Gamma_D$,
\begin{equation}\label{eq:buffered-robin-estimate}
 \|w\|_{C^{k+2,\alpha}(D')}
 \le C\left(\|w\|_{L^2(D)}+\|f\|_{C^{k,\alpha}(D)}
                    +\|f^{\partial}\|_{C^{k+1,\alpha}(\Gamma_D)}\right).
\end{equation}
The obliqueness constant has a uniform positive lower bound.
For $q\ge1$, the local linear maximum estimate and interpolation
on nested patches also imply
\begin{equation*}
 \|w\|_{L^\infty(D')}
 \le C_q\bigl(\|w\|_{L^q(D)}+\|f\|_{L^\infty(D)}
                         +\|f^{\partial}\|_{L^\infty(\Gamma_D)}\bigr).
\end{equation*}
This estimate applies to signed $w$. Together with
\eqref{eq:buffered-robin-estimate}, it gives local $C^1$ compactness
from local $L^1$ bounds and smooth bounded coefficients and data.
For positive solutions of the homogeneous equation with bounded
Robin coefficient, the boundary Harnack argument in
\cite[proof of Proposition 2.4]{LM} applies.
In particular, on local half-annuli,
\begin{equation}\label{eq:buffered-robin-harnack}
 \sup_{D'}w\le C\inf_{D'}w.
\end{equation}
The constants depend on a bound for the scaled Robin potential
and are uniform under multiplication of $w$ by a positive constant.

For a bounded positive semilinear solution, the initial boundary
H\"older estimate can be obtained from its local Neumann
representation. For a bounded face density $f$, split the kernel difference at
$|x-y|=2r$. The local kernel and outer gradient bounds yield
\[
 |G_\Gamma f(x)-G_\Gamma f(x')|
 \le C\|f\|_\infty\left(r+r\int_{2r}^1\frac{\dd s}{s}\right)
 \le C\|f\|_\infty r(1+|\log r|)
\]
whenever $|x-x'|\le r$. Volume potentials satisfy a stronger bound.
Hence $u\in C^\alpha$
on smaller patches for every $\alpha<1$.
If $0<c\le u\le C$, then $u^{n/(n-2)}\in C^\alpha$.
The boundary Schauder estimate first establishes $u\in C^{1,\alpha}$.
Iteration of \eqref{eq:buffered-robin-estimate} then yields uniform
$C^k$ bounds after Harnack normalization.

Let $u_j$ be weak solutions converging strongly in $H^1$ to
a smooth positive solution on a fixed patch. To obtain a uniform local
supremum bound, first consider a positive solution $u$ and test
with bounded truncations of
$\eta^2u^{2\beta-1}$, $\beta>1$. Ellipticity, Young's inequality,
and the trace Sobolev inequality imply
\begin{equation}\label{eq:critical-robin-test}
 \begin{aligned}
 &\|\nabla(\eta u^\beta)\|_{L^2(D)}^2
       +\|\eta u^\beta\|_{L^{\frac{2(n-1)}{n-2}}(\Gamma)}^2
       +\|\eta u^\beta\|_{L^{\frac{2n}{n-2}}(D)}^2\\
 &\quad\le C\beta^2\left\{
 \int_D(|\nabla\eta|^2+\eta^2)u^{2\beta}
 +\int_\Gamma(1+u^{\frac{2}{n-2}})(\eta u^\beta)^2\right\}.
 \end{aligned}
\end{equation}
For the positive-part equation before positivity has been proved,
use the signed test $\eta^2u|u|^{2\beta-2}$, with bounded
truncations. The same inequality holds with $u$ replaced by
$|u|$, because the nonlinear boundary contribution is bounded
by $C|u|^{\frac{2}{n-2}}(\eta|u|^\beta)^2$.
Here $C$ depends on the bounded lower-order geometric coefficients.
The $L^2$ trace inequality gives
\[
 \int_\Gamma(\eta u^\beta)^2
 \le\varepsilon\|\nabla(\eta u^\beta)\|_2^2
       +C_\varepsilon\|\eta u^\beta\|_2^2.
\]
This absorbs the bounded face term.
For one fixed $1<\beta\le n/(n-2)$, H\"older yields
\[
 \int_\Gamma u^{\frac{2}{n-2}}(\eta u^\beta)^2
 \le\|u\|_{L^{\frac{2(n-1)}{n-2}}(\supp\eta\cap\Gamma)}^{\frac{2}{n-2}}
                         \|\eta u^\beta\|_{L^{\frac{2(n-1)}{n-2}}(\Gamma)}^2.
\]
Choose the half-ball so that
\[
 C\beta^2\|u_j\|_{L^{2(n-1)/(n-2)}(\supp\eta\cap\Gamma)}^{2/(n-2)}<\frac12.
\]
The trace theorem transfers strong $H^1$ convergence to the
critical trace space. Since the limit is smooth, shrinking the
patch makes the displayed norm small uniformly for all late $j$.
Absorption in \eqref{eq:critical-robin-test} yields
\[
 u_j\in L^{\frac{2(n-1)\beta}{n-2}}(\Gamma'),
 \qquad u_j^{\frac2{n-2}}\in L^{(n-1)\beta}(\Gamma').
\]
For a coefficient $c\in L^p(\Gamma')$, $p>n-1$, interpolation yields
\[
 \int_{\Gamma'}|c|w^2
 \le\varepsilon\|w\|_{L^{\frac{2(n-1)}{n-2}}(\Gamma')}^2
                     +C_{\varepsilon,p,\|c\|_p}\|w\|_{L^2(\Gamma')}^2.
\]
Apply this inequality in the subsequent power tests on nested
half-balls. The ordinary Moser iteration now provides a uniform
local supremum bound. The preceding H\"older estimate ensures uniform
convergence. The positive limit then permits iteration of
\eqref{eq:buffered-robin-estimate}, yielding smooth convergence
on smaller patches. The
critical coefficient is absorbed at the first, fixed power.
Subsequent steps use the resulting subcritical coefficient.

\subsection{Comparison on expanding half-annuli}\label{app:annular-comparison}

Use the notation of Proposition~\ref{prop:annular-increment}, with
$1\le d<(n-2)/2$, $\sigma\to0$, $R\to\infty$ and $\rho R\to0$.
The constants are uniform on the compact bubble-parameter set.

\emph{The error estimate \eqref{eq:lower-error}.}
Let $G_0$ be the mixed Green kernel of $-\Delta_{g_\rho}$ on $\mathcal A_{\sigma,R}^+$ from
Appendix~\ref{app:green-kernels}, and put $G_0[F,F^\partial]:=G_0[F]+G_{0,\Gamma}[F^\partial]$.
After even reflection, the inner hemisphere becomes
$\partial\mathcal B_\sigma$. Here capacity is defined by the Dirichlet
energy of the reflected operator, with zero boundary values on
$\partial\mathcal B_R$. Choose a radial
$\zeta\in C_c^\infty(\mathcal B_2)$ equal to one near
$\overline{\mathcal B_1}$, and set $\zeta_\sigma(x)=\zeta(x/\sigma)$.
For $2\sigma<R$, uniform ellipticity and the variational definition
of capacity give
\[
 \begin{aligned}
 \operatorname{Cap}(\overline{\mathcal B_\sigma},\mathcal B_R)
 &\le C\int_{\mathcal B_R}|\nabla\zeta_\sigma|^2\dd x\\
 &=C\sigma^{n-2}\int_{\mathcal B_2}|\nabla\zeta|^2\dd y
 \le C\sigma^{n-2}.
 \end{aligned}
\]
The equilibrium measure on $\partial\mathcal B_\sigma$ has total
mass equal to this capacity. By \eqref{eq:mixed-green-upper}, its
Green potential is bounded by $C\sigma^{n-2}r^{2-n}$ for
$r\ge2\sigma$. On $\sigma\le r\le2\sigma$, the potential is at
most one by the maximum principle, which implies the same bound
after enlarging the constant.
Comparison with this potential controls the inner harmonic lift
of $w|_{\Sigma_\sigma}$. The maximum principle controls the outer
lift of $w|_{\Sigma_R}$. Thus
\[
 |h_\sigma(x)|\le C\sigma^{n-2}r^{2-n},\qquad
 |h_R(x)|\le CR^{2-n}.
\]
Since $d+1<n-2$, \eqref{eq:potential-powers} and
\eqref{eq:leading-annular-residual} imply
\[
 G_0[|E|,|E^\partial|]\le C\rho^{d+1}\phi_{d+1}.
\]

Fix $\eta>0$ small and $L>1$ large so that $\supp\chi\subset\{\eta<r<L\}$. Define the positive operator
\[
Sf := G_0\!\left[\frac{n-2}{4(n-1)}|R_{g_\rho}|f\right] +
G_{0,\Gamma}\!\left[\left(\frac{n-2}{2}|h_{g_\rho}|+c\mathbf
1_{\{r<\eta\}\cup\{r>L\}}\right)f\right].
\]
By the potential estimates in Appendix~\ref{app:green-kernels}, for each weight $\phi\in\{\phi_{d+1},r^{2-n},1,K_0\}$,
\[
 S\phi \le C\bigl(\eta+L^{-1}+\rho R+(\rho R)^2\bigr)\phi \le \tfrac12\phi
\]
for suitable fixed $\eta,L$ and large indices. Put $A:=\sup_{\eta<r<L,\,x_n=0}|w|$. The remaining
boundary source $c\mathbf 1_{\{\eta\le r\le L\}}|w|$ has potential bounded by $CAK_0$ (since
$K_0(r)=(1+r)^{2-n}$). Applying $(I-S)^{-1}$ (of norm at most $2$ on each weight) to
\eqref{eq:lower-remainder-equations} yields
\begin{equation}\label{eq:two-cap-preliminary}
 |w(x)| \le C\bigl(\rho^{d+1}\phi_{d+1}(x)+\sigma^{n-2}r^{2-n}+R^{2-n}\bigr) + CAK_0(x).
\end{equation}
Suppose along a subsequence that
\[
 \frac{A}{\rho^{d+1}+\sigma^{n-2}+R^{2-n}}\longrightarrow\infty.
\]
The local Robin estimates in Appendix~\ref{app:local-robin} imply
$w/A\to w_\infty$ in
$C^1_{\mathrm{loc}}(\overline{\R^n_+}\setminus\{0\})$, where
\[
 \begin{gathered}
 \begin{cases}
 -\Delta w_\infty=0&\text{in }\R^n_+\setminus\{0\},\\
 (\partial_\nu-nU_{\theta_*}^{\frac2{n-2}})w_\infty=0
     &\text{on }\Gamma\setminus\{0\},
 \end{cases}\\
 \sup_{\substack{\eta\le r\le L\\x_n=0}}|w_\infty|=1.
 \end{gathered}
\]
By \eqref{eq:two-cap-preliminary}, $|w_\infty|\le CK_0$.
Thus zero is removable and $w_\infty\in\mathcal D^{1,2}(\R^n_+)$.
Equation~\eqref{eq:ball-robin-form} implies
\[
 w_\infty=\sum_{\alpha=0}^{n-1}a_\alpha Z_\alpha(\theta_*),\qquad
 \sum_\alpha a_\alpha\int_\Gamma\chi Z_\alpha Z_\beta=0
       \quad\text{for }0\le\beta<n.
\]
The positive definiteness of the Gram matrix forces every
$a_\alpha$ to vanish, contradicting the normalization. Hence
\[
 A\le C(\rho^{d+1}+\sigma^{n-2}+R^{2-n}).
\]
Substitution into \eqref{eq:two-cap-preliminary}, using
$K_0\le\phi_{d+1},r^{2-n},1$, proves $|w|\le CW(r)$.

For $2\sigma\le r\le R/2$, rescaling $x=ry$ reduces
\eqref{eq:lower-remainder-equations} on $1/2\le|y|\le2$ to an
oblique system with uniformly bounded coefficients and data
\[
 r^2|E(ry)|+r|E^\partial(ry)|\le CW(r).
\]
By \eqref{eq:buffered-robin-estimate}, $r|\nabla w|\le CW(r)$.
Near $r=\sigma$ and $r=R$, \eqref{eq:annular-bounds} and the bounds
$W(\sigma)\ge c$, $W(R)\ge R^{2-n}$ yield \eqref{eq:lower-error}.

\emph{The comparison \eqref{eq:actual-virial-comparison}.}
Expand $\mathcal V_\rho(v)-\mathcal V_\rho(\Phi)$ using \eqref{eq:full-first-action-variation} and
$\log g=\sum_{k=d}^{2d}H_k+O(r^{2d+1})$. In the leading linear term $D_f J_1(H_d,U)[w]$,
integration by parts pairs $w$ with $(F_{H_d},F_{H_d}^\partial)$ of sizes $O((1+r)^{d-n})$ and
$O((1+r)^{d+1-n})$. For $r\ge 1$, \eqref{eq:lower-error} bounds both pairings by
\begin{align}\label{eq:lower-pairings}
 \rho^d\int_1^R r^{d-1}W(r)\dd r
 &\le C\rho^{2d+1}\int_1^R r^{2d+2-n}\dd r + C(\rho R)^d(\sigma^{n-2}+R^{2-n}) \notag\\
 &= o(\rho^{2d}) + C(\rho R)^d(\sigma^{n-2}+R^{2-n}),
\end{align}
Indeed, $n\ge2d+3$ implies
\[
 \rho\int_1^R r^{2d+2-n}\dd r
 \le C\rho(1+\log R)\longrightarrow0.
\]
The two capacity terms satisfy
\[
 \sigma^{n-2}\int_1^R r^{d+1-n}\dd r\le C\sigma^{n-2},\qquad
 R^{2-n}\int_1^R r^{d-1}\dd r\le CR^dR^{2-n}.
\]
On $r<1$, the integral is $O(\rho^d(\rho^{d+1}+\sigma^{n-2}+R^{2-n}))$.

The quadratic error in $w$ in the first metric variation is bounded by
\begin{equation}\label{eq:lower-quadratic-pairing}
 \rho^d\int_\sigma^R r^{n+d-3}W(r)^2\dd r \le o(\rho^{2d}) + C(\rho R)^d(\sigma^{n-2}+R^{2-n}),
\end{equation}
Indeed, the geometric summand contributes
\[
 \rho^{3d+2}\int_1^R r^{3d+3-n}\dd r=o(\rho^{2d}).
\]
A growing integral is bounded by $C\rho^{n-2}$ because $\rho R<1$ and
$n-2>2d$. The bounded and logarithmic cases lead to the same conclusion.
The two capacity summands satisfy
\[
 \begin{aligned}
 \rho^d\int_\sigma^R r^{n+d-3}
       (\sigma^{n-2}r^{2-n})^2\dd r
     &\le C(\rho\sigma)^d\sigma^{n-2},\\
 \rho^d\int_\sigma^R r^{n+d-3}R^{4-2n}\dd r
     &\le C(\rho R)^dR^{2-n}.
 \end{aligned}
\]
Cauchy's inequality controls the mixed products.
The spherical boundary terms from integration by parts are
$O(\rho^d\sigma^{n+d-1}+\rho^d R^{d+2-n})$, with the same capacity bound.

For the mixed term involving $\rho^d\psi$, divide by $\rho^{2d}$.
On $\kappa\le r\le K$, the normalized pairing tends to zero by local
$C^1$ convergence. On $r<\kappa$, smoothness of $\psi$ and
\eqref{eq:annular-bounds} bound it by $C\kappa^{n+d-2}$.
On $r>K$, the same annular bounds and the correction decay yield
\[
 C\int_K^\infty r^{2d+1-n}\dd r
 =\frac{C}{n-2-2d}K^{2d+2-n}\longrightarrow0.
\]
Taking the annulus limit first, then $\kappa\to0$ and $K\to\infty$,
yields $o(\rho^{2d})$.

Finally, the remaining geometric integrals satisfy
\[
 \begin{aligned}
 \rho^{k+d+1}\int_1^R r^{k+d+2-n}\dd r
     &=o(\rho^{2d})\quad\text{for }d<k\le2d,\\
 \rho^{2d+1}\left(1+\int_1^R r^{2d+2-n}\dd r\right)
     &=o(\rho^{2d}).
 \end{aligned}
\]
The higher-degree spherical boundary terms and quadratic errors obey
the same capacity bound because $(\rho R)^k\le(\rho R)^d$.
This proves \eqref{eq:actual-virial-comparison}.

\subsection{Scalar and boundary-tangent corrections on the ball}\label{app:ball-projections}

Use the tensors $h$, $\bar h$ and $k$ of
Lemma~\ref{lem:ball-hessian}.

\begin{proof}[Construction of the corrections]
The scalar variation satisfies
$DR_\delta[\sigma\delta]=-(n-1)\Delta\sigma$ and
$Dh_\delta[\sigma\delta]=(\partial_\nu\sigma-\sigma)/2$.
Thus $\sigma$ must solve
\begin{equation}\label{eq:scalar-completion}
 \begin{cases}
 \Delta\sigma=DR_\delta[h]/(n-1)&\text{in }\mathcal B_1,\\
 \partial_\nu\sigma-\sigma=-2Dh_\delta[h]&\text{on }\partial\mathcal B_1.
 \end{cases}
\end{equation}
The Robin kernel is $\Span\{x_1,\ldots,x_n\}$.
Fix the solution by $\int_{\partial\mathcal B_1}\sigma x_i=0$
for $1\le i\le n$. To verify compatibility with this kernel,
pair the curvature data with the coordinate functions.
Integration by parts yields, for smooth $h,f$,
\begin{equation}\label{eq:curvature-compatibility}
\begin{split}
 &\int_{\mathcal B_1}fDR_\delta[h]+2(n-1)\int_{\partial\mathcal B_1}fDh_\delta[h]
 \\
 &=\int_{\mathcal B_1}(\nabla\tr h-\diver h)\cdot\nabla f
 +\int_{\partial\mathcal B_1}\langle h(\nu,\cdot)^T,\nabla_Tf\rangle
 -\int_{\partial\mathcal B_1}(\tr_T h)f.
\end{split}
\end{equation}
For $f=x_i$, the first integral on the right equals
\[
 \int_{\partial\mathcal B_1}
 \bigl((\tr_T h)x_i
       -\langle h(\nu,\cdot)^T,\nabla_Tx_i\rangle\bigr).
\]
Here $h(\nu,\cdot)^T$ is the tangential one-form.
The three boundary terms therefore cancel, giving
\[
 \int_{\mathcal B_1}x_iDR_\delta[h]
       +2(n-1)\int_{\partial\mathcal B_1}x_iDh_\delta[h]=0.
\]
The right side of \eqref{eq:curvature-compatibility} is bounded by
$C\|h\|_{H^1}\|f\|_{H^1}$, with the tangential derivative interpreted
in $H^{1/2}$--$H^{-1/2}$ duality. Robin inversion, retaining the
invertible constant mode and fixing the coordinate kernel, provides
$\|\sigma\|_{H^1}\le C\|h\|_{H^1}$.

To obtain a divergence-free remainder with zero normal-tangential
entries, minimize $\|\bar h-\mathcal L_X\delta\|_2^2$ over
boundary-tangent fields orthogonal to rotations. For $X\cdot\nu=0$,
\[
 \int_{\mathcal B_1}|\mathcal L_X\delta|^2
 =2\int_{\mathcal B_1}|\nabla X|^2
  +2\int_{\mathcal B_1}(\diver X)^2
  -2\int_{\partial\mathcal B_1}|X|^2.
\]
For $X\cdot\nu=0$ and $\int_{\mathcal B_1}X\cdot Ax=0$ for every
$A^t=-A$, Korn's inequality implies
\[
 \|X\|_{H^1}\le C\|\mathcal L_X\delta\|_2.
\]
Indeed a Euclidean Killing field is $c+Ax$. Tangency forces $c=0$,
and orthogonality gives $A=0$.
Korn's inequality makes the minimizing sequence bounded in
$H^1$, and weak lower semicontinuity ensures the existence of a minimizer.
Its Euler equation is
$\int_{\mathcal B_1}k:\mathcal L_Y\delta=0$
for every boundary-tangent $Y$.
Integration by parts yields
\begin{equation}\label{eq:tangent-vector-system}
 \begin{cases}
 \Delta X+\nabla\diver X=\diver\bar h&\text{in }\mathcal B_1,\\
 X\cdot\nu=0&\text{on }\partial\mathcal B_1,\\
 (\mathcal L_X\delta)(\nu,T)=\bar h(\nu,T)&\text{on }\partial\mathcal B_1.
 \end{cases}
\end{equation}
Hence $\diver k=0$ and $k(\nu,T)=0$.
For the zero-data half-space system, tangential Fourier transformation
and the energy identity imply $\int|\mathcal L_X\delta|^2=0$.
A decaying nonzero-frequency Killing field vanishes.
The complementing condition and \cite{ADN} thus imply
\begin{equation*}
 \|X\|_{H^2}+\|k\|_{H^1}\le C\|\bar h\|_{H^1}.
\end{equation*}
Since $R_\delta=0$ and $h_\delta=1$,
\[
 DR_\delta[\mathcal L_X\delta]=0,\qquad
 Dh_\delta[\mathcal L_X\delta]=0.
\]
Thus $k$ satisfies both curvature constraints as well.

\end{proof}

\subsection{Adjoint function and local error estimate}\label{app:adjoint-estimates}

Use the normalized sequence, $W,w,Z$ and $I^\pm$ from the proof
of Proposition~\ref{prop:upper-criterion}, with $R=\delta/\eps$.
Thus $g_{ij}\to\delta_{ij}$ and $v\to U$ locally smoothly,
$U/2\le W\le3U/2$, and $w(0)=\nabla_{x'}w(0)=0$.
Volume and boundary integrals in this subsection are over
$\mathcal B_R^+$ and $B_R$, unless a domain is specified.

The rescaled puncture satisfies
\[
 \frac{|\Phi_k^{-1}(0)|}{\eps_k}
       \asymp\frac{|x_k|}{\eps_k}\longrightarrow\infty.
\]
Lemma~\ref{lem:no-source} gives the weak equation and
\eqref{eq:green-representation} on the full half-ball.
Let $G$ be the positive Green kernel of $(-L_g,B_g)$, with zero
Dirichlet value on $\Sigma_R$.
Its uniform bounds are established in Appendix~\ref{app:green-kernels}.

For the family $W_\lambda$, define $E_\lambda,E_\lambda^\partial$
by replacing $W$ with $W_\lambda$ in the residual equations.
All derivatives $\partial_\lambda E_\lambda$ and
$\partial_\lambda E_\lambda^\partial$ below are evaluated at $\lambda=1$.
Differentiating the residual equations, we obtain
\[
 \begin{cases}
 -L_gZ=-\partial_\lambda E_\lambda&\text{in }\mathcal B_R^+,\\
 (B_g-c_W)Z=-\partial_\lambda E_\lambda^\partial&\text{on }B_R.
 \end{cases}
\]

For volume or boundary data $f$, use the corresponding
measure in
\[
 [f]_{r_*}=\int_{r<2r_*}|f|
                         +\int_{r>2r_*}G(0,x)|f|.
\]
The Green potentials and local derivatives are measured by
$\norm{\mathcal R}_G$, while $\norm{\mathcal R}_*$ measures the
energy-dual sources.
\begin{equation}\label{eq:defined-core-residual}
 \begin{aligned}
 \norm{\mathcal R}_G={}&[E]_{r_*}+[E^{\partial}]_{r_*}
       +[\partial_\lambda E_\lambda]_{r_*}
       +[\partial_\lambda E_\lambda^{\partial}]_{r_*}\\
 &+\norm{E}_{C^1(\mathcal B_{2r_*}^+)}
       +\norm{\partial_\lambda E_\lambda}_{C^1(\mathcal B_{2r_*}^+)}\\
 &+\norm{E^{\partial}}_{C^2(B_{2r_*})}
       +\norm{\partial_\lambda E_\lambda^{\partial}}_{C^2(B_{2r_*})},\\
 \norm{\mathcal R}_*={}&\norm{E}_{L^{2n/(n+2)}(\mathcal B_R^+)}
                  +\norm{E^{\partial}}_{L^{2(n-1)/n}(B_R)}\\
 &+\norm{\partial_\lambda E_\lambda}_{L^{2n/(n+2)}(\mathcal B_R^+)}
       +\norm{\partial_\lambda E_\lambda^{\partial}}_{L^{2(n-1)/n}(B_R)}.
 \end{aligned}
\end{equation}
For local boundary estimates, take the norms on a fixed
neighborhood of $\overline{\mathcal B_{r_*}^+}$. The selected puncture
lies outside that neighborhood for sufficiently late rescalings.
The definitions can therefore use ordinary smooth local norms.

\begin{proof}[Construction of $Y$]
\emph{Step 1. Exterior boundary-value problem.}
For $\eta\in H^1(\{x_n>0,\ r_*/4<r<R\})$ vanishing on
$\Sigma_{r_*/4}\cup\Sigma_R$, the trace Sobolev inequality implies
\[
 \int_{B_R\setminus B_{r_*/4}}c_W\eta^2
\le C\|c_W\|_{L^{n-1}(B_R\setminus B_{r_*/4})}
             \|\nabla\eta\|_2^2
\le\frac{C}{r_*}\|\nabla\eta\|_2^2.
\]
Combining this with \eqref{eq:small-mixed-coercivity} yields,
for large $r_*$ and small $\delta$,
\begin{equation}\label{eq:exterior-coercivity}
 Q_g(\eta,\eta) - \int_{B_R\setminus B_{r_*/4}}c_W\eta^2
 \ge (c-C/r_*-C\delta-C\delta^2)\|\nabla\eta\|_2^2 \ge \frac{c}{2}\|\nabla\eta\|_2^2.
\end{equation}
By Lax--Milgram, there is a unique solution $\widetilde Y$ of
\[
\begin{cases}
 -L_g\widetilde Y=0&\text{in }\{x_n>0,\ r_*/4<r<R\},\\
 (B_g-c_W)\widetilde Y=0&\text{on }B_R\setminus\overline{B_{r_*/4}},\\
 \widetilde Y=Z&\text{on }\Sigma_{r_*/4},\\
 \widetilde Y=0&\text{on }\Sigma_R.
\end{cases}
\]
For large fixed $r_*$, the inner datum is negative because
\[
 Z\longrightarrow-Z_0
 =\frac{n-2}{2}(1-r^2)U^{\frac{n}{n-2}}<0
 \quad\text{on }\Sigma_{r_*/4}.
\]
Since $c_W\le C(1+r)^{-2}$, the potential estimates in Appendix~\ref{app:green-kernels} imply
\begin{equation}\label{eq:exterior-green-potential}
 \sup_x\int_{B_R\setminus B_{r_*/4}}G(x,z)c_W(z)\dd\sigma_g(z) \le \frac{C}{r_*} < \frac12.
\end{equation}
Hence $\widetilde Y<0$ and $-\widetilde Y\asymp G(0,\cdot)$ for $r\ge r_*/2$ (uniformly up to
$\Sigma_R$ with the boundary vanishing factor).

\emph{Step 2. Cutoff blending and outer flux.}
Choose $\chi\in C^\infty(\R)$ with $0\le\chi\le1$, equal to zero for
$r<r_*/2$ and to one for $r>3r_*/4$. Define
\[
 Y:=\begin{cases}
 Z,&r\le r_*/2,\\
 (1-\chi)Z+\chi\widetilde Y,&r_*/2<r<R.
 \end{cases}
\]
Both $Z$ and $\widetilde Y$ are negative on the blending annulus, so
$0\le Y_+\le C$ and $\supp Y_+\subset\overline{\mathcal B_{r_*}^+}$.
Since $Y=\widetilde Y$ near $\Sigma_R$,
\[
 \partial_{\nu_g}Y\asymp-\partial_{\nu_g}G(0,\cdot)>0.
\]
Let $\varphi$ be the rescaled positive gauge of
Appendix~\ref{app:green-kernels}, normalized by $\varphi(0)=1$.
Its bounds $c\le\varphi\le C$ are uniform for $0<\delta\le\delta_0$.
Green's identity yields
\[
 -\int_{\Sigma_R}\varphi\partial_{\nu_g}G(0,\cdot)\dd\sigma_g=1.
\]
Hence
\begin{equation}\label{eq:uniform-outer-flux}
 c_0\le\int_{\Sigma_R}\partial_{\nu_g}Y\dd\sigma_g\le C_0,
\end{equation}
with constants independent of $\delta$ and of the late sequence index. This proves \eqref{eq:outer-flux}.

\emph{Step 3. Defect and local error bounds.}
The difference $\widetilde Y-Z$ solves
\[
 (-L_g,B_g-c_W)(\widetilde Y-Z)
 =(\partial_\lambda E_\lambda,\partial_\lambda E_\lambda^\partial)
\]
on $\{r_*/4<r<R\}$ with datum $0$ on $\Sigma_{r_*/4}$ and $-Z$ on $\Sigma_R$.
Appendix~\ref{app:green-kernels} and \eqref{eq:buffered-robin-estimate} imply on
$\{r_*/2<r<3r_*/4\}$,
\begin{equation}\label{eq:blend-defect}
 \|\widetilde Y-Z\|_{C^2(\{r_*/2<r<3r_*/4\})} \le C\bigl(\|\mathcal R\|_G + R^{2-n}\bigr).
\end{equation}
Choose a lift $\ell$ supported in $R/2<r<2R$ with
$\ell|_{\Sigma_R}=-Z$ and $\|\nabla\ell\|_2^2\le CR^{2-n}$.
Test the difference equation with $\eta=\widetilde Y-Z-\ell$ and
use \eqref{eq:exterior-coercivity} to obtain
\[
 \|\nabla(\widetilde Y-Z)\|_{L^2(\{r_*/4<r<R\})}
 \le C\bigl(\|\mathcal R\|_*+R^{(2-n)/2}\bigr).
\]
Extend $Y-Z$ by $\ell$ outside $\mathcal B_R^+$.
Poincar\'e on $\{r_*/4<r<r_*\}$, together with the Sobolev and trace
inequalities, implies
\[
 |\mathcal R[Y-Z]| \le C\|\mathcal R\|_*\bigl(\|\mathcal R\|_*+R^{(2-n)/2}\bigr) \le C\|\mathcal R\|_*^2 + CR^{2-n},
\]
proving \eqref{eq:patch-comparison}.
Finally, $L_gY=\partial_\lambda E_\lambda$ on $\{r\le r_*/2\}$.
For $r>r_*/2$, the blending formula and $\partial_{\nu_g}\chi=0$ on $B_R$ give
\[
 \begin{aligned}
 L_gY={}&(1-\chi)\partial_\lambda E_\lambda
       +2\langle\nabla\chi,\nabla(\widetilde Y-Z)\rangle_g\\
       &+(\Delta_g\chi)(\widetilde Y-Z),\\
 (B_g-c_W)Y={}&-(1-\chi)\partial_\lambda E_\lambda^\partial.
 \end{aligned}
\]
These defects vanish for $r\ge r_*$.
Estimate~\eqref{eq:blend-defect} therefore yields
\begin{equation}\label{eq:adjoint-local-errors}
 \|L_gY\|_{L^\infty(\mathcal B_R^+)}
 +\|(B_g-c_W)Y\|_{L^\infty(B_R)}
 \le C(\norm{\mathcal R}_G+R^{2-n}).
\end{equation}
\end{proof}

\begin{lemma}\label{lem:core}
Let $w$ be the error in \eqref{eq:error-equations}, with
$w(0)=0$, $\nabla_{x'}w(0)=0$, and
$|w|\le C\eps^{(n-2)/2}$ on $\Sigma_R$.
Then \eqref{eq:core-estimate} holds for each sufficiently large
fixed $r_*$, with a constant independent of the late sequence index.
\end{lemma}
\begin{proof}
Since $-Y\asymp G(0,\cdot)$ for $r\ge r_*$,
\[
 \int_{B_R\setminus B_{r_*}}G(0,z)N_W(w)\dd\sigma_g\le CI^-.
\]
Write $N_W(w)=\alpha_\eps w$ on $B_{r_*}$.
Then $\|\alpha_\eps\|_{C^j(B_{r_*})}\to0$ and
\[
 c_W+\alpha_\eps
 =n\int_0^1(W+sw)^{\frac2{n-2}}\dd s\ge0.
\]
Define
\[
 \begin{aligned}
 S_\eps f(x)&:=\int_{B_R\setminus B_{r_*}}G(x,z)c_W(z)f(z)\dd\sigma_g(z),\\
 T_\eps\eta(x)&:=\int_{B_{r_*}}G(x,z)(c_W+\alpha_\eps)(z)\eta(z)\dd\sigma_g(z).
 \end{aligned}
\]
By \eqref{eq:exterior-green-potential} and \eqref{eq:green-three},
for large $r_*$,
\[
 \|S_\eps\|_{L^\infty\to L^\infty}\le\frac12,\qquad
 S_\eps G(\cdot,y)\le\frac12G(\cdot,y).
\]
The kernel also satisfies
\[
 \int_{B_{r_*}}G(x,y)\dd\sigma_g(x)
 \le C_{r_*}\begin{cases}
 1,&|y|\le2r_*,\\
 G(0,y),&|y|>2r_*.
 \end{cases}
\]
Symmetry and Fubini therefore imply, for every volume or boundary source $f$,
\begin{equation}\label{eq:core-source-L1}
 \|S_\eps^j G[f]\|_{L^1(B_{r_*})} \le C_{r_*}2^{-j}[f]_{r_*},
\end{equation}
and similarly in local volume $L^1$.

By \eqref{eq:green-representation},
\[
 w=S_\eps w+T_\eps(w|_{B_{r_*}})+h_\eps,
\]
where
\[
 \|h_\eps\|_{L^1(\mathcal B_{2r_*}^+)}
 +\|h_\eps\|_{L^1(B_{2r_*})}
 \le C\bigl(\eps^{(n-2)/2}+I^-+\|\mathcal R\|_G\bigr).
\]
For each fixed $\eps$, $v^{n/(n-2)}\in L^1(B_R)$ and hence
$\int_{B_R}c_W|w|<\infty$. Estimate~\eqref{eq:core-source-L1} implies
\[
 \|S_\eps^Jw\|_{L^1(B_{r_*})}\longrightarrow0
 \quad\text{as }J\to\infty.
\]
Summing the Neumann series, we obtain
\begin{equation}\label{eq:core-fredholm-equation}
 \begin{gathered}
 w=\mathcal K_\eps(w|_{B_{r_*}})+\widetilde h_\eps,\\
 \|\widetilde h_\eps\|_{L^1(B_{r_*})}
 \le C\bigl(\eps^{\frac{n-2}{2}}+I^-+\|\mathcal R\|_G\bigr).
 \end{gathered}
\end{equation}
where
\[
 \mathcal K_\eps:=\sum_{j=0}^\infty S_\eps^jT_\eps.
\]
Removing the singularity $|x-z|<h$ and truncating at $J$ iterates
produces $\mathcal K_\eps^{h,J}$ with equicontinuous kernels and
\[
 \sup_\eps\|\mathcal K_\eps-\mathcal K_\eps^{h,J}\|_{L^1\to L^1}
 \le C(h+2^{-J}).
\]
Arzel\`a--Ascoli yields collective compactness of $\mathcal K_\eps$
on $L^1(B_{r_*})$.
For $|x|>2r_*$, \eqref{eq:green-three} also implies
\begin{equation}\label{eq:linear-extension-tail}
 |\mathcal K_\eps\eta(x)| \le C\|\eta\|_{L^1(B_{r_*})}G(0,x) \le C\|\eta\|_{L^1(B_{r_*})}|x|^{2-n}.
\end{equation}

If \eqref{eq:core-estimate} fails, set $m_\eps=\|w\|_{L^1(B_{r_*})}$ and
$\widehat w_\eps=w/m_\eps$. Then
\[
 \frac{\eps^{(n-2)/2}+I^-+\|\mathcal R\|_G}{m_\eps}\longrightarrow0.
\]
Equation~\eqref{eq:core-fredholm-equation} and collective compactness give
\[
 \widehat w_\eps\longrightarrow w_\infty
 \quad\text{in }L^1(B_{r_*}),\qquad
 \|w_\infty\|_{L^1(B_{r_*})}=1,
\]
and $\widehat w_\eps-\mathcal K_\eps\widehat w_\eps\to0$ in local $L^1$.
Passing to the limit in \eqref{eq:error-equations} and applying the
local Robin estimates in Appendix~\ref{app:local-robin} yields
$\widehat w_\eps\to w_\infty$ in $C^1(\overline{\mathcal B_{r_*}^+})$, where
\begin{equation}\label{eq:core-jacobi-limit}
 \begin{cases}
 \Delta w_\infty=0&\text{in }\R^n_+,\\
 \partial_\nu w_\infty=nU^{\frac2{n-2}}w_\infty&\text{on }\Gamma,\\
 w_\infty(0)=0,\quad\nabla_{x'}w_\infty(0)=0.
 \end{cases}
\end{equation}
By \eqref{eq:linear-extension-tail} and \eqref{eq:buffered-robin-estimate},
\[
 |w_\infty|+r|\nabla w_\infty|\le Cr^{2-n}
 \qquad\text{for }r\ge1.
\]
Thus $w_\infty\in\mathcal D^{1,2}(\R^n_+)$. By
\eqref{eq:ball-robin-form}, $w_\infty=\sum_{\alpha=0}^{n-1}a_\alpha Z_\alpha$.
The conditions at zero in \eqref{eq:core-jacobi-limit} force all
coefficients to vanish, contradicting the normalization.
\end{proof}

\subsection{Taylor estimates for corrected bubbles}\label{app:taylor-estimates}

\begin{lemma}\label{lem:taylor-budget}
Let $W_\lambda$ be the corrected family in the proof of
Proposition~\ref{prop:upper-criterion}, with the orders $k,L$ fixed
there, $Z=-\partial_\lambda W_\lambda|_{\lambda=1}$ and
$R=\delta/\eps$. With the residual norms of
\eqref{eq:defined-core-residual}, for sufficiently small fixed
$\delta$ its residuals satisfy
\begin{equation}\label{eq:budget-residuals}
 \norm{\mathcal R}_G\le C(|H|_\eps^2+\eps^{\frac{n-2}{2}}).
\end{equation}
In the energy-dual norm,
\[
 \norm{\mathcal R}_*^2\le C\eps^2|H|_\eps^2+
 C\eps^{n-2}(1+|\log\eps|)^2=o(|H|_\eps^2+\eps^{\frac{n-2}{2}}).
\]
The dilation projection satisfies
\begin{equation}\label{eq:budget-projection}
 \mathcal R[Z]=\cQ\left(\sum_{d=1}^k\eps^dH_d\right)
                 +o(|H|_\eps^2+\eps^{\frac{n-2}{2}}).
\end{equation}
The bounds are uniform at the selected centers and under any fixed
number of bubble-parameter derivatives.
\end{lemma}
\begin{proof}
Since $\frac12 U\le W\le\frac32 U$, expand the metric operators and normalized boundary
nonlinearity through a fixed weight $M>2L+n+4$. Each residual term $(E_\ell,E_\ell^\partial)$ of
weight $\ell$ satisfies
\begin{equation}\label{eq:residual-monomial}
 |E_\ell| \le C\eps^\ell|c_\ell|(1+r)^{\ell-n},\qquad
 |E_\ell^\partial| \le C\eps^\ell|c_\ell|(1+r)^{\ell+1-n},
\end{equation}
Here $c_\ell$ is linear in $H_\ell$ or is a nonlinear coefficient satisfying
\[
 |c_\ell|\le C\prod_{a=1}^q|H_{d_a}|,\qquad
 q\ge2,\qquad\sum_{a=1}^q d_a=\ell.
\]
The Taylor remainder has weight $M+1$, using $\eps r\le\delta$. By
\eqref{eq:jet-correction}, all terms linear in $H_\ell$ for $\ell\le L$ vanish.

\emph{Step 1. Proof of \eqref{eq:budget-residuals}.}
Integrating against $G(0,x)\le Cr^{2-n}$, we obtain
\begin{equation}\label{eq:weighted-radial-budget}
 \eps^\ell|c_\ell|\left(1+\int_1^R r^{\ell+1-n}\dd r\right)
 \le C\eps^\ell|c_\ell|\begin{cases}1,&\ell<n-2,\\ 1+\log R,&\ell=n-2,\\ R^{\ell+2-n},&\ell>n-2.\end{cases}
\end{equation}
For $\ell\ge n-2$, the bound is
\[
 O\bigl(\eps^{n-2}(1+|\log\eps|)\bigr)
 =o(\eps^{(n-2)/2}).
\]
Surviving linear terms have $\ell\ge L+1\ge(n-2)/2$.
For quadratic terms ($q=2$), if $i,j\le k$, then
\[
 \eps^{i+j}|H_i||H_j|\le|H|_\eps^2.
\]
If $i\le k<j$, then $j>(n-2)/4$ and
\begin{equation}\label{eq:cross-budget}
 \eps^{i+j}|H_i||H_j|
 \le C\eps^{j-\frac{n-2}{4}}
       \bigl(|H|_\eps^2+\eps^{\frac{n-2}{2}}\bigr)
 =o\bigl(|H|_\eps^2+\eps^{\frac{n-2}{2}}\bigr).
\end{equation}
If $i,j>k$, then $i+j\ge2(k+1)>(n-2)/2$.
For $q\ge3$, order the degrees as $d_1\le\dots\le d_q$.
When $\sum d_a\le(n-2)/2$, we have $d_1,d_2\le k$ and
\[
 \eps^{\sum d_a}\prod|H_{d_a}|\le C\eps|H|_\eps^2.
\]
Higher weights contribute
$O(\eps|H|_\eps^2)+o(\eps^{(n-2)/2})$.
The same bounds apply to $\partial_\lambda E_\lambda$,
$\partial_\lambda E_\lambda^\partial$ and the local
$C^1,C^2$ norms, proving the $\norm{\mathcal R}_G$ bound.

For $\norm{\mathcal R}_*$, the volume and boundary dual exponents satisfy
\[
 \ell-n+\frac{n+2}{2}
 =\ell+1-n+\frac n2
 =\ell-\frac{n-2}{2}.
\]
Consequently,
\begin{equation}\label{eq:dual-radial-budget}
 \|E_\ell\|_{L^{\frac{2n}{n+2}}(\mathcal B_R^+)}^2 + \|E_\ell^\partial\|_{L^{\frac{2(n-1)}{n}}(B_R)}^2
\le C\eps^{2\ell}|c_\ell|^2\begin{cases}1,&\ell<\frac{n-2}{2},\\ (1+\log
R)^2,&\ell=\frac{n-2}{2},\\ 1+R^{2\ell-(n-2)},&\ell>\frac{n-2}{2}.\end{cases}
\end{equation}
Surviving linear terms have $\ell\ge(n-2)/2$ and contribute
$O(\eps^{n-2}(1+|\log\eps|)^2)$.
Any nonlinear term with $\ell<(n-2)/2$ has a factor $H_i$ with
$i\le k$ and $\ell-i\ge1$, so
\[
 \eps^{2\ell}|c_\ell|^2\le C\eps^2|H|_\eps^2.
\]
Summing proves the $\norm{\mathcal R}_*^2$ bound.

\emph{Step 2. Proof of \eqref{eq:budget-projection}.}
Integrating by parts yields
\[
 -\left.\partial_\lambda
 J_{g,\mathcal B_R^+}(W_\lambda)\right|_{\lambda=1}
 =\mathcal R[Z]+\int_{\Sigma_R}Z\partial_{\nu_g}W\dd\sigma_g.
\]
For $d<n-2$, \eqref{eq:full-first-action-variation} and
$J_1(\chi_TH_d,U)=0$ give
\[
 |J_1(H_d,U)-J_1(\chi_TH_d,U)|
 \le C|H_d|T^{d+2-n}\longrightarrow0
 \quad\text{as }T\to\infty.
\]
Thus $J_1(H_d,U)=0$.
Terms with $d\ge n-2$ contribute $o(\eps^{(n-2)/2})$ by
\eqref{eq:weighted-radial-budget}.
The quadratic contribution is
\[
 -\left.\lambda\partial_\lambda
       [\lambda^{i+j}\mathcal C(H_i,H_j)]\right|_{\lambda=1}
 =\cQ(H_i,H_j).
\]
Terms with $\max(i,j)>k$ are controlled by \eqref{eq:cross-budget}
and \eqref{eq:weighted-radial-budget}.
On $\Sigma_R$,
\[
 |W|+r|\nabla W|+|Z|+r|\nabla Z|\le Cr^{2-n}.
\]
The spherical boundary term is therefore
\[
 O(R^{2-n})=O_\delta(\eps^{n-2})=o(\eps^{(n-2)/2}),
\]
which proves \eqref{eq:budget-projection}.
\end{proof}

\emph{The eighth-order terms in dimension eighteen.}
Here $n=18$, $L=7$ and $R=\delta/\eps$. The determined products of
lower-order coefficients are included among the nonlinear monomials.
Since $|H_4|\le C$,
\[
 |H|_\eps^2=|H|_{\eps,3}^2+\eps^8|H_4|^2
 \le C\bigl(|H|_{\eps,3}^2+\eps^8\bigr).
\]
Thus \eqref{eq:budget-residuals} implies the $\|\mathcal R\|_G$ bound
in \eqref{eq:eighteen-residual-budget}.
Linear terms with $\ell\le7$ vanish by \eqref{eq:jet-correction}.
Any nonlinear monomial with ordered degrees
$1\le d_1\le\dots\le d_q$, $q\ge2$ and $\sum d_i\le7$ has
$d_1\le3$. Estimate~\eqref{eq:dual-radial-budget} implies
\[
 \eps^{2\sum d_i}\prod|H_{d_i}|^2
 \le C\eps^2|H|_{\eps,3}^2.
\]
At weight eight,
\[
 |E_8|\le C\eps^8(1+r)^{-10},\qquad
 |E_8^\partial|\le C\eps^8(1+r)^{-9}.
\]
The corresponding dual norms satisfy
\[
 \begin{aligned}
 \|E_8\|_{L^{9/5}}^2
     &\le C\eps^{16}(1+\log R)^{10/9},\\
 \|E_8^\partial\|_{L^{17/9}}^2
     &\le C\eps^{16}(1+\log R)^{18/17}.
 \end{aligned}
\]
The Green-weighted radial exponent is $8+1-18=-9$.
Higher weights follow from \eqref{eq:weighted-radial-budget} and
\eqref{eq:dual-radial-budget}, proving \eqref{eq:eighteen-residual-budget}.

For the projection, the remaining quadratic pairs satisfy
$1\le i\le j$, $j\ge4$ and $i+j\le8$.
Those with $j=4$ contribute $M_\eps$ and $\eps^8\cQ(H_4)$.
Those with $j\ge5$ have $i\le3$ and contribute at most
\[
 C\eps^5|H|_{\eps,3}
 \le C\eps\bigl(|H|_{\eps,3}^2+\eps^8\bigr).
\]
Terms with $q\ge3$ and $\sum d_i\le8$ have $d_1,d_2\le3$ and
contribute $O(\eps|H|_{\eps,3}^2)$.
The independent linear term $H_8$ satisfies $J_1(H_8,U)=0$ with
unscaled cutoff error $O(R^{-8})$, so its contribution is
\[
 \eps^8O(R^{-8})=O_\delta(\eps^{16}).
\]
The Taylor remainder is $O(\eps^9)$ and the spherical boundary terms
are $O_\delta(\eps^{16})$.
Combining these bounds with \eqref{eq:patch-comparison} proves
\eqref{eq:eighteen-projection}.

\section{Geometric quadratic forms}

\subsection{Integration and projection formulas}\label{app:integration-formulas}

Write $\Gamma(z)$ and $\mathrm B(z,w)$ for Euler's gamma and beta functions.
Euler's beta integral \cite[Section 1.1]{AAR} evaluates the
direct metric terms. Gegenbauer projections below evaluate the
Robin correction and boundary shape terms.
For multi-indices $\alpha$ and $p>(n-1)/2+|\alpha|$, the tangential
moments used throughout the paper are
\begin{equation}\label{eq:beta-rule}
 \int_{\R^{n-1}}{x'}^{2\alpha}(|{x'}|^2+(1+{x_n})^2)^{-p}\dd {x'}
 =\frac{\prod_i\Gamma(\alpha_i+\tfrac12)
          \Gamma(p-\tfrac{n-1}{2}-|\alpha|)}{\Gamma(p)}
    (1+{x_n})^{n-1+2|\alpha|-2p}.
\end{equation}
The remaining normal integration uses
\[
 \int_0^\infty {x_n}^r(1+{x_n})^{-p}\dd {x_n}
 =\frac{\Gamma(r+1)\Gamma(p-r-1)}{\Gamma(p)}\quad\text{for }r>-1,\ p>r+1.
\]
Odd tangential moments vanish by reflection.

We use the Gegenbauer normalization and orthogonality of
\cite[Section 4.7]{Szego}. Write the factorials
\[
 (z)_k=\prod_{j=0}^{k-1}(z+j),\qquad
 k^{\underline r}=\prod_{j=0}^{r-1}(k-j),
\]
for $k,r\in\mathbb N_0$, with empty products equal to one.
For tangential harmonic degree $\ell$, put
$\lambda=\ell+(n-2)/2$ and $A=\lambda+1/2$.
For $A>\nu$, the projection and norm formulas are
\begin{equation}\label{eq:gegenbauer-beta}
 \begin{aligned}
 I_\nu(k)&=\int_{-1}^1(1-t^2)^{\lambda-1/2}
                 (1+t)^{-\nu}C_k^\lambda(t)\dd t\\
 &=(-1)^k2^{2\lambda-\nu}\mathrm B(A-\nu,A)
       \frac{(2\lambda)_k(\nu)_k}{k!(2\lambda-\nu+1)_k},\\
 h_k&=\mathrm B(1/2,A)\frac{(2\lambda)_k}{k!}\frac\lambda{k+\lambda}.
 \end{aligned}
\end{equation}
Rodrigues' formula \cite[equation (4.7.12)]{Szego}, followed by $k$
integrations by parts, reduces $I_\nu(k)$ to Euler's beta integral.
The endpoint terms vanish since $A>\nu$. The expression for $h_k$
is the norm in \cite[equation (4.7.15)]{Szego}.
For $r\in\mathbb N_0$, $b>0$ and $c>b+r+1$, Gauss summation
and its differentiated form \cite[Sections 2.1--2.2]{AAR}
give
\begin{equation}\label{eq:factorial-moment}
 \sum_{k=0}^\infty k^{\underline r}\frac{(b)_k}{(c)_k}
 =\frac{(c-1)r!(b)_r}
        {\prod_{j=0}^r(c-b-1-r+j)}.
\end{equation}

Taking $b=1$, $c=v+1$ and using
$\mathrm B(k+1,v)=v^{-1}(1)_k/(v+1)_k$, together with
$k^{\underline r}(k+\alpha)=k^{\underline{r+1}}+(r+\alpha)k^{\underline r}$,
also yields, for $v>r+2$,
\begin{equation*}
 \sum_{k=0}^\infty k^{\underline r}(k+\alpha)\mathrm B(k+1,v)
 =r!\{(r+\alpha)\mathrm B(r+1,v-r-1)
                    +(r+1)\mathrm B(r+2,v-r-2)\}.
\end{equation*}
These formulas are used for the metric Robin pairings and the
Euclidean shape energy. After cancellation of the common positive
factor, the entries below are rational functions of the dimension.
All dimensions in the following finite comparisons are integers.
We use the recurrence $\mathrm B(x+1,y)=x\mathrm B(x,y)/(x+y)$
together with the displayed integral formulas.

\subsection{The quadratic and cubic metric forms}\label{app:general-forms}

Write $D=|x'|^2+(1+x_n)^2$ and $I_n=\int U^2$.
The harmonic Robin solution $h_q$ is defined in
\eqref{eq:harmonic-robin}.

For \eqref{eq:first-correction-upper}, with $n>5$, put
$Q_A={x'}^tAx'$ and choose the interior particular solution
\[
 \phi=\frac{n-2}{2}(x_n-1)Q_AD^{-n/2}.
\]
It has the source of $\psi_{H_1}$, with Robin residual $-q_A$,
where $q_A=\left.\frac{n-2}{2}Q_AD^{-n/2}\right|_\Gamma$.
Indeed,
\[
 \begin{cases}
 -\Delta\phi=2n(n-2)x_nQ_AD^{-(n+2)/2}=F_{H_1}
                                            &\text{in }\R^n_+,\\
 (\partial_\nu-nU^{2/(n-2)})\phi=-q_A&\text{on }\Gamma.
 \end{cases}
\]
Thus $\psi_{H_1}=\phi+h_{q_A}$.
The tangential harmonic degree is two, so all compactified source
degrees are at least two.
The Robin inverse formula \eqref{eq:robin-tail-control}, including
all modes, implies
\[
 Q(h_{q_A})
 =\sum_{j\ge0}\frac{\int_\Gamma D|\Pi_jq_A|^2}{2(j+1)}
 \le\frac12\int_\Gamma Dq_A^2.
\]
Testing \eqref{eq:harmonic-robin} with $\phi$ yields
$Q(\phi,h_{q_A})=\int_\Gamma q_A\phi$.
Since the Robin residual of $\phi$ is $-q_A$,
\[
 Q(\phi)=\int F_{H_1}\phi-\int_\Gamma q_A\phi.
\]
Consequently
\[
 \begin{aligned}
 Q(\psi_{H_1})
 &=Q(\phi)+2Q(\phi,h_{q_A})+Q(h_{q_A})\\
 &=\int F_{H_1}\phi+\int_\Gamma q_A\phi+Q(h_{q_A}).
 \end{aligned}
\]
The integration formula \eqref{eq:beta-rule} yields
\[
 \int F_{H_1}\phi=\frac{6-n}{4(n-1)}I_n|A|^2.
\]
On $\Gamma$, $\phi=-q_A$, so the boundary pairing is
\[
 \int_\Gamma q_A\phi=-\int_\Gamma q_A^2
       =-\frac{(n-2)(n-4)}{8(n-1)}I_n|A|^2.
\]
The remaining inverse bound is
\[
 \frac12\int_\Gamma Dq_A^2
       =\frac{(n-2)(n-4)}{8(n-5)}I_n|A|^2.
\]
Adding these three values proves \eqref{eq:first-correction-upper}.

For $14\le n\le18$, the representation table in the proof of
Proposition~\ref{prop:metric-coercivity} leaves five mixed matrices.
We compute their entries and control the Robin remainders by
\eqref{eq:spectral-tail-matrix}.

For the $[2]$ representation, use the same trace-free $A$ and define
\[
 L_A=(A{x'})\otimes{x'}+{x'}\otimes(A{x'})-\frac2{n-1}Q_AI.
\]
The five coefficient tensors are
\[
 (h_1,h_2,h_3,h_4,h_5)
 =(-2x_nA,-2x_n^2A,-2x_n^3A,-2x_n|x'|^2A,-2x_nL_A).
\]

Write $\psi_i=\psi_{h_i}$ for their linear corrections.
The homogeneous degrees are $d=(1,2,3,3,3)$. Define
$M^0_{ij}=-(d_i+d_j)J_2(h_i,h_j)/(I_n|A|^2)$ using
\eqref{eq:direct-form}.
Interior particular solutions are
$\phi_i=Q_A[P_i({x_n})D^{-n/2}+R_i({x_n})D^{-{\frac{n-2}{2}}}]$, with
\[
\begin{array}{c|c|c}
 i&P_i({x_n})&R_i({x_n})\\\hline
 1&(n-2)({x_n}-1)/2&0\\
 2&(n-2)({x_n}^2-{x_n}+1)/3&1/6\\
 3&(n-2)({x_n}-1)({x_n}^2+1)/4&({x_n}-1)/4\\
 4&-(n-2)({x_n}+1)^2(3{x_n}-1)/12&((4n-11){x_n}-2n+3)/12\\
 5&-(n-2)^2({x_n}+1)^2(3{x_n}-1)/(6(n-1))&
 ((2n^2-11n+16){x_n}-n^2+3n-3)/(6(n-1))
\end{array}
\]
For a tangential homogeneous harmonic polynomial $Q_\ell$ of
degree $\ell$, these formulas are checked using
\[
 \Delta[Q_\ell P({x_n})D^{-p}]
 =Q_\ell\{P''D^{-p}-4p({x_n}+1)P'D^{-p-1}
              +2p(2p+2-n-2\ell)PD^{-p-1}\}.
\]
Set $q_0=Q_AD^{-n/2}$, $q_1=Q_AD^{-{\frac{n-2}{2}}}$ on the face and
$u_i=P_i'(0)+2R_i(0)$, $v_i=R_i'(0)$.
The Robin residual is $-u_iq_0-v_iq_1$, so
$\psi_i=\phi_i+h_{u_iq_0+v_iq_1}$, where $h_q$ solves
\eqref{eq:harmonic-robin}. The integration formulas in Appendix~\ref{app:integration-formulas} evaluate the moment as
\[
 K_p=|A|^{-2}\int_{{x_n}=0}Q_A^2D^{-p}
       =\frac{\pi^{(n-1)/2}\Gamma(p-(n-1)/2-2)}{2\Gamma(p)}.
\]

Here is a finite lower bound for the complete inverse contribution.
For the two-dimensional source space spanned by $q_0,q_1$,
use indices $i,j\in\{0,1\}$ and mode index $k\ge0$.
Expand
$C_k^{(n+2)/2}(1-2x)=\sum_rb_{kr}x^r$ and define
\[
 a_{ik}=4\sum_rb_{kr}K_{n+1-i+r},\qquad
 N_k=16\sum_{r,t}b_{kr}b_{kt}K_{n+3+r+t}.
\]
The projected Gram matrix has entries $(P_k)_{ij}=a_{ik}a_{jk}/N_k$.
The full source Gram matrix has entries $T_{ij}=K_{n-1-i-j}$.

The source mode has total spherical degree $k+2$, with inverse
eigenvalue $1/[2(k+1)]$. Put
$H_{\le K}=\sum_{k=0}^KP_k/[2(k+1)]$.
Orthogonality of the projections and the decreasing inverse
eigenvalues imply
\begin{equation}\label{eq:spectral-tail-matrix}
 H=\sum_{k=0}^\infty\frac{P_k}{2(k+1)},\qquad
 0\le H-H_{\le K}\le
 R_*:=\frac{T-\sum_{k=0}^KP_k}{2(K+2)}.
\end{equation}
Retain the first eight modes, so $K=7$ throughout.
Let $V$ have rows
$(u_i,v_i)$, and let
\[
 G^0_{ij}=|A|^{-2}
  \left(\int F_{h_i}\phi_j+\int(u_jq_0+v_jq_1)\phi_i\right),
 \qquad \ell=(2u_1,u_2,0,0,0)^t.
\]
Green symmetry gives
$Q(\psi_i,\psi_j)/|A|^2=G^0_{ij}+(VHV^t)_{ij}$.
Apply the weighted square identity \eqref{eq:spectral-completion}
to the nonnegative remainder $H-H_{\le K}$.
Since $v_1=v_2=0$, dropping the positive square produces the
following lower matrix for the full dilation form,
\begin{equation}\label{eq:five-copy-lower}
 M_{ij}=M^0_{ij}
 +\frac{d_i+d_j}{2I_n}(G^0+VH_{\le K}V^t)_{ij}
 -\frac{(R_*)_{00}}{12I_n}\ell_i\ell_j.
\end{equation}
For a coefficient vector $x=(x_1,\ldots,x_5)^t$, the
degree-one and degree-two source vectors are
$u_1x_1(1,0)^t$ and $u_2x_2(1,0)^t$, since $v_1=v_2=0$.
The negative square in \eqref{eq:spectral-completion} is therefore
\[
 -\frac1{12}(H-H_{\le K})_{00}(2u_1x_1+u_2x_2)^2
 \ge-\frac1{12}(R_*)_{00}(\ell^tx)^2.
\]
The remaining spectral contribution is the nonnegative square
in that identity. This proves the rank-one remainder bound
in \eqref{eq:five-copy-lower}.

For the lower matrix in \eqref{eq:five-copy-lower}, put
$p_j=\det M_{[1,j]}/\det M_{[1,j-1]}$, where
$M_{[1,j]}$ is the leading $j\times j$ submatrix and
$\det M_{[1,0]}=1$. Substitution of the integral formulas yields, for
each integer $14\le n\le18$, the componentwise bounds
\begin{equation}\label{eq:metric-schur-bounds}
 (p_1,p_2,p_3,p_4,p_5)>
 \left(\frac16,\frac1{120},\frac1{1000},\frac1{25},\frac16\right).
\end{equation}
The leading principal minors are therefore positive. Sylvester's
criterion establishes $M>0$ and positivity of the complete five-copy form.

The following representatives, together with
\eqref{eq:direct-form} and the integration formulas, specify
all the smaller matrices. For the Weyl block take
$K_{12}=1$, $K_{13}=-1$, $K_{24}=-1$, $K_{34}=1$ and zero
otherwise, with
$W_{ii}=\sum_{j\ne i}K_{ij}{x'}_j^2$ and
$W_{ij}=-K_{ij}{x'}_i{x'}_j$ for $i\ne j$. The basis is $(W,{x_n}W)$,
with degrees $(2,3)$, and $F_W=F_{x_nW}=0$.

For the hook block, let $T_{ij,k}$ be symmetric in $i,j$, with
zero contractions and zero cyclic sum. Take
$T_{12,3}=T_{21,3}=1$, $T_{13,2}=T_{31,2}=-1$ and zero
otherwise. Set
\[
 P_{ij}=T_{ij,k}{x'}_k,\qquad
 C_T=|{x'}|^2P-{x'}\otimes(P{x'})-(P{x'})\otimes {x'}.
\]
The basis $(-2{x_n}P,-2{x_n}^2P,C_T)$ has degrees $(2,3,3)$.
Their sources and mixed coefficients satisfy
\[
 F_{-2x_nP}=F_{-2x_n^2P}=F_{C_T}=0,
 \qquad \mathcal C(-2x_nP,C_T)=\mathcal C(-2x_n^2P,C_T)=0.
\]

For the symmetric cubic block take all six permutations of
$T_{123}$ equal to one, put $P_{ij}=T_{ij,k}{x'}_k$ and
$Q_3=T_{ij,k}{x'}_i{x'}_j{x'}_k$, and use $(-2{x_n}P,-2{x_n}^2P)$.
For its first corrector a particular solution and residual datum
are
\[
 \phi=\frac{n-2}{6}(2{x_n}-1)Q_3D^{-n/2},\qquad
 q_*=\frac{n-2}{3}Q_3D^{-n/2}.
\]
The angular degree is three, so
$Q(h_{q_*})\le\tfrac14\int Dq_*^2$.
Define
\[
 Q_*:=I_n^{-1}\left(\int_{\R^n_+}F_{-2x_nP}\phi+\int_\Gamma q_*\phi
                                  +\tfrac14\int_\Gamma Dq_*^2\right).
\]
Equation~\eqref{eq:spectral-completion}, with the degree-one source
zero, shows that subtracting $Q_*/12$ from the first diagonal of the
direct matrix results in a lower matrix for the full form.

For the vector block, with $|v|=1$, put $L=v\cdot {x'}$,
$P_v={x'}\otimes v+v\otimes {x'}-2LI/(n-1)$, and use
$(-2{x_n}P_v,-2{x_n}^2P_v)$. The first particular solution is
\[
 \phi=L[P(x_n)D^{-n/2}+R(x_n)D^{-\frac{n-2}{2}}],
\]
where the two coefficient polynomials are
\[
 P(x_n)=-\frac{(n-2)^2}{n-1}(\tfrac23x_n^3+x_n^2),\qquad
 R(x_n)=\frac{(n-2)(n-3)}{2(n-1)}x_n-\frac{n-2}{2}.
\]
Its Robin residual is determined by
\[
 q_*=L\left(-(n-2)D^{-n/2}
        +\frac{(n-2)(n-3)}{2(n-1)}D^{-\frac{n-2}{2}}\right).
\]

Its degree-one spherical projection vanishes since
\[
 -(n-2)\int L^2D^{-n}
       +\frac{(n-2)(n-3)}{2(n-1)}\int L^2D^{-(n-1)}=0.
\]
The remaining inverse bound is $Q(h_{q_*})\le\tfrac12\int Dq_*^2$.
Here set
\[
 Q_*:=I_n^{-1}\left(\int_{\R^n_+}F_{-2x_nP_v}\phi+\int_\Gamma q_*\phi
                                  +\tfrac12\int_\Gamma Dq_*^2\right).
\]
Subtracting $Q_*/12$ from the first diagonal again yields a lower
matrix for the full form.

Apply the same determinant-ratio test to these smaller lower matrices.
The integral formulas and the stated inverse bounds imply, uniformly for
the three two-dimensional blocks and the hook block,
\[
 p_1>\frac1{10},\qquad p_2>\frac1{200}.
\]
For the hook block the remaining pivot satisfies $p_3>1$.
Every mixed block is therefore positive.
Orthogonal invariance and polarization extend the representative
calculations to their full isotypic components.

\subsection{The finite-jet form for Euclidean domains}

We evaluate the radial Gram matrices of the shape form
\eqref{eq:shape-form}. The ball Hessian computation
\eqref{eq:shape-hessian} and its kernel are derived in
Section~\ref{subsec:geometric-upper}.

For $f_j=H_\ell(z)|z|^{2j}$, with $\Delta H_\ell=0$ and
$\deg H_\ell=\ell$, write
$\widehat f_{\ell j}=(1+t)f_j(z)$ on $\Sph^{n-1}$ and put
$\kappa=n+2\ell-2$.
In stereographic coordinates $|z|^2=(1-t)/(1+t)$ and
$H_\ell(z)=((1-t)/(1+t))^{\ell/2}Y_\ell$, so
\[
 \widehat f_{\ell j}=(1-t^2)^{\ell/2}Y_\ell
                    (1-t)^j(1+t)^{1-\ell-j}.
\]
Take $Y_\ell$ from an $L^2(\Sph^{n-2})$-orthonormal
basis of degree-$\ell$ harmonics.

In \eqref{eq:gegenbauer-beta}, take $\lambda=\kappa/2$.
Then $I_1(k)=(-1)^kI_1(0)$, with
$I_1(0)=\kappa h_\ell/(\kappa-1)$ and
$h_\ell=\int_{-1}^1(1-t^2)^{\ell+(n-3)/2}\dd t$.
The same projection formula implies, for integers
$1\le p<(\kappa+1)/2$,
\[
 \frac{I_p(k)}{I_1(k)}=
 \frac{(k+1)_{p-1}(k+\kappa+1-p)_{p-1}}
 {2^{p-1}(p-1)!((\kappa+1)/2-p)_{p-1}}.
\]
Expand the radial factor as
\[
 (1-t)^j(1+t)^{1-\ell-j}
 =\sum_{b=0}^j(-1)^b\binom jb2^{j-b}
                                  (1+t)^{-(\ell+j-1-b)}.
\]
For $k+\ell\ge2$, the projection of the radial factor, divided by
$I_1(k)$, is
\[
 P_{\ell j}(k)=\sum_{\substack{0\le b\le j\\\ell+j-1-b\ge1}}
 (-1)^b\binom jb2^{j-b}\frac{I_{\ell+j-1-b}(k)}{I_1(k)}.
\]
The omitted powers contribute only degree-zero and degree-one
spherical harmonics, whose shape energy vanishes.

The total spherical degree of mode $k$ is $k+\ell$.
Divide the product of the projections by the Gegenbauer norm
in \eqref{eq:gegenbauer-beta} and multiply by the shape eigenvalue
$(k+\ell)(k+\ell-1)(k+\ell+n-1)$.
To sum the resulting series, define $M_\kappa$ linearly on
polynomials by its values on falling factorials,
\[
 M_\kappa(k^{\underline r})
          =\frac{(r!)^2(\kappa-1)}{\prod_{b=0}^r(\kappa-r-2+b)}.
\]
This is \eqref{eq:factorial-moment} with $b=1$, $c=\kappa$.
The convergence condition is $\kappa>r+2$ for every
falling-factorial degree $r$ that occurs. The radial Gram entries
per unit angular norm are
\begin{equation}\label{eq:euclidean-gram-evaluator}
 G^{(\ell)}_{ij}
 =\frac{2\kappa h_\ell}{(\kappa-1)^2}
 M_\kappa\bigl((k+\kappa/2)(k+\ell)(k+\ell-1)
                (k+\ell+n-1)P_{\ell i}(k)P_{\ell j}(k)\bigr).
\end{equation}
Here $i,j$ are the radial powers and $k$ is the Gegenbauer mode index.
The polynomial inside $M_\kappa$ has degree at most
$4\ell+2i+2j-4$. The retained graph degrees imply
\[
 \kappa-(4\ell+2i+2j-4)-2=n-2\ell-2i-2j>0,
\]
which verifies the convergence condition.
For each row of \eqref{eq:euclidean-radial-indices}, take the
corresponding principal submatrix and reindex its members in
increasing radial degree.
Substitution gives \eqref{eq:shape-correlation-bounds} and the
remaining correlation and Schur-pivot bounds in
Proposition~\ref{prop:euclidean-shape-form}.
The common positive factor $2\kappa h_\ell/(\kappa-1)^2$ cancels
in the correlations and in the normalization by $G^{(2)}_{00}$.

\par\medskip
\noindent\textbf{Concurrent work.}
In a preprint posted before ours, Bogo~\cite[Theorem~D]{BogoAsymptotics} established Fowler asymptotics with a positive power rate for smooth background metrics, assuming a priori two-sided critical growth of the boundary trace. In contrast, our Theorem~\ref{thm:conditional} assumes only the critical upper bound and derives the matching lower bound for nonremovable singularities, together with the Fowler asymptotics. We obtained the asymptotic result independently of Bogo's work, using the Jacobi-field classification established by Jin, Liao, and Xiong~\cite{JLXBoundary}.

\par\medskip
\noindent\textbf{Acknowledgements.}
The author thanks Professor Jingang Xiong for his continued advice
and support in choosing the research topic and writing this paper,
and for his guidance in the author's studies more generally.

\par\medskip
\noindent\textbf{AI usage statement.}
The author selected the overall approach and used GPT-6.1 Sol to carry
out extensive calculations. The model also assisted with language editing.
The author reviewed, verified, and substantially revised all adopted
material and takes full responsibility for the correctness and final
content of this paper.

\end{document}